\documentclass[12pt]{report}

\usepackage{amssymb,amsmath,amsthm}
\usepackage{xy}
\usepackage{enumerate}

\newtheorem{lemma}{Lemma}[section]
\newtheorem{definition}[lemma]{Definition}

\newtheorem{theorem}[lemma]{Theorem}
\newtheorem{remark}[lemma]{Remark}
\newtheorem{notation}[lemma]{Notation}

\newtheorem{corollary}[lemma]{Corollary}
\newtheorem{example}[lemma]{Example}

\begin{document}
\pagenumbering{roman}
\begin{center}
\thispagestyle{empty}
\mbox{}\\[2cm]
{\large FUSION ALGEBRAS,}\\[0.5cm]
{\large SYMMETRIC POLYNOMIALS,} \\[0.5cm]
{\large ORBITS OF N-GROUPS, } \\[0.5cm]
{\large AND RANK-LEVEL DUALITY}\\[3cm]
{\large BY}\\[0.5cm]
{\large Omar Saldarriaga}\\[1cm]
B.S., Antioquia University, 1997\\
M.A., State University of New York at Binghamton, 2000\\[6cm]
{\large DISSERTATION}\\[0.6cm]
Submitted in partial fulfillment of the requirements for\\
the degree of Doctor of Philosophy in Mathematics\\
in the Graduate School of\\
Binghamton University\\
State University of New York\\
2004
\end{center}
\newpage
\setcounter{page}{2}
\begin{center}
\hspace{1cm}\\[7cm]
\copyright \ Copyright by Omar Saldarriaga, 2004.\\
All rights reserved.
\newpage
\hspace{1cm}\\[6cm]
Accepted in partial fulfillment of the requirements for\\
the degree of Doctor of Philosophy in Mathematics\\
in the Graduate School of\\
Binghamton University\\
State University of New York\\
2004\\[6cm]
\end{center}
Alex Feingold \hrulefill \, June 1, 2004\\
Department of Mathematics\\[0.4cm]
Fernando Guzman \hrulefill \, June 1, 2004\\
Department of Mathematics\\[0.4cm]
Benjamin Brewster \hrulefill \, June 1, 2004\\
Department of Mathematics\\[0.4cm]
Charles Nelson \hrulefill \, June 1, 2004\\
Department of Physics\\[0.4cm]
\newpage
\renewcommand{\baselinestretch}{1.5}
\small\normalsize
\addcontentsline{toc}{chapter}{Abstract}
{\center\bf Abstract\\[0.5cm]}
A method of computing fusion coefficients for Lie algebras of type 
$A_{n-1}$ on level $k$ was recently developed by A. Feingold and M. Weiner~\cite{FW} using 
orbits of $\mathbb{Z}_n ^k$ under the permutation action of $S_k$ on $k$-tuples. 
They got the fusion coefficients only for n = 2 and 3. 
We will extend this method to all $n \geq 2$ and all $k \geq 1$. First we 
show a connection between Young diagrams and $S_k$-orbits of $\mathbb{Z}_n ^k$, 
and using Pieri rules we prove that this method works 
for certain specific weights that generate the fusion algebra. 
Then we show that the orbit method does not work in general, 
but with the help of the Jacobi-Trudi determinant, we give an iterative method to 
reproduce all type A fusion products.

\newpage
\addcontentsline{toc}{chapter}{Acknowledgements}
{\center\bf Acknowledgments\\[0.5cm]}
First of all I would like to thank my advisor Alex Feingold for all his help during my years at Binghamton, for all the time 
he devoted to me. I will never forget the many afternoons we expent talking about math. It has been the greatest experience 
in my life and made me grow mathematically a lot. I owe him my success in this program. There are no words to describe how 
grateful I am toward him.

I want to thank my friends for all the moments. They have been an integral part of my life during the last years. I wish the 
best to all of them. I want to thank my wife Diana for all her support and patience and my family in Colombia for all their 
support from a distance. It has been very difficult to be away from them.

I want to thank Professor B. Brewster, Professor C. Nelson and Professor F. Guzman for taking the time to be in my committee,
adding extra work to their already busy schedule. I am especially grateful to Professor Guzman whose proof reading and 
comments have made this work more accesible to non-specialists. I also want to thank him for all his support during my years 
at Binghamton and for giving me the opportunity to be in this Ph.D program. I enrolled in this program thanks to him.

I finally want to thank Professor L.-C. Kappe. Her hard work in advising and helping students to get jobs is incredible. I 
have never seen anyone so devoted to others. Her hard work always pays off because everyone in Binghamton gets a job thanks 
to her. She is, if not the most, one of the most valuable people in the mathematics department.
\newpage
\tableofcontents

\newpage
\pagenumbering{arabic}
\chapter{Introduction}

Fusion algebras play a very important role in conformal field theory \cite{Fu} and in the theory of vertex operator
algebras \cite{FZ}. Many equivalent interpretations have been found in other mathematical contexts such as quantum 
groups, Hecke algebras at roots of unity \cite{GW}, quantum cohomology of the Grassmanian \cite{BCF} and quite a few other 
areas.

There is a finite dimensional fusion algebra associated with an affine Kac-Moody algebra and a level $1\leq k\in \mathbb{Z}$,
with distinguished basis indexed by the level k weights. The structure constants $N_{a,b}^{(k)c}$ for such an algebra, also 
known as the fusion coefficients, can be expressed by the Kac-Walton algorithm (\cite{Kac} and \cite{Wa}) as an alternating 
sum of tensor product coefficients, or as a modification of the Racah-Speiser algorithm for tensor products.

This algorithm does not give any insight into why the structure constants for the fusion algebra are non-negative. This fact 
comes from other theories such as vertex operator algebras  where the fusion coefficients are the dimension of spaces of 
intertwining operators.

Many authors have tried to find a closed formula to express the fusion coefficients, but despite all the work of 
mathematicians and physicists during the 1990's the results in this direction are few, and closed formulas are known only 
for small ranks. Gepner and Witten \cite{GeWi} proved that for $A_1$
\[ N_{a,b}^{(k)c}=\begin{cases}
1, &\text{ if } c=a+b\mod 2 \text{  and  } |a-b|\leq c\leq\min\{a+b,2k-a-b\}\\
0, &\text{otherwise.}
\end{cases}
\]
For $A_2$ and $A_3$ closed formulas were given by Begin, Kirillov, Mathieu and Walton \cite{BKMW} and \cite{BMW}. They used 
combinatorial objects called Berenstein-Zelevinsky triangles. Another attempt at finding a closed formula was done by 
Feingold and Weiner \cite{FW}. They used $S_k$-orbits of $\mathbb{Z}_2^k$ and $\mathbb{Z}_3^k$. For $a\in\mathbb{Z}_N^k$, let
$[a]$ denote the $S_k$-orbit of $a$. There is a bijection between the level $k$ weights $\lambda$ of a type $A_{N-1}$ affine 
Kac-Moody algebra and the $S_k$-orbits $[\lambda]$ of $\mathbb{Z}_N^k$. For $A_1$ and $k\geq 1$ they proved 
\[N_{\mu ,\lambda}^{(k)\nu}=M_{[\mu] ,[\lambda]}^{(k)[\nu]}, \]
and for $A_2$ and $k\geq 1$ they proved
\[ M_{[\mu] ,[\lambda]}^{(k)[\nu]}=\binom{N_{\mu ,\lambda}^{(k)\nu}+1}{2},\]
where $M_{[\mu] ,[\lambda]}^{(k)[\nu]}$ is the number of orbits of 
$T([\mu],[\lambda],[\nu])=\{(x,y,z)\in [\mu]\times [\lambda]\times [\nu] \mid x+y=z \}$ under the action of $S_k$ on 
$T([\mu],[\lambda],[\nu])$ given by $S_k$ simultaneously acting on each $k$-tuple.

This work is an extension of the method of Feingold and Weiner to all ranks ($A_N$ for $N\geq 3$) and all levels. 
Although this method does not give a general closed formula for high ranks ($A_N$ for $N\geq 3$), it brings some progress to
the theory of type A fusion algebras. Using $S_k$-orbits of $\mathbb{Z}_N^k$, the Feingold-Weiner algorithm gives a closed 
formula for fusion coefficients $N_{\lambda,\mu}^{(k)\nu}=M_{[\lambda],[\mu]}^{(k)[\nu]}$ where $\mu$ is a weight 
of the form \eqref{E:weight}. Weights of this form include a set of generators for the fusion algebra, so we can compute the 
rest of the structure constants iteratively by using the Jacobi-Trudi determinant.

The structure of this thesis is as follows. Chapter 2 gives the definition of an abstract fusion algebra due to Fuchs 
\cite{Fu}, as well as preliminary material that leads to this definition. We also present a brief description  of the 
Racah-Speiser and Kac-Walton algorithms for tensor and fusion products and we finish the chapter by setting up the notation 
we will follow for the rest of the thesis. In Chapter 3 we describe the method of Feingold and Weiner \cite{FW} and give an 
interesting characterization for orbits with only zeroes and ones. In Chapter 4 we give a full description of the algebra of 
symmetric polynomials and its relation to the type A fusion algebras. We also describe the work of Goodman and Wenzl 
\cite{GW} that shows that the fusion algebra of type A is isomorphic to a quotient of the algebra of symmetric polynomials . 
We also include a proof of the type A rank-level duality that says that $\mathcal{F}'(A_{N-1},k)\cong\mathcal{F}'(A_{k-1},N)$
where $\mathcal{F}'(A_{N-1},k)$ is the quotient of the fusion algebra $\mathcal{F}(A_{N-1},k)$ by the ideal that identifies 
the set of simple currents with the identity. We finish the chapter with some generalities about simple currents. Chapter 5 
exhibits the connection between arithmetic of $S_k$-orbits of $\mathbb{Z}_N^k$ and fusion Pieri rules for symmetric 
polynomials and shows how the orbits method yields an algorithm for computing fusion coefficients for special weights. In
Chapter 6 we show an application of the orbits method to tensor products by showing that classical Pieri rules are 
equivalent to the orbit method for orbits with only zeroes and ones, a result that is an extension of the one for fusion
products. We finish the thesis by exhibiting a different version of each the Racah-Speiser and Kac-Walton algorithms from the
Young tableaux point view which are more practical than the classical algorithms for making computations by hand since, even
for $N\geq 3$, computations can be done on a piece of paper.

\chapter{Preliminaries } \label{S:pre}

In this chapter we give the definition of Lie algebras, the construction and classification of Kac-Moody Lie algebras, the 
definition of twisted and untwisted affine Lie algebras and its irreducible modules. We also set up the notation that we 
will follow for the rest of this work.

\section{Kac-Moody Lie algebras }

We begin by giving some important definitions that we will be using throughout the chapter.

\begin{definition}
a) A \textbf{Lie algebra} is a vector space L over a field $K$ (usually $K=\mathbb{R}$ or $\mathbb{C}$) equipped with a 
bilinear function $[\cdot,\cdot]:L\times L\rightarrow L$, called bracket, satisfying the conditions\

{\ 1)  $[x,x]=0, \quad$ for all $x\in L$. }

{\ 2) $[x,[y,z]]+[y,[z,x]]+[z,[x,y]]=0, \quad$ for all $x,y,z\in L$. }

\noindent
b) A subspace $I$ of a Lie algebra $L$ is called an \textbf{ideal} if for all $x\in L$ and for all $y\in I$ we have that
$[x,y]\in I$.

\noindent
c) We say that a Lie algebra $L$ is \textbf{simple} if the only ideals of $L$ are ${0}$ and $L$ and $\dim(L)>1$.

\noindent
d) The center $\mathcal{C}$ of a Lie algebra $L$ is the ideal 
\[ \mathcal{C}=\{ y\in L\mid [x,y]=0 \text{ for all } x\in L\}. \]

\noindent
e) We say that a vector space $V$ is module for a Lie algebra $L$, if there is a a bilinear function 
$\cdot:L\times V\rightarrow V$ satisfying
\[ [x,y]\cdot v=x\cdot(y\cdot v)-y\cdot(x\cdot v), \]
for all $x,y\in L$ and for all $v\in V$. If $V$ is a module for a Lie algebra $L$, we say that $V$ is an $L$-module.

\end{definition}

\begin{example}
Consider the vector space $sl_N$ of $N\times N$ matrices of trace 0 with the bracket defined by $[A,B]=AB-BA$. Under this 
bracket $sl_N$ becomes a Lie algebra, and it can be proved that this algebra is a simple Lie algebra. The vector space 
$V=K^N$ is a module for $sl_N$ under the action $x\cdot v=xv$ of matrix multiplication for $x\in sl_N$ and $v\in K^N$. We 
will see a construction of this Lie algebra in a later example in this chapter. 
\end{example}

We need the following definition before we present the definition of a Kac-Moody algebra. 

\begin{definition}
Given a matrix $A=(a_{ij})_{i,j=1}^t$ of rank $l$, we define a \textbf{realization} of $A$ to be a triple 
$\left(\mathfrak{h},\Pi,\check\Pi\right)$ where $\mathfrak{h}$ is a complex vector space of dimension $2t-l$, 
$\Pi=\{\alpha_1,...,\alpha_t\}\subseteq \mathfrak{h}^* $ and 
$\check\Pi=\{\check\alpha_1,...,\check\alpha_t\}\subseteq \mathfrak{h} $ satisfying 

\noindent
1) $\Pi$ and $\check\Pi$ are linearly independent,

\noindent
2) $\alpha_j(\check\alpha_i)=a_{ij}$. 
\end{definition}

We define a Lie algebra $\tilde{\mathfrak{g}}(A)$ with generators $e_i,f_i\quad (i=1,...,t)$ and $\mathfrak{h}$ 
satisfying the following relations
\begin{align} \label{E:relations}
& [h,h']=0, \qquad \text{for all } h,h'\in\mathfrak{h}; \notag\\
& [e_i,f_j]=\delta_{ij}\check{\alpha}_i, \qquad \text{for } i,j=1,...,t; \\
& [h,e_i]= \alpha_i(h)e_i,\qquad [h,f_i]=-\alpha_i(h)f_i, \quad \text{for all } h\in\mathfrak{h}\text{ and } i=1,...,t.
\notag 
\end{align}
The elements $e_i,f_i\quad (i=1,...,t)$  and the set $\mathfrak{h}$ are called the \textbf{Chevalley generators} of 
$\tilde{\mathfrak{g}}(A)$.

It can be proved that the Lie algebra $\tilde{\mathfrak{g}}(A)$ has a unique maximal ideal $\tau$ intersecting 
$\mathfrak{h}$ trivially. We can now define the Kac-Moody algebra associated with $A$.

\begin{definition}
Given a matrix A, we define the associated \textbf{Kac-Moody algebra} as the Lie algebra 
$\mathfrak{g}(A)=\tilde{\mathfrak{g}}(A)/\tau$ where $\tau$ is the unique maximal ideal intersecting $\mathfrak{h}$ 
trivially.
\end{definition}

\begin{example} \label{E:sln}
Consider the $(N-1)\times(N-1)$ matrix 
\[ A =  
\left( \begin{matrix}
    2 & -1 & 0 & \dots & 0 & 0 \\
    -1 & 2 & -1 &\dots & 0 & 0 \\
   0 & -1 & 2 & \dots & 0 & 0 \\
    \vdots & \vdots & \vdots & \ddots &  \vdots & \vdots \\
    0 & 0 & 0 & \dots & 2 & -1\\
0 & 0 & 0 & \dots & -1 & 2\\ 
 \end{matrix} \right)
\]
Let $E_{i,j}$ be the $N\times N$ matrix with 1 the $(i,j)$-entry and 0 everywhere else and set
\[ \check\alpha_i=E_{i,i}-E_{i+1,i+1}\quad i=1,...,N-1. \]

Let $\mathfrak{h}=Span\{\check\alpha_1,...,\check\alpha_{N-1}\}$. 

Define the linear functional $\epsilon_i$ by
\[ \epsilon_i(diag(a_1,...,a_N))=a_i, \quad i=1,...,N ,\]
and set 
\[ \alpha_i=\epsilon_i-\epsilon_{i+1}, \quad i=1,...,N-1, \]
then $\mathfrak{h}^*=Span\{\alpha_1,...,\alpha_{N-1}\}$. Now, set $e_i=E_{i,i+1}$ and $f_i=E_{i+1,i}$ for $i=1,...,N-1$. It 
can be verified that $e_i$, $f_i$ and $\mathfrak{h}$ satisfy the relations \eqref{E:relations}. It can also be proved that 
the algebra $\tilde{\mathfrak{g}}(A)$ is simple, so in particular it has no non-zero ideals which intersect $\mathfrak{h}$ 
trivially. Therefore $\mathfrak{g}(A)=\tilde{\mathfrak{g}}(A)$.

This Lie algebra is isomorphic to the Lie algebra of traceless $N\times N$ matrices defined in the previous example.
\end{example}Next
we  give the definition of generalized Cartan matrix.

The classification of Kac-Moody algebras is known for a special type of matrices known as generalized Cartan matrices, which
we now define. 

\begin{definition}
An integral $N\times N$ matrix $A=(a_{ij})_{i,j=1}^N$ is called a \textbf{generalized Cartan matrix} if satisfies the 
following conditions:

1) $a_{ii}=2$ for $1\leq i\leq N$,

2) $a_{ij}\leq 0$ for $i\ne j$,

3) $a_{ij}=0$ implies $a_{ji}=0$.
\end{definition}

Before stating the classification of Kac-Moody algebras associated to generalized Cartan matrices, we need  the following 
definition.

\begin{definition}
A matrix is called \textbf{decomposable} if, after some reordering of the indexes, the matrix has the form
\[ \left( \begin{matrix} A_1 & 0 \\ 0 & A_2 \end{matrix} \right). \]

If a matrix is not decomposable will be called \textbf{indecomposable}. 
\end{definition}

\begin{notation}
Let $u=\left[\begin{smallmatrix} u_1 \\ \vdots \\ u_N \end{smallmatrix} \right]\in \mathbb{R}^N$. We write $u>0$, if $u_i>0$
for $i=1,...,N$ and $u<0$ if $u_i<0$ for $i=1,...,N$.
\end{notation}

There are finite and infinite dimensional Kac-Moody Lie algebras. Every finite dimensional simple Lie algebra 
over $\mathbb{C}$ is a Kac-Moody algebra, and their classification, developed at the begining of the previous century,
is well known. The infinite dimensional Kac-Moody algebras are classified into two types, affine and indefinite. 
The affine are divided into twisted and untwisted and are well understood. 
There is less known about the indefinite type. Given a matrix A, one can tell if the Kac-Moody Lie algebra
$\mathfrak{g}=\mathfrak{g}(A)$ is finite, affine or indefinite. The classification is due to Vinberg \cite{Vin}. 
(See also \cite{Wan} and \cite{Kac}.)

\begin{theorem} \cite{Vin}
Let $A$ be an indecomposable generalized Cartan matrix, and $\mathfrak{g}=\mathfrak{g}(A)$ be the Kac-Moody algebra built 
from $A$. Then $\mathfrak{g}$ is one of the following three types.

a) $\mathfrak{g}$ is finite if and only if $\det A\ne 0$; $Au>0$ for some $u>0$; $Au\geq 0$ implies $u=0$ or $u>0$.

b) $\mathfrak{g}$ is affine if and only if corank$(A)=1$; $Au=0$ for some $u>0$; $Au\geq 0$ implies $Au=0$.

c) $\mathfrak{g}$ is indefinite if and only if $Au<0$ for some $u>0$; $Au\geq 0$ and $u\geq 0$ implies $u=0$.
\end{theorem}

\begin{definition}
Let $A$ be an indecomposable generalized Cartan matrix. We say that $A$ is of finite, affine or indefinite type, if the 
Kac-Moody algebra $\mathfrak{g}(A)$ built from $A$ is of finite, affine or indefinite type, respectively.
\end{definition}


\begin{example} \label{E:slnhat}
Consider the $N\times N$ matrix
\begin{equation} \label{E:cartan}
 \hat{A} =  [a_{ij}]_{i,j=0}^{N-1} = 
\left( \begin{matrix}
    2 & -1 & 0 & \dots & 0 & -1 \\
    -1 & 2 & -1 &\dots & 0 & 0 \\
   0 & -1 & 2 & \dots & 0 & 0 \\
    \vdots & \vdots & \vdots & \ddots &  \vdots & \vdots \\
    0 & 0 & 0 & \dots & 2 & -1\\
-1 & 0 & 0 & \dots & -1 & 2\\ 
 \end{matrix} \right).
\end{equation}
The Kac-Moody algebra $\mathfrak{g}(\hat A)$ is affine since 
$\hat{A}\cdot\left[\begin{smallmatrix} 1 \\ \vdots \\ 1 \end{smallmatrix} \right]=0$. Note also that by removing the first 
row and column we get the matrix $A$ from Example \ref{E:sln}, which is a matrix of finite type.

This matrix can be realized on a complex vector space $\mathfrak{H}$ of dimension $2N-(N-1)=N+1$ as follows.

Consider the matrix  
\[  A^E = [a_{ij}]_{i,j=0}^{N} = \left( \begin{matrix} 
\hat{A} & c \\
b & 2
\end{matrix} \right) \]
where $c^t=(0,\dots,0,1)=b$. For $i=0,...,N-1$, let $\check\alpha_{i}$ be the $i$th row of $A^E$ and let $\alpha_{i}$ be 
the $i$th coordinate function, i.e., the linear functional defined by
\[  \alpha_{i}(a_0,\dots,a_{N})=a_i. \]
Set
\[ \Pi=\{ \alpha_0,...,\alpha_{N-1}\} \quad \text{ and } \quad \check\Pi=\{\check\alpha_0,...,\check\alpha_{N-1} \}. \]
Then it is easy to check that $\{ \mathfrak{H}=\mathbb{C}^{N+1},\Pi,\check\Pi\}$ is a realization of $\hat A$. 

The Lie algebra $\mathfrak{g}(\hat A)$ with Chevalley generators $e_i,f_i$ for $i=1,...,N-1$ and $\mathfrak{H}$ defined by
\eqref{E:relations} is denoted by $\widehat{sl}_N$.
\end{example}

\section{Finite vs untwisted affine }

Let $A=(a_{ij})_{i,j=1}^{N-1}$ be an indecomposable generalized Cartan matrix of finite type. Then 
$\mathfrak{g}=\mathfrak{g}(A)$ is a finite dimensional simple Lie algebra of rank $N-1$. The finite dimensional simple Lie 
algebras over $\mathbb{C}$ are classified into four infinite families, $A_n$, $n\geq 1$, $B_n$, $n\geq 2$, $C_n$, $n\geq 3$ 
and $D_n$, $n\geq 4$, and five exceptional types, $E_6$, $E_7$, $E_8$, $F_4$ and $G_2$.

Let $(\mathfrak{h},\Pi,\check\Pi)$ be a realization of the matrix $A$ where 
$\Pi=\{\alpha_1,\dots,\alpha_{N-1}\}\subseteq\mathfrak{h}^*$ and 
$\check\Pi=\{\check\alpha_1,...,\check\alpha_{N-1}\}\subseteq\mathfrak{h}$. 
The subspace $\mathfrak{h}$ is called the Cartan sub-algebra of $\mathfrak{g}$. The sets $\Pi$ and $\check\Pi$ are called 
the set of simple roots and simple co-roots of $\mathfrak{g}$. We can define a symmetric bilinear form $(\cdot,\cdot)$ on 
$\mathfrak{g}$ which is invariant in the sense that $([x,y],z)=(x,[y,z])$, for all $x,y,z\in \mathfrak{g}$. The restriction 
of the form $(\cdot,\cdot)$ to $\mathfrak{h}$ is nondegenerate, so induces a form on $\mathfrak{h}^*$ which is determined by
\[ a_{ij}=\langle \alpha_i,\alpha_j \rangle=\dfrac{2(\alpha_i,\alpha_j)}{(\alpha_j,\alpha_j)},\quad 1\leq i,j\leq N-1  \]
and the normalization $(\theta,\theta)=2$ where $\theta$ is a special vector called ``the highest root" (we will give this 
vector explicitly and explain why it is called ``the highest root").

The fundamental weights of $\mathfrak{g}$ are linear functionals $\lambda_1,...,\lambda_{N-1}\in \mathfrak{h}^*$ satisfying 
\[   \langle \lambda_i,\alpha_j \rangle= \delta_{ij}, \quad 1\leq i,j\leq N-1. \]
The weight lattice of $\mathfrak{g}$ is the set 
\[ P=\left\{ \sum_{i=1}^{N-1}a_i\lambda_i \in\mathfrak{h}^*\bigg| a_i\in\mathbb{Z} \text{ for } 1\leq i\leq N-1 \right\} \]
 and the set of dominant integral weights of $g$ is
\[ P^+=\left\{ \sum_{i=1}^{N-1} a_i\lambda_i \ \bigg|\ 0\leq a_i \in \mathbb{Z} \right\}. \]
It is known that every irreducible finite dimensional module $V$ for $\mathfrak{g}$ is determined up to isomorphism by its 
highest weight $\lambda=\sum_{i=1}^{N-1}a_i\lambda_i$ where $0\leq a_i\in\mathbb{Z}$, for $i=1,...,N-1$. We now proceed to 
describe the structure of a $\mathfrak{g}$-module, define its highest weight and define the highest root of $\mathfrak{g}$. 
Any irreducible finite dimensional $\mathfrak{g}$-module $V$ has a direct sum decomposition $V=\bigoplus_{\mu\in P}V_\mu$ 
into weight spaces $V_\mu=\{v\in V\mid h\cdot v=\mu(h)v, \mbox{ for all } h\in \mathfrak{h} \}$. Here we use $x\cdot v$ for 
the action of $x\in \mathfrak{g}$ on $v\in V$. Each weight space $V_\mu$ is a generalized eigenspace for the simultaneous 
action of the abelian Cartan sub-algebra $\mathfrak{h}$, and we set $\Pi(V)=\{\mu\in P\mid V_\mu\ne 0\}$.

The simple Lie algebra $\mathfrak{g}$ becomes a $\mathfrak{g}$-module under the adjoint action
\[ ad_x(y)=x\cdot y=[x,y] \quad\text{for}\quad x,y\in\mathfrak{g}. \]
The weight space decomposition of $\mathfrak{g}$, called the Cartan decomposition, is 
\[ \mathfrak{g}=\mathfrak{h}\oplus\bigoplus_{0\ne\alpha\in P}\mathfrak{g}_\alpha, \]
where the root spaces for $0\ne\alpha\in P$ are
\[ \mathfrak{g}_\alpha=\{x\in\mathfrak{g}\mid [h,x]=\alpha(h)x, \mbox{ for all } h\in \mathfrak{h} \} \]
and we define the set of roots by
\[ \Phi=\{0\ne\alpha\in P\mid \mathfrak{g}_\alpha\ne 0\}.\]
It is well known that $\dim(\mathfrak{g}_\alpha )=1$ for $\alpha\in\Phi$ and since, by definition, $[h,e_i]=\alpha_i(h)e_i$ 
and $[h,f_i]=-\alpha_i(h)f_i$, for $1\leq i\leq N-1$, we have $\alpha_i,-\alpha_i\in\Phi$. One has, in fact, that for any 
$\alpha\in\Phi$, $\alpha=\sum_{i=1}^{N-1}n_i\alpha_i$ where either all $0\leq n_i\in\mathbb{Z}$ or $0\geq n_i\in\mathbb{Z}$.
That is, $\Phi=\Phi^+\cup \Phi^-$ decomposes into positive and negative roots (with all simple roots positive). This gives a 
triangular  decomposition $\mathfrak{g}=\mathfrak{g}^-\oplus\mathfrak{h}\oplus\mathfrak{g}^+$ where 
$\mathfrak{g}^\pm=\oplus_{\alpha\in\Phi^\pm}\mathfrak{g}_\alpha$.

Define a partial order, $\leq$, on $P$ by 
\[ \lambda_1\leq\lambda_2 \quad\text{if and only if}\quad \lambda_2-\lambda_1=\sum_{i=1}^{N-1}n_i\alpha_i,\]
where $0\leq n_i\in \mathbb{Z}$ for $1\leq i\leq N-1$.

Let $V$ be a finite dimensional irreducible $\mathfrak{g}$-module. There is a unique weight $\lambda\in \Pi(V)$ such that 
$\mathfrak{g}^+\cdot V_\lambda=0$, $\lambda\geq\beta$ for all $\beta\in \Pi(V)$,  and 
$\lambda=\sum_{i=1}^{N-1}a_i\lambda_i\in P^+$ ($0\leq a_i\in\mathbb{Z}$) determines $V$ up to  isomorphism. We call 
$\lambda$ the highest weight of $V$ and write $V=V^\lambda$, and 
$\Pi^\lambda=\Pi(V^\lambda)=\{\beta\in\mathfrak{h}^*\mid V_\beta^\lambda\ne 0\}\subset P$. The highest weight $\theta$ of 
$\mathfrak{g}$ as a $\mathfrak{g}$-module is always a root and it is called the highest root of $\mathfrak{g}$.
 
The Weyl group $W$ of $\mathfrak{g}$ is the subgroup of $GL(\mathfrak{h}^*)$ generated by the simple reflections 
\[ r_i(\lambda)=\lambda-\lambda(\check\alpha_i)\alpha_i, \quad \text{for}\quad i=1,...,N-1. \]
Every element $w\in W$ can be written as a word in the $r_i$'s and the \textbf{length $l(w)$ of $w$} is defined as the number
of simple reflections in a reduced word for $w$.

A matrix of affine type can be obtained from a matrix of finite type by adding one column and a row so that the resulting 
matrix has a 1-dimensional kernel. The affine Lie algebra corresponding to a matrix obtained in this way has the following 
construction. 

Let $\hat A=\left(\begin{matrix}  2 & b \\ c & A\end{matrix} \right)$ be an $N\times N$ 
indecomposable generalized Cartan matrix of affine type, where $A$ is the matrix of finite type considered at the begining of
this section. There is a unique vector  $u=\left[\begin{smallmatrix} u_0 \\ \vdots \\ u_{N-1} \end{smallmatrix} \right]$ such
that $u>0$, $\hat{A}u=0$, $u_i\in \mathbb{Z}$ for $i=0,\dots,N-1$ and $\gcd{(u_0,\dots,u_{N-1})}=1$. The vectors 
$b=(b_1,\dots,b_{N-1})$ and $c=(c_1,\dots,c_{N-1})^t$ are determined by 
\[ c_i=-( u_1a_{i1}+\dots+ u_{N-1}a_{i,N-1}),\quad b_i=-(\hat{u}_1a_{1i}+\dots+\hat{u}_{N-1}a_{N-1,i}),
\text{ for } 1\leq i\leq N-1, \]
where $\hat{u}_1,\dots,\hat{u}_{N-1}$ are the last $N-1$ components of the unique vector 
$\hat{u}=\left[\begin{smallmatrix}\hat u_0 \\ \vdots \\ \hat u_{N-1} \end{smallmatrix} \right]$ satisfying 
\[ \hat{A}^t \hat{u}=0,\quad \hat u_0,\dots,\hat u_{N-1}\in \mathbb{Z}\quad\text{and}\quad 
\gcd{(\hat u_0,\dots,\hat u_{N-1})}=1, \]
and $a_{i1},a_{i2}\dots,a_{i,N-1}$ for $i=1,\dots N-1$ are the components of the $i$th row of the matrix A and $a_{1i},
a_{2i},\dots,a_{N-1,i}$ for $i=1,\dots N-1$ are the components of the $i$th column of A.

Let $\hat{\mathfrak{g}}=\mathfrak{g}(\hat A)$ be the Kac-Moody algebra built from $\hat A$. Then  $\hat{\mathfrak{g}}$ is an 
affine Lie algebra with realization $(\mathfrak{H},\Psi,\hat{\Psi})$ where $\dim\mathfrak{H}=N+1$. The number 
\[ h=\sum_{i=0}^{N-1}u_i\] 
 is called the Coxeter number and 
\[ \check{h}=\sum_{i=0}^{N-1}\check{u}_i \]
 is called the \textbf{dual Coxeter number}.
The Lie algebra $\hat{\mathfrak{g}}$ has a one dimensional center spanned by 
\[ c=\sum_{i=0}^{N-1}\check{u}_i\check\alpha_i. \]
The element $c$ is called the \textbf{canonical central element}.
The highest root of the finite dimensional Lie algebra $\mathfrak{g}$ is explicitly given by
\[ \theta=\sum_{i=1}^{N-1}u_i\alpha_i. \]

The realization of $\hat{\mathfrak{g}}$, $(\mathfrak{H},\Psi,\check\Psi)$, can be constructed so that 
$\mathfrak{h}\subseteq\mathfrak{H}$, $\Pi\subseteq\Psi$
and $\check\Pi\subseteq\check\Psi$. This can be done as follows. 

There exists an element in $d\in\mathfrak{H}$ satisfying the conditions $\alpha_0(d)=1$ and $\alpha_i(d)=0$ for 
$i=1,\dots,N-1$. The vector $d$ is uniquely determined up to multiple of $c$ and 
$B=\{c,\check\alpha_1,\dots,\check\alpha_{N-1},d\}$ is a basis for $\mathfrak{H}$.

Let $d^*\in\mathfrak{H}^*$ be the linear functional dual to $d$ with respect to the basis $B$. This linear functional is 
called the \textbf{null root} of $\hat{\mathfrak{g}}$ and it is given explicitly by 
\[ d^*=\sum_{i=0}^{N-1}u_i\alpha_i\in\mathfrak{H}^*. \]

Let $c^*\in\mathfrak{H}^*$ be the element dual to the canonical central element $c$ with respect to the basis $B$. Then we 
have $c^*(c)=1$, $c^*(d)=0$ and $d^*(c)=0$.

We take $\Psi=\{\alpha_0,\alpha_1,\dots,\alpha_{N-1}\}$ and 
$\check\Psi=\{\check\alpha_0,\check\alpha_1,\dots,\check\alpha_{N-1}\}$. The sets $\Psi$ and $\check\Psi$ are called the set 
of simple roots and simple co-roots of $\hat{\mathfrak{g}}$. 

We also identify linear functionals in $\mathfrak{h}^*$ with linear functionals in $\mathfrak{H}^*$ by having the same values
on $\mathfrak{h}$ and being zero on $\alpha_0$ and $d$.

The fundamental weights of $\hat{\mathfrak{g}}$ are linear functionals $\Lambda_0,\dots,\Lambda_{N-1}\in \mathfrak{h}^*$ 
satisfying 
\[   \langle \Lambda_i,\alpha_j \rangle= \delta_{ij}, \quad 0\leq i,j\leq N-1, \]
and it can be checked that the fundamental weights are given by 
\[  \Lambda_0=c^*,\quad \Lambda_i=u_i\dfrac{(\alpha_i,\alpha_i)}{2}c^*+\lambda_i, \quad \text{for}\quad 1\leq i\leq N-1, \]
where $\Lambda_0=c^*+rd^*$ is normalized by taking $r=0$ and $\lambda_1,\dots,\lambda_{N-1}$ are the weights for the 
finite dimensional Lie algebra $\mathfrak{g}$.

The set of integral weights for $\hat{\mathfrak{g}}$ is defined by 
\[ \hat P=\left\{\Lambda=\sum_{i=0}^{N-1}n_i\Lambda_i \bigg|  n_i\in\mathbb{Z} \right\}. \]

The set of dominant integral weights for $\hat{\mathfrak{g}}$ is defined by 
\[ \hat P^+=\left\{\Lambda=\sum_{i=0}^{N-1}n_i\Lambda_i \bigg| 0\leq n_i\in\mathbb{Z} \right\}. \]

We can also define a partial order $\leq$ on $\mathfrak{H}^*$ by
\[ \Lambda\leq \Lambda' \text{ if and only if } \Lambda'-\Lambda=\sum_{i=0}^{N-1}a_i\alpha_i, \quad\text{for } 0\leq a_i\in
\mathbb{Z}. \]

The affine Lie algebras are classified into 2 types, twisted and untwisted. For the untwisted affine Lie algebras 
$\hat{\mathfrak{g}}$ we get an isomorphism of Lie algebras
\[ \hat{\mathfrak{g}}\cong \mathbb{C}[t,t^{-1}]\otimes \mathfrak{g}\oplus\mathbb{C}c\oplus\mathbb{C}d, \]
where $\mathfrak{g}$ is the finite dimensional Lie algebra and the bracket on 
$\mathbb{C}[t,t^{-1}]\otimes\mathfrak{g}\oplus\mathbb{C}c\oplus\mathbb{C}d$ is defined by 
\[ [x(m),y(n)]=[x,y](m+n)+m\delta_{m+n,0}(x,y)c,\]
\[ [d,x(m)]=mx(m),\quad \text{and}\quad [c,x(m)]=[c,d]=0, \]
where $x(m)=t^m\otimes x$, and $(\cdot,\cdot)$ is the normalized invariant symmetric bilinear form from $\mathfrak{g}$, 
$x,y\in\mathfrak{g}$. The restriction of $(\cdot,\cdot)$ to $\mathfrak{h}$ is non-degenerate, so it induces a form on 
 $\mathfrak{h}^*$. The normalization is taken so that $(\theta,\theta)=2$. We say that $\hat{\mathfrak{g}}$ is of type 
$X_n^{(1)}$ if $\mathfrak{g}$ is of type $X_n$, where $X=A,B,C,D,E,F$ or $G$.

We can extend the symmetric bilinear form $(\cdot,\cdot)$ on $\mathfrak{g}$ to $\hat{\mathfrak{g}}$ by 
\[ \left( x(m),y(n)\right)=\delta_{m+n,0}(x,y),\quad \left(x(m),c\right)=\left(y(n),d\right)=(c,c)=(d,d)=0,\]
and
\[ (c,d)=1. \]

The Weyl group $\widehat{W}$ of $\mathfrak{\hat{g}}$ is the subgroup of $GL(\mathfrak{H}^*)$ generated by the simple 
reflections 
\[ r_i(\Lambda)=\Lambda-\Lambda(\check\alpha_i)\alpha_i, \quad \text{for}\quad i=0,...,N-1, \]
and $W=\langle r_i\mid 1\leq i\leq N-1\rangle$ is a subgroup of $\widehat{W}$. For each $0\leq k\in \mathbb{Z}$ the affine 
Weyl group $\widehat{W}$ acts on the weight lattice $P$ of $\mathfrak{g}$ with the usual action of the simple reflections of 
$W$ and with 
\begin{equation} \label{E:affineaction}
 r_0(\beta)=r_{\theta}(\beta)+(k+\check{h})\theta,
\end{equation}
where $r_{\theta}(\lambda)=\lambda-\dfrac{2(\lambda,\theta)}{(\theta,\theta)}\theta$ is the reflection with respect to the 
highest root $\theta$ and $\check{h}$ is the dual Coxeter number.
The set of weights of level $k$ of $\mathfrak{g}$ is the set
\[ P_k^+=\left\{ \sum_{j=1}^{N-1} a_j\lambda_j\in P^+ \ \bigg|\ \sum_{j=1}^{N-1}a_j\leq k \right\} ,\]
and the fundamental region for the action of $\widehat{W}$ on $P$ is $P_{k+\check{h}}^+$.

The representation theory of $\hat{\mathfrak{g}}$ is very similar to the representation theory of ${\mathfrak{g}}$ if we 
restrict only to $\hat{\mathfrak{g}}$-modules $V$ that satisfy the following, where 
 \[ V_{\Lambda}=\{v\in V\mid h\cdot v= \Lambda(h)v, \quad \text{for all}\quad h\in\mathfrak{H}\},\qquad
\Pi(V)=\{\Lambda\in \mathfrak{H}^*\mid V_{\Lambda}\ne 0\}\] and \[D(\lambda)=\{\mu\in\mathfrak{H}^*\mid\mu\leq\lambda\}. \]

i) $V$ is $\mathfrak{H}$-diagonalizable, i.e., $V=\bigoplus_{\Lambda\in\mathfrak{H}^*}V_{\Lambda}$.

ii) $\dim V_{\Lambda}<\infty$, for all $\Lambda\in\mathfrak{H}^*$.

iii) There exits a finite number of elements $\beta_1,\dots,\beta_s$ so that $\Pi(V)\subseteq\cup_{i=1}^s D(\beta_i)$.

As a consequence, every $\hat{\mathfrak{g}}$-module $V$ satisfying the above conditions has a highest weight 
$\Lambda\in\hat P$ and if $V$ is irreducible then it would be determined, up to isomorphism, by its highest weight $\Lambda$.
We use the notation $V^{\Lambda}$ to denote the irreducible highest weight module for $\hat{\mathfrak{g}}$. The canonical 
central element acts on $V^{\Lambda}$ as a scalar $k$ and we have
\[ k=\Lambda(c)=\sum_{i=0}^{N-1}\Lambda_i(c)=\sum_{i=0}^{N-1}u_i\dfrac{(\alpha_i,\alpha_i)}{2} \]
so for a fixed $k$, there are only finitely many $\Lambda\in\hat P^+$ with $\Lambda(c)=k$.

\begin{example}
The untwisted affine algebra $\widehat{sl}_N$ defined in Example \ref{E:slnhat} has simple roots
\[ \alpha_0,\dots,\alpha_{N-1}. \]
The vector $u=\left[\begin{smallmatrix} 1\\ \vdots \\1 \end{smallmatrix}\right]$ satisfies $\hat A\cdot u=0$ where $\hat A$ 
is the matrix given in \eqref{E:cartan} and since $\hat A$ is symmetric we also have $\hat A^t\cdot u=0$. Then the canonical 
central element $c$ and the null root $d^*$ are given by 
\[ c=\sum_{i=0}^{N-1}\check\alpha_i \quad\text{and}\quad  d^*=\sum_{i=0}^{N-1}\alpha_i.\]
The fundamental weights are given by the equations
\begin{equation} \label{E:affineweights}
 \Lambda_0=c^*,\quad \Lambda_i=c^*+\lambda_i, \quad \text{for}\quad 1\leq i\leq N-1, 
\end{equation}
where $\lambda_1,\dots,\lambda_{N-1}$ are the fundamental weights for the Lie algebra $\mathfrak{g}=sl_N$. 

Every irreducible highest weight module $\hat V$ for $\widehat{sl}_N$ on level $k$ is determined, up to isomorphism, by a 
weight $\Lambda=\sum_{i=0}^{N-1}n_i\Lambda_i$ satisfying the condition $k=\sum_{i=0}^{N-1}n_i$. From 
\eqref{E:affineweights} we see that 
\[ \Lambda=n_0c^*+\sum_{i=1}^{N-1}(n_ic^*+n_i\lambda_i)=kc^*+\sum_{i=1}^{N-1}n_i\lambda_i, \]
where $\sum_{i=1}^{N-1}n_i\leq k$.

The Coxeter and dual Coxeter numbers of $\hat{\mathfrak{g}}$ are $h=\check{h}=N$ and the level $k$ action of the affine 
reflection $r_0$ on the weight lattice $P$ of $\mathfrak{g}$ is given by
\[ r_0(\lambda)=r_{\theta}(\lambda)+(k+\hat{h})\theta, \]
where $\theta=\sum_{i=1}^{N-1}\alpha_i$ is the highest root of $\mathfrak{g}=sl_N$.
\end{example}

\section{Fusion algebras }

Let us begin with the definition of fusion algebra, due to J. Fuchs \cite{Fu}. (Also see \cite{Fe}.)

\begin{definition} \label{d:fusionalgebra}
A fusion algebra is a finite dimensional commutative associative algebra $\mathcal{F}$ over $\mathbb{Q}$ satisfying the 
following. 

1) There is a distinguished basis 
\[ B=\{ x_a\mid a\in A \},\quad\text{for some  finite index set }A,  \]
so that the product of the basis elements is given by 
\[ x_a \cdot x_b=\sum_{c\in A}N_{a,b}^c x_c, \]
where $0\leq N_{a,b}^c \in\mathbb{Z}$.

2) There is an element $\omega\in A$ such that the map
\begin{equation} \label{E:conjugation}
 \mathcal{C}(x_a)=\sum_{b\in A}N_{a,b}^{\omega} x_b 
\end{equation}
is an involution. That is, the matrix $C=(C_{ab})_{a,b\in A}$ for $\mathcal{C}$, where $C_{ab}=N_{a,b}^{\omega}$, satisfies 
the equation $C^2=I$ (the identity matrix). The map $\mathcal{C}$ is called the \textbf{conjugation map}.
\end{definition}

We have several consequences from the definition. Since the entries of $C$ are non-negative integers and $C^2=I$, then it 
follows that either $C=I$ or $C$ is a permutation matrix of order 2. Hence, there is a permutation $\sigma:A\rightarrow A$
satisfying $\sigma^2=1$ and so that
\[ C_{a,b}=\delta_{a,\sigma(b)}. \]
The conjugation map $\mathcal{C}:\mathcal{F}\rightarrow\mathcal{F}$ defined in \eqref{E:conjugation}, is then given by
\[ \mathcal{C}(x_a)=\sum_{b\in A}\delta_{a,\sigma(b)}x_b=x_{\sigma(a)} \]
and since it is an automorphism, we also get
\begin{align}
 x_{\sigma(a)}\cdot x_{\sigma(b)}&=\mathcal{C}(x_a)\mathcal{C}(x_b)=\mathcal{C}(x_a\cdot x_b)=
\mathcal{C}\left(\sum_{c\in A}N_{a,b}^c x_c\right)\notag\\&=\sum_{c\in A}N_{a,b}^c \mathcal{C}(x_c)=
\sum_{c\in A}N_{a,b}^c x_{\sigma(c)}. \notag
\end{align}
Therefore we get
\[ N_{a,b}^c=N_{\sigma(a),\sigma(b)}^{\sigma(c)}. \]
Set
\[ N_{a,b,c}=N_{a,b}^{\sigma(c)}. \]
It follows from commutativity and associativity of $\mathcal{F}$ that the constants $N_{a,b,c}$ are totally symmetric in
$a,b$ and $c$. Then we also get 
\[ N_{\omega,b}^c=N_{\omega,b,\sigma(c)}=N_{b,\sigma(c),\omega}=N_{b,\sigma(c)}^{\sigma(\omega)}=
N_{\sigma(b),c}^{\omega}=C_{\sigma(b),c}=\delta_{b,c}, \]
hence
\[ x_{\omega}x_b=\sum_{c\in A}N_{\omega,b}^c x_{c}=\sum_{c\in A}\delta_{b,c}x_{c}=x_b. \]
So $x_{\omega}$ is a multiplicative identity and $\sigma(\omega)=\omega$. 

\begin{example}
Consider a 3-dimensional algebra $\mathcal{F}$ over $\mathbb{Q}$ with basis \newline$B=\{x_0,x_1,x_2\}$ and product given by 
the table (where the blank entries  are determined by commutativity.)
\begin{center}
\begin{tabular}{ |c || c|c|c|}\hline 
 & $x_0$ & $x_1$ & $x_2$ \\ \hline \hline
 $x_0$ & $x_0$ & $x_1$ & $x_2$ \\ \hline
$x_1$ &  & $x_0$ & $x_2$ \\ \hline
$x_2$ & &  & $x_0+x_1$ \\ \hline
\end{tabular}
\end{center}
It can be checked that the product defined by this table is associative with identity element $x_0$, and the conjugation map 
$\mathcal{C}$ is the identity map with $C=I$. Therefore $\mathcal{F}$ is a fusion algebra.
\end{example}

\section{Fusion algebras from affine Lie algebras } \label{S:racaspeicer}

In this section we discuss the decomposition of the tensor product of two irreducible finite dimensional modules for a 
finite dimensional Lie algebra $\mathfrak{g}$ and use it to define a ``level $k$ truncated tensor product" of two irreducible
highest weight modules for the untwisted affine Lie algebra $\hat{\mathfrak{g}}$. This truncated product defines a fusion 
algebra with basis given by irreducible highest weight modules of $\hat{\mathfrak{g}}$ on level $k$.

Let $\mathfrak{g}$ be a finite dimensional algebra and let $V^{\lambda}$ and $V^{\mu}$ be two irreducible finite dimensional 
modules for $\mathfrak{g}$. The tensor product $V^{\lambda}\otimes V^{\mu}$ is completely reducible. Therefore it can be 
expressed as a direct sum of irreducible modules
\[ V^{\lambda}\otimes V^{\mu}=\bigoplus_{\nu\in P^+}Mult_{\lambda,\mu}(\nu)V^{\nu}. \]
The coefficients $Mult_{\lambda,\mu}(\nu)$ can be computed by means of the Racah-Speiser algorithm which can be described as 
follows. For any $\nu\in P^+$ we have
\begin{equation} \label{E:racah}
 Mult_{\lambda,\mu}(\nu)=\sum_{w\in W}(-1)^{l(w)}Mult_{\lambda}(w(\nu+\rho)-\mu-\rho), 
\end{equation}
where $Mult_{\lambda}(\beta)=dim(V_{\beta}^\lambda)$ and $\rho=\sum_{i=1}^{N-1}\lambda_i$. This formula yields the following 
geometric algorithm.

\noindent
Step 1. 

Shift the weight diagram $\Pi^{\lambda}$ of $V^{\lambda}$ by adding $\mu +\rho$ to each weight.

\noindent
Step 2.

 Use the finite Weyl group to move all shifted weights $\Pi^{\lambda}+\mu +\rho$ into the dominant chamber, $P^+$, where they
accumulate as an alternating sum of inner multiplicities of $V^{\lambda}$, adding if the required $w$ is even, subtracting if
it is odd.

\noindent
Step 3.

 The resulting pattern of numbers will be non-negative integers, zero if on a chamber wall, and after shifting the pattern 
back by subtracting $\rho$ we get the tensor product multiplicities.

\begin{example}
Consider the modules $V^{\lambda}$ and $V^{\mu}$ of $sl_3$ where $\lambda=\lambda_1+\lambda_2$ and $\mu=2\lambda_1$. The 
weight decomposition of $V^{\lambda}$  is given by
\[ V^{\lambda}=V_{\lambda_1+\lambda_2}\oplus V_{2\lambda_1-\lambda_2}\oplus V_{-\lambda_1+2\lambda_2}\oplus 2V_{0}\oplus 
V_{\lambda_1-2\lambda_2}\oplus V_{-2\lambda_1+\lambda_2}\oplus V_{-\lambda_1-\lambda_2}. \]
So the weight diagram of $V^{\lambda}$ is given by 
\vskip10pt
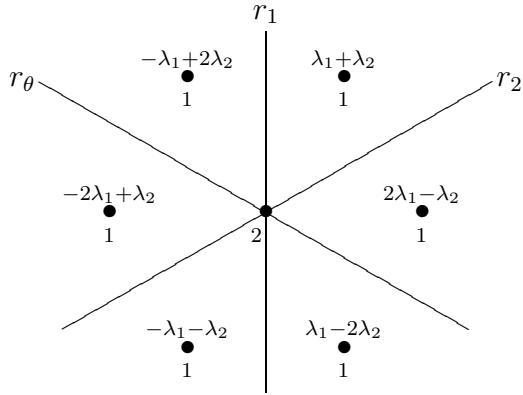
\begin{xy}
<1.2cm,0cm>:
(-6.62,0); (0,0) **{ };
(.866,1.5) *{\bullet};
(.866,1.7) *{_{\lambda_1+\lambda_2}};
(.866,1.25) *{_1};
(1.73,0) *{\bullet};
(1.73,.2) *{_{2\lambda_1-\lambda_2}};
(1.73,-.25) *{_1};
(-.866,1.5) *{\bullet};
(-.866,1.7) *{_{-\lambda_1+2\lambda_2}};
(-.866,1.25)*{_1};
(-1.73,0) *{\bullet};
(-1.73,.2) *{_{-2\lambda_1+\lambda_2}};
(-1.73,-.25) *{_1};
(.866,-1.5) *{\bullet};
(.866,-1.3) *{_{\lambda_1-2\lambda_2}};
(.866,-1.75) *{_1};
(-.866,-1.5) *{\bullet};
(-.866,-1.3) *{_{-\lambda_1-\lambda_2}};
(-.866,-1.75) *{_1};
(0,0) *{\bullet};
(-.1,-.25) *{_2};
(2.511,1.45); (-2.251,-1.3) **@{-};
(2.7,1.45) *{r_2};
(0,2); (0,-2) **@{-};
(0,2.2) *{r_1};
(-2.511,1.45); (2.251,-1.3) **@{-};
(-2.7,1.45) *{r_\theta};
\end{xy}
\noindent
where the number below each dot is the inner multiplicity, or the dimension of the weight space, and the two lines labeled 
$r_1$ and $r_2$ are the fixed lines of the simple reflections.
\vskip10pt
Now, adding $\mu +\rho=2\lambda_1+(\lambda_1+\lambda_2)=3\lambda_1+\lambda_2$ to all weights in the above diagram, we get
\vskip10pt
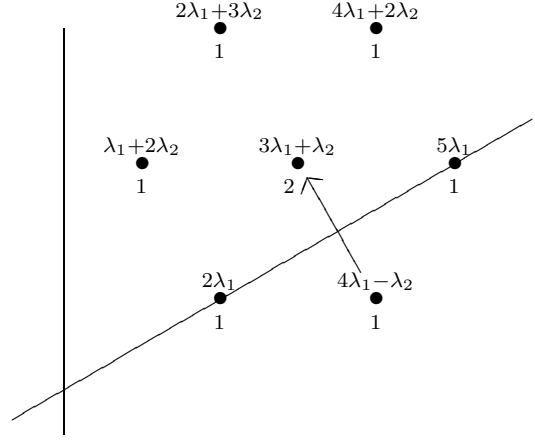
\begin{xy}
<1.2cm,0cm>:
(-6.62,0); (0,0) **{ };
(3.46,4) *{\bullet};
(3.46,3.75) *{_1};
(3.46,4.2)  *{_{4\lambda_1+2\lambda_2}};
(4.33,2.5) *{\bullet};
(4.33,2.7) *{_{5\lambda_1}};
(4.33,2.25) *{_1};
(1.73,4) *{\bullet};
(1.73,3.75) *{_1};
(1.73,4.2) *{_{2\lambda_1+3\lambda_2}};
(2.59,2.5) *{\bullet};
(2.5,2.25) *{_2};
(2.59,2.7) *{_{3\lambda_1+\lambda_2}};
(3.46,1) *{\bullet};
(3.46,1.2) *{_{4\lambda_1-\lambda_2}};
(3.46,.75) *{_1};
(.866,2.5) *{\bullet};
(.866,2.7) *{_{\lambda_1+2\lambda_2}};
(.866,2.25) *{_1};
(1.73,1) *{\bullet};
(1.73,1.2) *{_{2\lambda_1}};
(1.73,.75) *{_1};
(5.19,3); (-.57,-.333) **@{-};
(0,4); (0,-.5) **@{-};
(3.29,1.3); (2.7,2.35) **@{-};
(2.87,2.33);  (2.7,2.35)  **@{-};
(2.63,2.17); (2.7,2.35)  **@{-};
\end{xy}
\vskip10pt
\noindent
where the arrow indicates that the weight outside the fundamental chamber gets reflected onto the weight inside the 
fundamental chamber, reducing its multiplicity by 1, and since the weights on the chamber wall do not 
count for the tensor product, we are left with the pattern
\vskip10pt
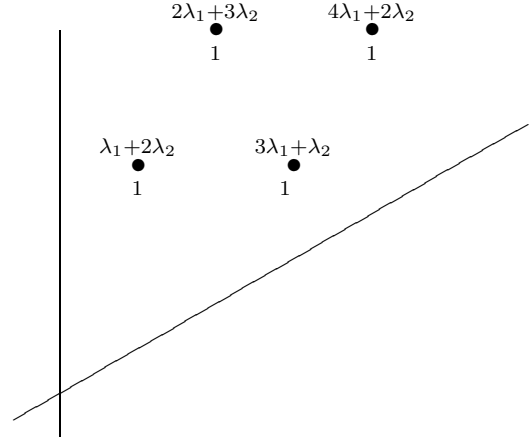
\begin{xy}
<1.2cm,0cm>:
(-6.62,0); (0,0) **{ };
(3.46,4) *{\bullet};
(3.46,3.75) *{_1};
(3.46,4.2)  *{_{4\lambda_1+2\lambda_2}};
(1.73,4) *{\bullet};
(1.73,3.75) *{_1};
(1.73,4.2) *{_{2\lambda_1+3\lambda_2}};
(2.59,2.5) *{\bullet};
(2.5,2.25) *{_1};
(2.59,2.7) *{_{3\lambda_1+\lambda_2}};
(.866,2.5) *{\bullet};
(.866,2.25) *{_1};
(.866,2.7) *{_{\lambda_1+2\lambda_2}};
(5.19,2.96); (-.51,-.325) **@{-};
(0,4); (0,-.5) **@{-};
\end{xy}
\vskip10pt
Subtracting $\rho=\lambda_1+\lambda_2$ from all these weights we get the tensor product decomposition of 
$V^{\lambda_1+\lambda_2}\otimes V^{2\lambda_1}$, that is
\[ V^{\lambda_1+\lambda_2}\otimes V^{2\lambda_1}=V^{3\lambda_1+\lambda_2}\oplus V^{\lambda_1+2\lambda_2}\oplus
V^{2\lambda_1}\oplus V^{\lambda_2}. \]
\end{example}

The Racah-Speiser algorithm was modified by Kac and Walton for computing the fusion coefficients 
$N_{\mu,\lambda}^{(k)\nu}$ for $\lambda,\mu,\nu\in P^+_k$, as follows. Under the level $k$ action of $\widehat{W}$ on $P$ we
have 
\[ N_{\mu,\lambda}^{(k)\nu}=\sum_{w\in \widehat{W}}(-1)^{l(w)}Mult_{\lambda}(w(\nu+\rho)-\mu-\rho), \]
where $Mult_{\lambda}(\beta)=dim(V_{\beta}^\lambda)$ and $\rho=\sum_{i=1}^{N-1}\lambda_i$. This formula now involving a 
summation over the affine Weyl group, $\widehat{W}$, also yields a 
geometric algorithm as follows.

\noindent
Step 1. 

Shift the weight diagram of $V^{\lambda}$ by adding $\mu +\rho$ to each weight.

\noindent
Step 2. 

Use the level $k$ action of the affine Weyl group $\widehat W$ on $P$ to move all shifted weights into the fundamental 
domain, $P_k^+$, bounded by the reflection walls of all simple reflections $r_i$, for $1\leq i\leq N-1$, and of the affine 
reflection $r_0$. That affine reflection is a hyperplane perpendicular to $\theta$ going through the point 
$(k+\check h)\theta$. The reflected weights counted with inner multiplicities accumulate as an alternating sum of inner 
multiplicities of $V^{\lambda}$, adding if the required $w$ is even, subtracting if it is odd.

\noindent 
Step 3. 

The resulting pattern of numbers will be non-negative integers, zero if on a reflection wall, and after shifting the pattern 
back by subtracting $\rho$, we get the level $k$ fusion product multiplicities.

\begin{example}
Consider the modules $V^{\lambda}$ and $V^{\mu}$ of $sl_3$ where $\lambda=\lambda_1+\lambda_2$ and $\mu=2\lambda_1$ and let 
$k=2$. The weight decomposition and weight diagram of $V^{\lambda}$ is given in the previous example.
Adding $\mu +\rho=2\lambda_1+\lambda_1+\lambda_2=3\lambda_1+\lambda_2$ to all weights in its diagram, we get
\vskip10pt
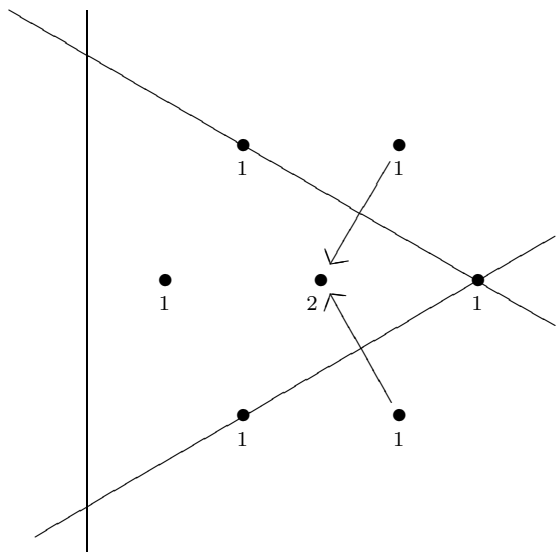
\begin{xy}
<1.2cm,0cm>:
(-6.62,0); (0,0) **{ };
(3.46,4) *{\bullet};
(3.46,3.75) *{_1};
(4.33,2.5) *{\bullet};
(4.33,2.25) *{_1};
(1.73,4) *{\bullet};
(1.73,3.75) *{_1};
(2.59,2.5) *{\bullet};
(2.5,2.25) *{_2};
(3.46,1) *{\bullet};
(3.46,.75) *{_1};
(.866,2.5) *{\bullet};
(.866,2.25) *{_1};
(1.73,1) *{\bullet};
(1.73,.75) *{_1};
(5.19,3); (-.57,-.333) **@{-};
(0,5.5); (0,-.5) **@{-};
(3.38,1.15); (2.7,2.35) **@{-};
(2.87,2.33);  (2.7,2.35)  **@{-};
(2.63,2.17); (2.7,2.35)  **@{-};
(-.86,5.5); (5.19,2) **@{-};
(3.36,3.82); (2.7,2.69) **@{-};
(2.7,2.69); (2.64,2.85) **@{-};
(2.7,2.69); (2.9,2.72) **@{-};
\end{xy}
\vskip10pt
\noindent
where the arrows indicate each of the weights that are outside the affine fundamental domain, and get reflected onto a 
weight inside the affine fundamental domain. In this case, 2 outside weights get reflected onto the weight 
$3\lambda_1+\lambda_2$, reducing its multiplicity by 2, and since the weights on the reflection wall do not contribute to the
fusion product, we are left with the pattern
\vskip10pt
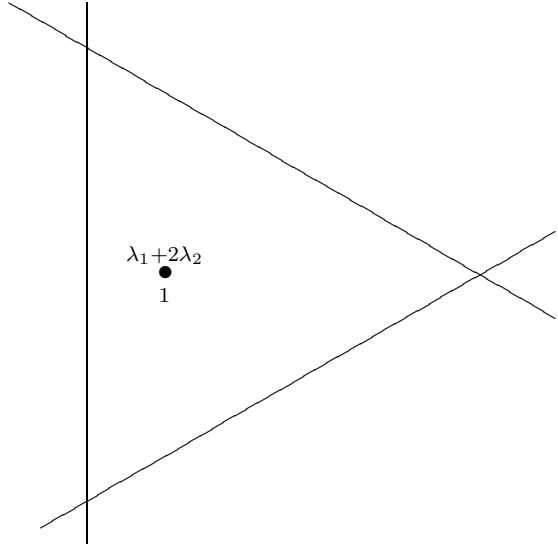
\begin{xy}
<1.2cm,0cm>:
(-6.62,0); (0,0) **{ };
(.866,2.5) *{\bullet};
(.866,2.25) *{_1};
(.866,2.7) *{_{\lambda_1+2\lambda_2}};
(5.19,2.96); (-.51,-.325) **@{-};
(0,5.5); (0,-.5) **@{-};
(-.86,5.5); (5.19,2) **@{-};
\end{xy}
\vskip10pt
Subtracting $\rho=\lambda_1+\lambda_2$ from this weight we get the fusion product  
\[ V^{\lambda_1+\lambda_2}\otimes_2 V^{2\lambda_1}=V^{\lambda_2}. \]

\end{example}

\section{Notation for the rest of the thesis } \label{S:intro}

Let $\hat{\mathfrak{g}}$ be the affine algebra of type $A_{N-1} ^{(1)}$ built from $\mathfrak{g}=sl_N$. Let 
$\lambda _1,\lambda _2,...,\lambda _{N-1}$ denote the fundamental weights of $\mathfrak{g}$ and let:
\begin{equation} \label{E:domintwts}
 P^+=\left\{ \sum_{j=1}^{N-1} a_j\lambda_j \ \bigg|\ 0\leq a_j \in \mathbb{Z} \right\} 
\end{equation}
denote the dominant integral weights of $\mathfrak{g}$. The irreducible modules $V^\Lambda$ for $\hat{\mathfrak{g}}$ of level 
$k$ are indexed by dominant highest weights $\Lambda=kc+\lambda$ where $\lambda$ is in 
\[ P_k^+=\left\{ \sum_{j=1}^{N-1} a_j\lambda_j\in P^+ \ \bigg|\ \sum_{j=1}^{N-1}a_j\leq k \right\} .\]
However, for a fixed level $k$, we will denote these modules by $V^\lambda$ instead, where $\lambda\in  P_k^+$. These weights
are in one-to-one correspondence with the set of $N$-tuples $(a_0,a_1,...,a_{N-1})$ whose sum is $k$, where 
$a_0=k-\sum_{j=1}^{N-1}a_j$.

Let G be the group $\mathbb{Z}_N ^k$ and let $S_k$ act on it by permuting the $k$-tuples, so every orbit of $\mathbb{Z}_N ^k$
under this action has a unique standard representative in the form
\[ \left( (N-1)^{a_{N-1}},...,1^{a_1},0^{a_0}\right) \]
where the exponent indicates the number of repetitions of the base. 

We get a one-to-one correspondence between $S_k$-orbits of $\mathbb{Z}_N ^k$ and $N$-tuples whose sum is $k$, by:
\begin{equation} \label{E:partit}
 \left( (N-1)^{a_{N-1}},...,1^{a_1},0^{a_0}\right)\longleftrightarrow (a_0,a_1,...,a_{N-1}) 
\end{equation}
therefore we get a correspondence between weights $\lambda\in P_k^+$ and orbits of $\mathbb{Z}_N ^k$ under the described 
action of $S_k$. The orbit corresponding to $\lambda$ will be denoted by $[\lambda]$ and the correspondence is given by
\begin{equation} \label{E:ow}
\lambda=\sum_{j=1}^{N-1}a_j\lambda_j \in P_k^+ \longmapsto [\lambda]=
\left[\left( (N-1)^{a_{N-1}},...,1^{a_1},0^{a_0}\right)\right],
\end{equation}
where $a_0=k-\sum_{j=1}^{N-1}a_j$ and $\left[\left( (N-1)^{a_{N-1}},...,1^{a_1},0^{a_0}\right)\right]$ denotes the orbit of
$\mathbb{Z}_N ^k$ whose standard representative is $\left( (N-1)^{a_{N-1}},...,1^{a_1},0^{a_0}\right)$.

A partition is a finite sequence of non-negative integers $(\mu_1,\dots,\mu_n,\dots)$ so that 
$\mu_1\geq\mu_2\geq\dots\geq\mu_n\geq\dots$. The length of the partition $\mu$, $l(\mu)$, is the number of non-zero 
$\mu_i$'s.

There is also a map between partitions of length at most $N$ and dominant integral weights of $A_{N-1}$  given by: 
\begin{equation} \label{E:wp}
 (\mu_1,...,\mu_N)\longmapsto \lambda=\sum_{j=1}^{N-1}(\mu_j-\mu_{j+1})\lambda_j \in P^+
\end{equation}
and if $\mu_1-\mu_N\leq k$ then $\lambda\in P_k^+$.

Note that this map is not one-to-one since given a partition $(\mu)=(\mu_1,...,\mu_N)$ with $\mu_N\ne 0$, then the partitions
$(\mu)$ and $(\mu_1-\mu_N,\mu_2-\mu_N,...,\mu_{N-1}-\mu_N)$ lead to the same weight. However, this correspondence has a right
inverse which can be described as follows. If $\lambda=\sum_{j=1}^{N-1}a_j\lambda_j\in P^+$ then the partition associated to 
the weight $\lambda$, which will be denoted by $(\lambda)$, is given by the map
\begin{equation} \label{E:pw}
 \lambda=\sum_{j=1}^{N-1}a_j\lambda_j \longmapsto (\lambda)=\left(\sum_{j=1}^{N-1} a_j,\sum_{j=2}^{N-1} a_j,...,a_{N-1}
\right),
\end{equation}
and if $\lambda\in P_k^+$ then $(\lambda)$ has largest part at most $k$.

Another way of defining ~\eqref{E:wp} is by defining an equivalence relation $\thicksim$ on partitions of length at most $N$ 
by 
\begin{equation} \label{E:er}
 (\mu_1,...,\mu_N)\thicksim(\nu_1,...,\nu_N) \text{ if and only if } \mu_i-\mu_{i+1}=\nu_i-\nu_{i+1} \text{ for } 
1\leq i\leq N-1.
\end{equation}

The equivalence class for $(\mu_1,...,\mu_N)$  has a unique representative in the set of partitions of length at most $N-1$ 
given by $(\mu_1-\mu_N,\mu_2-\mu_N,...,\mu_{N-1}-\mu_N)$ and it is clear from ~\eqref{E:wp} that both can be used to get the 
weight associated to $(\mu_1,...,\mu_N)$. Furthermore two such partitions $(\mu)$ and $(\mu')$ are in the same equivalence 
class under $\thicksim$ if and only if $\mu'_i=\mu_i+c$ for some integer $c$ and for $1\leq i\leq N$.

It is clear that the left hand side of \eqref{E:partit} is also a partition, as the sequence is non-increasing, which will be
proved to be the conjugate of the partition on the right hand side of \eqref{E:pw}. Combining \eqref{E:ow} and \eqref{E:pw}
gives a map from $S_k$-orbits of $\mathbb{Z}_N ^k$ to partitions of length at most $N-1$ with largest part at most $k$. And 
from \eqref{E:ow} and \eqref{E:wp} we get a map from partitions $(\mu_1,...,\mu_N)$ of length at most $N$ with 
$\mu_1-\mu_N\leq k$ to $S_k$-orbits of $\mathbb{Z}_N ^k$. This will be useful in proving that for specific $\lambda$ to be 
described in section 4 below, the operation among $S_k$-orbits defined by Feingold and Weiner \cite{FW}, matches level $k$ 
fusion Pieri rules for computing fusion coefficients $N_{\lambda ,\mu}^{(k)\nu}$. We will prove that these specific weights
include generators of the fusion ring of type $A_{N-1}$ of level $k$, $\mathcal{F}(A_{N-1},k)$. Therefore, using the 
Jacobi-Trudi determinant, a simple iteration allows one to compute the rest of the structure constants 
$N_{\lambda ,\mu}^{(k)\nu}$ for arbitrary $\mu$ and $\nu$.

\noindent
This describes a new method of computing type $A$ coefficients which uses only arithmetic in $\mathbb{Z}_N ^k$ and an
iteration to obtain the whole multiplication table for $\mathcal{F}(A_{N-1},k)$.

\chapter{Orbits of $\mathbb{Z}_N ^k$ under the action of $S_k$ } \label{C:oua}

In this chapter we follow Feingold and Weiner \cite{FW}. Let $G$ be the group $\mathbb{Z}_N ^k$, and for each $N$-tuple of 
nonnegative integers 
\[ (a_0,a_1,...,a_{N-1}) \qquad \text{ such that} \qquad a_0+a_1+...+a_{N-1}=k \]
we define the subset of $G$
\begin{equation} \label{E:Part}
[(a_0,a_1,...,a_{N-1})]=\left\{ x\in\mathbb{Z}_N ^k \mid j \text{ occurs }  a_j \text{ times in } x \text{, } 
0\leq j\leq N-1 \right\}.
\end{equation}
Then $G$ is a disjoint union of these subsets. 

Note that the symmetric group $S_k$ acts on $G$ by permuting $k$-tuples, the \textbf{set of orbits under this action},
$\mathcal{O}=\mathcal{O}(N,k)$,  consists of the subsets \eqref{E:Part} defined above, and each orbit contains a unique 
representative in \textbf{standard form} 
\begin{equation} \label{E:orb}
 \left( (N-1)^{a_{N-1}},...,1^{a_1},0^{a_0}\right) 
\end{equation}
where the exponent indicates the number of repetitions of the base. 

\begin{notation}
Given $x\in \mathbb{Z}_N ^k$, we will denote the \textbf{orbit of x} by $[x]$ and the \textbf{representative in standard 
form of this orbit} will be denoted by ${\hat x}$.
\end{notation}

For orbits $[a]$,$[b]$ and $[c]$, we define the set:
\[ T([a],[b],[c])=\{(x,y,z)\in [a]\times [b]\times [c] \mid x+y=z \}. \]
Note that $\sigma\in S_k$ acts on $(x,y,z)\in T([a],[b],[c])$ by $\sigma(x,y,z)=(\sigma x,\sigma y,\sigma z)$.
\begin{definition}
Denote by $\mathbf{M_{[a],[b]}^{(k)[c]}}$\textbf{ the number of }$\mathbf{S_k}$-\textbf{orbits of }$\mathbf{T([a],[b],[c])}$.
\end{definition}

This suggests that we could use these numbers as the structure constants of an algebra 
(depending on $N$ and $k$ and over any field of characteristic 0) with basis $\mathcal{O}$, by defining a bilinear product 
as follows:
\begin{equation}  \label{E:orbprod}
[a]\times [b]=\sum_{c\in \mathcal{O}}M_{[a],[b]}^{(k)[c]} [c].
\end{equation}
In \cite{FW}, it was shown that for $N=2$ and all $k\geq 1$, this product coincides with the product in the associative 
fusion algebra $\mathcal{F}(A_1,k)$, but for $N=3$ and all $k\geq 1$ it was shown that 
\begin{equation} \label{E:FW}
 M_{[\mu] ,[\lambda]}^{(k)[\nu]} =\binom{N_{\mu,\lambda}^{(k)\nu}+1}{2} . 
\end{equation}

Now we proceed to describe how to compute the product of two $S_k$-orbits of $\mathbb{Z}_N ^k$. This description will be 
useful in proving some theorems in later chapters.

\begin{definition}
Let $[a]$, $[b]\in\mathcal{O}$, and assume that $[b]=\{y_1,...,y_t \}$. For $1\leq i\leq t$ set:
\begin{equation} \label{E:list}
 z_i=\hat a +y_i .
\end{equation}
We say that the equation $z_j=\hat a +y_j$ in the list \eqref{E:list} is \textbf{redundant}, if for some $i<j$ and 
$\sigma\in S_k$ we have
\[ \sigma\hat{a} =\hat{a} \text{, }\sigma y_j=y_i\text{ and } \sigma z_j=z_i ,\]
that is, if the triples $(\hat{a} ,y_i,z_i)$ and $(\hat{a} ,y_j,z_j)$ are in the same $S_k$-orbit of $T([a],[b],[z_i])$.
\end{definition}

Now we can describe how the product of two orbits can be computed. Let $[a]$, $[b]\in\mathcal{O}$ and fix a representative 
from the orbit $[a]$, say the representative in standard form, $\hat a$, and assume that $[b]=\{y_1,...,y_t \}$. For every
$y_i\in [b]$, set 
\[ z_i=\hat a +y_i , \qquad 1\leq i\leq t . \]
Remove all redundant equations from this list, and without loss of generality, assume that after removing all redundancies, 
we are left with the first $s$ equations, for some $s\leq t$. That is, the list:
\[ z_i=\hat a +y_i , \qquad 1\leq i\leq s,  \]
has no redundancies. Then 
\[ [a]\times [b]=[z_1]+[z_2]+...+[z_s] . \]
Note that several $z_i$'s could be in the same orbit, and for every $[c]\in\mathcal{O}$,  
\begin{equation} \label{E:M1}
 M_{[a],[b]}^{(k)[c]}=Card\{1\leq i\leq s \mid z_i\in [c] \} . 
\end{equation}

\begin{remark} \label{R:M_[a],[b]^[c]}
From Equation \eqref{E:M1}, we also get that $M_{[a],[b]}^{(k)[c]}$ can be computed by removing all redundancies from the 
list of equations 
\[  z=\hat a +y, \]
where $y\in[b]$ and $z\in[c]$.
\end{remark}

\begin{example} \label{E:or1}
Let $a=(2,1,0)$ and $b=(1,1,0)$ be elements in $\mathbb{Z}_3 ^3$.  Then we have that $\hat a =(2,1,0)$ and 
$[b]=\{ (1,1,0), (1,0,1), (0,1,1) \}$. Now to compute $[a]\times [b]$, we remove all redundancies from the list:
\begin{align}
& (2,1,0) + (1,1,0) = (0,2,0) ,\notag\\
           &(2,1,0)+ (1,0,1) = (0,1,1) ,\notag\\
           &(2,1,0)+ (0,1,1) = (2,2,1) .\notag  
\end{align}
Since there are no redundant equations in that list, we get that:
\begin{equation} \label{E:orbpieri}
 [(2,1,0)]\times [(1,1,0)]=  [(2,2,1)]+[(1,1,0)]+[(2,0,0)] .
\end{equation}
\end{example}

\begin{example} \label{E:or}
Now let $a=(2,2,1)$ and $b=(1,0,0)$ be elements in $\mathbb{Z}_3 ^3$.  Then we have that $\hat a =(2,2,1)$ and 
$[b]=\{ (1,0,0), (0,1,0), (0,0,1)\}$. To compute $[a]\times [b]$, we remove all redundancies from the list:
\begin{align}
  &(2,2,1)+ (1,0,0) = (0,2,1) , \notag\\
           &(2,2,1)+ (0,1,0) = (2,0,1) ,\notag\\
           &(2,2,1)+ (0,0,1) = (2,2,2) .\notag  
\end{align}
We can see that the second equation is redundant, so we can remove it from the list and we get
\[ [(2,2,1)]\times [(1,0,0)]=  [(2,2,2)]+[(2,1,0)]. \]
\end{example}

\begin{example} \label{E:or2}
Now let $a=(3,2,1)$ and $b=(1,1,0)$ be elements in $\mathbb{Z}_4 ^3$.  Then we have that $\hat a =(3,2,1)$ and 
$[b]=\{ (1,1,0), (1,0,1), (0,1,1)\}$. To compute $[a]\times [b]$, we remove all redundancies from the list:
\begin{align} 
  &(3,2,1)+ (1,1,0) = (0,3,1) , \notag\\
           &(3,2,1)+ (1,0,1) = (0,2,2) ,\notag\\
           &(3,2,1)+ (0,1,1) = (3,3,2) .\notag  
\end{align}
There are no redundant equations in this list, so we get:
\begin{equation} \label{E:example2}
 [(3,2,1)]\times [(1,1,0)]=[(3,1,0)]  +[(2,2,0)]+[(3,3,2)]. 
\end{equation}
\end{example}

In our next example we show how this product of orbits fails to be associative.

\begin{example} \label{E:nonassos}
Let $a=(2,1,0)$, $b=(1,1,0)$ and $c=(1,0,0)$ in $\mathbb{Z}_3 ^3$. An easy calculation using the results of Examples 
\ref{E:or1} and \ref{E:or} shows that: 
\[ \left( [a]\times [b]\right)\times [c] = [(2,2,2)]+[(1,1,1)]+3[(2,1,0)]+[(0,0,0)] \]
and 
\[ [a]\times\left( [b]\times [c]\right)=[(2,2,2)]+[(1,1,1)]+4[(2,1,0)]+[(0,0,0)] .\]
Note that multiplicities of $[(2,1,0)]$ do not match.
\end{example}

Our aim is to prove, once we establish the connection between Young diagrams and orbits, that the operation 
\eqref{E:orbprod} matches fusion Pieri rules when $[b]$ is the orbit corresponding to a weight of a certain form. 

In Examples ~\ref{E:or1}, ~\ref{E:or} and ~\ref{E:or2} we computed the product of two orbits where one of the orbits 
consisted of only zeroes and ones and we always got multiplicity 1 for every one of the orbits in the product. The next 
lemma, which is essential in meeting our aim, shows that this is always true. (In general, multiplicities could be
higher than one. See Example~\ref{E:last} at the end of this chapter.)

\begin{theorem} \label{T:orbitproduct}
Let $[a],[c]\in\mathcal{O}$ and $[b]=[(1^m,0^{k-m})]$ for some $m\leq k$. Suppose that $M_{[a],[b]}^{(k)[c]}\ne 0$ then 
$M_{[a],[b]}^{(k)[c]}=1$.
\end{theorem}
\begin{proof}
Let's assume that the representative in standard form, $\hat{a}$, of the the orbit $[a]$ is
\begin{equation} \label{E:hat}
\hat{a}=\left( (N-1)^{a_{N-1}},...,1^{a_1},0^{a_0}\right).  
\end{equation}
The following claim is the key of the proof.

\noindent
Claim: If   $z\in[c]$ and $y\in[b]$ satisfy the equation
\begin{equation} \label{E:nonred}
z=\hat{a}+y,
\end{equation}
then there exists $\sigma\in S_k$ so that $\sigma\hat a =\hat a$, 
\[ \sigma y  =(1^{m_{N-1}},0^{a_{N-1}-m_{N-1}},...,1^{m_1},0^{a_1-m_1},1^{m_0},0^{a_0-m_0}) \] 
and
\[ \sigma z=\left( 0^{m_{N-1}},(N-1)^{a_{N-1}-m_{N-1}+m_{N-2}},...,1^{a_1-m_1+m_0},0^{a_0-m_0}\right),\]
where $m_0$, $m_1$,...,$m_{N-1}$ are integers such that $0\leq m_i\leq a_i$, for $0\leq i\leq N-1$ and \newline
$m_0+m_1+...+m_{N-1}=m$.

\noindent
Proof of the claim:

Consider the bijection:
\begin{align} 
\alpha :  \mathbb{Z}_N^k & \longrightarrow\mathbb{Z}_N^{a_{N-1}}\times...\times \mathbb{Z}_N^{a_1}\times\mathbb{Z}_N^{a_0} 
\notag\\
 p &\longmapsto (p_{N-1},...,p_1,p_0) \notag
\end{align}
where $p_{N-1}$ consists of the first $a_{N-1}$ entries of $p$, $p_{N-2}$ contains the next $a_{N-2}$ entries of $p$ and 
so on, until finally $p_0$ contains the last $a_0$ entries of $p$.
Let $(y_{N-1},...,y_1,y_0)\in\mathbb{Z}_N^{a_{N-1}}\times...\times \mathbb{Z}_N^{a_1}\times\mathbb{Z}_N^{a_0}$ be the 
image of $y$ under this map.
Note that for every $0\leq j\leq N-1$, $y_j$ is an $a_j$-tuple consisting of only 0's and 1's. Let $m_j$ be the number of 
1's in $y_j$ (hence $\sum m_j=m$), and let $\sigma_j$ be a permutation in the symmetric group $S_{a_j}$ so that 
\begin{equation} \label{E:y}
\sigma_j y_j=(1^{m_j},0^{a_j -m_j})
\end{equation}
Now, let each of the $\sigma_j$ act on $p\in\mathbb{Z}_N^k$ by 
\[ \bar{\sigma}_j p=\alpha^{-1}(p_{N-1},...,\sigma_j p_j,...,p_1,p_0), \]
where $(p_{N-1},...,p_1,p_0)=\alpha (p)$.
and set $\sigma=\bar{\sigma}_0\dots \bar{\sigma}_{N-1}$.
Let $(a_{N-1},...,a_1,a_0)=\alpha(\hat a)$ so equation \eqref{E:hat} gives $a_j=(j^{a_j})$ and $\sigma$ stabilizes 
$\hat a$, that is,
\begin{equation} \label{E:stab}
\sigma\hat a=\hat a.
\end{equation}
From \eqref{E:y} and the definition of $\sigma$, we get that:
\begin{align} \label{E:sig}
 \sigma y &  =\sigma (y'_{N-1},...,y'_1,y'_0) \\
& =(1^{m_{N-1}},0^{a_{N-1}-m_{N-1}},...,1^{m_1},0^{a_1-m_1},1^{m_0},0^{a_0-m_0}). \notag
\end{align}
We can write $\hat a$ as follows:
\begin{align} \label{E:mua}
 \hat a & =\left( (N-1)^{a_{N-1}},...,1^{a_1},0^{a_0}\right) \\
& =\left( (N-1)^{m_{N-1}},(N-1)^{a_{N-1}-m_{N-1}},...,1^{m_1},1^{a_1-m_1},0^{m_0},0^{a_0-m_0}\right). \notag
\end{align}
Now since the action of $S_k$ on $\mathbb{Z}_N^k$ is linear, we get from equations \eqref{E:nonred}, \eqref{E:stab},
\eqref{E:sig} and \eqref{E:mua} that:
\begin{align} \label{E:z2}
\sigma z  & =\hat a + \sigma y \\
& =\left( N^{m_{N-1}},(N-1)^{a_{N-1}-m_{N-1}},...,2^{m_1},1^{a_1-m_1},1^{m_0},0^{a_0-m_0}\right) \notag\\
& =\left( N^{m_{N-1}},(N-1)^{a_{N-1}-m_{N-1}+m_{N-2}},...,1^{a_1-m_1+m_0},0^{a_0-m_0}\right) \notag\\
& =\left( 0^{m_{N-1}},(N-1)^{a_{N-1}-m_{N-1}+m_{N-2}},...,1^{a_1-m_1+m_0},0^{a_0-m_0}\right) \notag
\end{align}
which finishes the proof of the claim.

Now, to complete the proof of the theorem we see that from Remark~\ref{R:M_[a],[b]^[c]} we have that 
\[ M_{[a],[b]}^{(k)[c]}= \mbox{ the number of non-redundant equations of the form \eqref{E:nonred} } , \]
where $y\in[b] $ and $z\in[c]$. 

Since $M_{[a],[b]}^{(k)[c]}\ne 0$, there exists $y\in[b]$ and $z\in[c]$ satisfying \eqref{E:nonred}. We have to prove that if
$y'\in[b]$ and $z'\in[c]$ also satisfy \eqref{E:nonred}, then there exists $\gamma\in S_k$ so that
\[ \gamma\hat{a}=\hat{a}, \qquad \gamma y=y' \qquad\mbox{ and }\qquad \gamma z=z'. \] 
From the claim we have that since $z=\hat{a}+y$, then there exits $\sigma\in S_k$ such that $\sigma\hat a=\hat a$,
\begin{equation} \label{E:sigma y}
\sigma y  =(1^{m_{N-1}},0^{a_{N-1}-m_{N-1}},...,1^{m_1},0^{a_1-m_1},1^{m_0},0^{a_0-m_0}) 
\end{equation}
and
\begin{equation} \label{E:sigmaz}
 \sigma z=\left( 0^{m_{N-1}},(N-1)^{a_{N-1}-m_{N-1}+m_{N-2}},...,1^{a_1-m_1+m_0},0^{a_0-m_0}\right),
\end{equation}
for some integers $m_0$, $m_1$,...,$m_{N-1}$ such that $0\leq m_i\leq a_i$, for $0\leq i\leq N-1$ and \newline
$m_0+m_1+...+m_{N-1}=m$.

From Equation \eqref{E:sigmaz} it follows that the representative in standard form, $\hat{c}$, from the orbit 
$[c]=[\sigma z]=[z]$ has the form
\begin{equation} \label{E:hat{c}}
 \hat{c} =\left((N-1)^{a_{N-1}-m_{N-1}+m_{N-2}},...,1^{a_1-m_1+m_0},0^{a_0-m_0+m_{N-1}}\right) .
\end{equation}
Now, since $z'=\hat{a}+y'$, the claim also guarantees the existence of a permutation $\sigma'\in S_k$ so that 
$\sigma'\hat a=\hat a$,
\begin{equation} \label{E:sigma y'}
 \sigma' y'  =(1^{t_{N-1}},0^{a_{N-1}-t_{N-1}},...,1^{t_1},0^{a_1-t_1},1^{t_0},0^{a_0-t_0}) 
\end{equation}
and 
\begin{equation} \label{E:cilin}
 \sigma' z'=\left( 0^{t_{N-1}},(N-1)^{a_{N-1}-t_{N-1}+t_{N-2}},...,1^{a_1-t_1+t_0},0^{a_0-t_0}\right) 
\end{equation}
for some integers $t_0$, $t_1$,...,$t_{N-1}$ such that $0\leq t_i\leq a_i$, for $0\leq i\leq N-1$ and \newline
$t_0+t_1+...+t_{N-1}=m$.
From Equation \eqref{E:cilin} it follows that
\begin{equation} \label{E:hat{c1}}
 \hat{c} =\left((N-1)^{a_{N-1}-t_{N-1}+t_{N-2}},...,1^{a_1-t_1+t_0},0^{a_0-t_0+t_{N-1}}\right). 
\end{equation}
Since the representative in standard form $\hat{c}$ of the orbit $[c]$ is unique, from Equations \eqref{E:hat{c}} and
\eqref{E:hat{c1}}, we get the system of equations 
\begin{align}
 a_{N-1}-t_{N-1}+t_{N-2} & =a_{N-1}-m_{N-1}+m_{N-2}\notag\\
 a_{N-2}-t_{N-2}+t_{N-3} & =a_{N-2}-m_{N-2}+m_{N-3}\notag\\
&\vdots\notag\\
a_1-t_1+t_0 & = a_1-m_1+m_0\notag\\
a_0-t_0+t_{N-1} & = a_0-m_0+m_{N-1}\notag\\
t_0+t_1+...+t_{N-1} &=m=m_0+m_1+...+m_{N-1}.\notag
\end{align}
Which is equivalent to the system 
\begin{align}
 -t_{N-1}+t_{N-2} & =-m_{N-1}+m_{N-2}\notag\\
 -t_{N-2}+t_{N-3} & =-m_{N-2}+m_{N-3}\notag\\
&\vdots\notag\\
-t_1+t_0 & = -m_1+m_0\notag\\
-t_0+t_{N-1} & = -m_0+m_{N-1}\notag\\
t_0+t_1+...+t_{N-1} &=m=m_0+m_1+...+m_{N-1},\notag
\end{align}
and it can be easily seen that this system has as unique solution $t_i=m_i$ for $i=0,1,...,N-1$. 

Therefore we have that $\sigma\hat{a}=\hat{a}=\sigma'\hat{a}$, from Equations \eqref{E:sigma y} and \eqref{E:sigma y'} 
we have that $\sigma y=\hat{a}=\sigma' y'$ and from Equations \eqref{E:sigmaz} and \eqref{E:cilin} we have that
$\sigma z=\hat{a}=\sigma' z'$. Hence by letting $\gamma=\sigma^{-1}\sigma'$ we get 
\[  \gamma\hat{a}=\hat{a} \qquad \gamma y'=y \qquad \text{and}\qquad\gamma z'=z. \]
Thus we have $M_{[a],[b]}^{(k)[c]}=1$.
\end{proof}

\begin{corollary} \label{C:char}
Let $[a]\in\mathcal{O}$, assume that $\hat a=\left( (N-1)^{a_{N-1}},...,1^{a_1},0^{a_0}\right)$ and let 
$[b]=[(1^m,0^{k-m})]$. Then for $[c]\in\mathcal{O}$ we have
\[ M_{[a],[b]}^{(k)[c]}=
\begin{cases} 
1, &\text{if $\hat{c} =\left((N-1)^{a_{N-1}-m_{N-1}+m_{N-2}},...,1^{a_1-m_1+m_0},0^{a_0-m_0+m_{N-1}}\right)$, }\\ 
&\text{ for some integers $m_0$, $m_1$,...,$m_{N-1}$ such that $\displaystyle\sum_{i=0}^{N-1}m_i=m$}\\
&  \text{ and  } \quad 0\leq m_i\leq a_i\quad\text{, for  } \quad 0\leq i\leq N-1 \\
0, &\text{otherwise.}
\end{cases}
\]
\end{corollary}

We finish the chapter with an example showing that the coefficient $M_{[a],[b]}^{(k)[c]}$ could be more than 1.

\begin{example} \label{E:last}
Let $a=(2,1,0) \in \mathbb{Z}_3 ^3$.  Then we have that $\hat a =(2,1,0)$ and 
$[a]=\{ (2,1,0), (2,0,1), (0,2,1), (0,1,2), (1,0,2), (1,2,0)\}$. To compute $[a]\times [a]$, we remove all redundancies from
 the list:
\begin{align} 
  &(2,1,0)+ (2,1,0) = (1,2,0) , \notag\\
           &(2,1,0)+ (2,0,1) = (1,1,1) ,\notag\\
           &(2,1,0)+ (0,2,1) = (2,0,1) ,\notag \\ 
           &(2,1,0)+ (0,1,2) = (2,2,2) ,\notag \\
           &(2,1,0)+ (1,0,2) = (0,1,2) ,\notag \\
           &(2,1,0)+ (1,2,0) = (0,0,0) .\notag 
\end{align}
We can see that there are no redundant equations in that list, so we get
\begin{equation} \label{E:example23}
 [(2,1,0)]\times [(2,1,0)]=[(2,2,2)]+[(1,1,1)]  +3[(2,1,0)]+[(0,0,0)]. 
\end{equation}
Note that the orbit $[(2,1,0)]$ has multiplicity 3 in the product.
\end{example}

\chapter{Young diagrams and Pieri rules }

\section{Symmetric polynomials and fusion algebras } \label{S:sympol}

Given a dominant integral weight of $\mathfrak{g}=sl_N$, $\mu=\sum_{j=1}^{N-1} a_j\lambda_j$, we have associated to it a 
partition, denoted by  $(\mu)$, given in \eqref{E:pw} by:
\begin{equation} \label{E:tab}
(\mu)= \left(\sum_{j=1}^{N-1} a_j,\sum_{j=2}^{N-1} a_j,...,a_{N-1}\right)=\left(\mu_1,\mu_2,...,\mu_{N-1}\right)
\end{equation}
where $\mu_t$, $1\leq t\leq N-1$, is the $t^{th}$ part of the partition $(\mu)$. To such partition we can associate a Young 
diagram which is defined as the set of unit squares centered at the points $(s,t)\in \mathbb{Z}^2$ for 
$1\leq s\leq\mu _t$ and $1\leq t\leq l(\mu )$, where $l(\mu )$ denotes the length of $(\mu)$, the largest value of $t$ such 
that $\mu_t\ne 0$.  The Young diagram associated to the partition 
\eqref{E:tab} is given below.
\vskip10pt
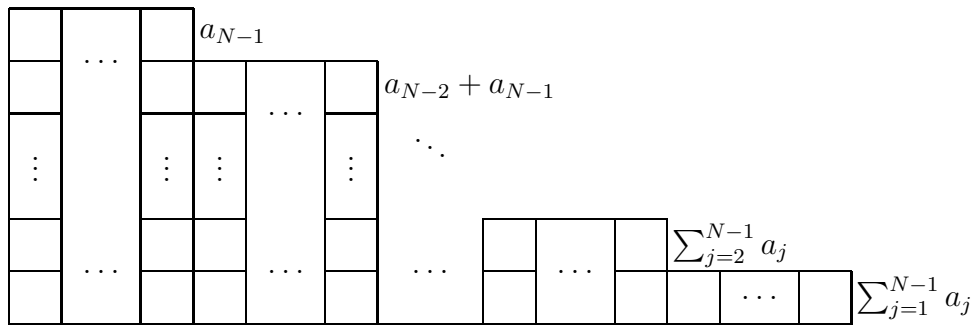
 \begin{xy}
<0.7cm,0cm>:
(-1.28,0); (0,0) **{ };
(0,0); (0,6) **@{-}, (1,0) **@{-};
(1,6); (0,6) **@{-}, (1,0) **@{-};
(0,1); (1,1) **@{-},
(0,2); (1,2) **@{-},
(0,4); (1,4) **@{-},
(0,5); (1,5) **@{-},
(1,0); (2.5,0) **@{-}, 
(1,6); (2.5,6) **@{-}, 
(1.75,5) *{\dots};
(1.75,1) *{\dots};
(0.5,3.1) *{\vdots};
(3,3.1) *{\vdots};
(4,3.1) *{\vdots};
(6.5,3.1) *{\vdots};
(4.25,5.5) *{ a_{N-1}};
(8.75,4.5) *{ a_{N-2}+a_{N-1}};
(13.7,1.5) *{\sum_{j=2}^{N-1} a_j};
(17.2,0.5) *{\sum_{j=1}^{N-1} a_j};
(2.5,0); (2.5,6) **@{-}, (4.5,0) **@{-};
(3.5,6); (2.5,6) **@{-}, (3.5,0) **@{-};
(2.5,1); (4.5,1) **@{-},
(2.5,2); (4.5,2) **@{-},
(2.5,4); (4.5,4) **@{-},
(2.5,5); (4.5,5) **@{-},
(4.5,0); (4.5,5) **@{-},
(5.25,1) *{\dots};
(5.25,4) *{\dots};
(4.5,0); (6,0) **@{-},
(4.5,5); (6,5) **@{-},
(6,0); (6,5) **@{-}, (7,0) **@{-};
(7,5); (6,5) **@{-}, (7,0) **@{-};
(6,1); (7,1) **@{-},
(6,2); (7,2) **@{-},
(6,4); (7,4) **@{-},
(8,1) *{\dots}; 
(7,0); (9,0) **@{-},
(8,3.5) *{\ddots};
(9,0); (9,2) **@{-}, (10,0) **@{-};
(10,2); (9,2) **@{-}, (10,0) **@{-};
(9,1); (10,1) **@{-},
(10,0); (11.5,0) **@{-},
(10,2); (11.5,2) **@{-},
(10.75,1) *{\dots};
(11.5,0); (11.5,2) **@{-}, (13.5,0) **@{-};
(12.5,2); (11.5,2) **@{-}, (12.5,0) **@{-};
(11.5,1); (13.5,1) **@{-},
(13.5,0); (13.5,1) **@{-},
(13.5,0); (15,0) **@{-},
(13.5,1); (15,1) **@{-},
(14.25,0.5) *{\dots};
(15,0); (15,1) **@{-}, (16,0) **@{-};
(16,1); (15,1) **@{-}, (16,0) **@{-};
\end{xy}
\vskip10pt
\noindent
where the label at the end of every row means the length of the row.

\begin{definition}
Given a Young diagram associated to a partition $(\mu)=\\ \left(\mu_1,\mu_2,...,\mu_N\right)$, we define the 
\textbf{conjugate Young diagram} to be the diagram obtained from the Young diagram of $(\mu)$ by interchanging rows and 
columns. The partition associated to the conjugate Young diagram will be denoted by 
$\mathbf{(\tilde\mu)=\left(\tilde\mu_1,\tilde\mu_2,...,\tilde\mu_t\right)}$ where $t=\mu_1$ and will be called
\textbf{the conjugate of $(\mu)$}.
\end{definition}

\begin{example}
Let $(\mu)=(5,4,1,1)$. Below we have the diagrams of $(\mu)$ and its conjugate $(\tilde\mu)$.
\vskip20pt
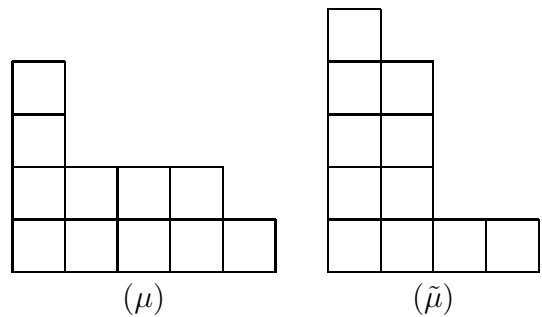
 \begin{xy}
<0.7cm,0cm>:
(-5.88,0); (0,0) **{ };
(0,0); (0,4) **@{-}, (5,0) **@{-};
(1,4); (1,0) **@{-}, (0,4) **@{-};
(5,1); (5,0) **@{-}, (0,1) **@{-};
(4,2); (4,0) **@{-}, (0,2) **@{-};
(0,3); (1,3) **@{-};
(2,2); (2,0) **@{-};
(3,2); (3,0) **@{-};
(6,0); (6,5) **@{-}, (10,0) **@{-};
(7,5); (6,5) **@{-}, (7,0) **@{-};
(10,1); (6,1) **@{-}, (10,0) **@{-};
(8,4); (6,4) **@{-}, (8,0) **@{-};
(6,2); (8,2) **@{-};
(6,3); (8,3) **@{-};
(9,0); (9,1) **@{-};
(2.5,-0.5) *{(\mu)};
(8,-0.5) *{(\tilde\mu)};
\end{xy}
\vskip10pt
Then the conjugate of $(\mu)$ is the partition $(\tilde\mu)=(4,2,2,2,1)$. Note that the length of $(\tilde\mu)$ equals 
$\mu_1$.
\end{example}

The following definitions are important to describe the product of symmetric polynomials.

\begin{definition}
If $(\nu)$ and $(\mu)$ are partitions so that $\mu_i\leq\nu_i$ for all $i$, the set difference between $(\nu)$ and
 $(\mu)$ is called a \textbf{skew partition} denoted by $(\nu) /(\mu)$ and its diagram is called a \textbf{skew diagram}.
\end{definition}

\begin{example} 
Let $(\mu)=(3,2,1,0)$ and  $(\nu)=(5,4,1,1)$, the Young diagrams for $(\nu)$, $(\mu)$ and $(\nu) /(\mu)$ are below

\vskip20pt
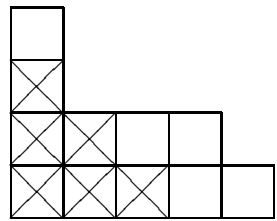
 \begin{xy}
<0.7cm,0cm>:
(-8.38,0); (0,0) **{ };
(0,0); (0,4) **@{-}, (5,0) **@{-}, (2,2) **@{-};
(1,4); (1,0) **@{-}, (0,4) **@{-};
(5,1); (5,0) **@{-}, (0,1) **@{-};
(4,2); (4,0) **@{-}, (0,2) **@{-};
(0,3); (1,3) **@{-}, (3,0) **@{-};
(2,2); (2,0) **@{-};
(3,2); (3,0) **@{-};
(2,1); (1,0) **@{-};
(3,1); (2,0) **@{-};
(1,2); (0,1) **@{-}; 
(1,3); (0,2) **@{-};
(0,2); (2,0) **@{-};
(0,1); (1,0) **@{-};
\end{xy}
\vskip10pt\noindent
where $(\nu)$ is the whole diagram, $(\mu)$ is the diagram formed by the crossed boxes and the skew diagram $(\nu) /(\mu)$ 
is formed by the empty boxes.
\end{example}

\begin{definition}
A skew diagram is called an \textbf{m-column strip} if it has m boxes with at most one box in each row, and is called an 
\textbf{m-row strip} if it has m boxes with at most one box in each column.
\end{definition}

\begin{example}
Let $(\mu)=(3,2,1,0)$ and  $(\nu)=(4,3,1,1)$

\vskip10pt
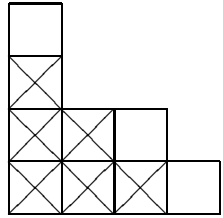
 \begin{xy}
<0.7cm,0cm>:
(-8.88,0); (0,0) **{ };
(0,0); (0,4) **@{-}, (4,0) **@{-}, (2,2) **@{-};
(1,4); (1,0) **@{-}, (0,4) **@{-};
(4,1); (4,0) **@{-}, (0,1) **@{-};
(3,2); (3,0) **@{-}, (0,2) **@{-};
(0,3); (1,3) **@{-}, (3,0) **@{-};
(2,2); (2,0) **@{-};
(2,1); (1,0) **@{-};
(3,1); (2,0) **@{-};
(0,1); (1,2) **@{-}, (1,0) **@{-};
(0,2); (1,3) **@{-}, (2,0) **@{-};
\end{xy}
\vskip10pt
The skew diagram $(\nu) /(\mu)$ in this example is both a 3-row strip and a 3-column strip.
\end{example}

\begin{definition}
A \textbf{tableau} of shape $(\nu) /(\mu)$ is a filling of the diagram $(\nu) /(\mu)$ with positive integers nondecreasing 
in rows and strictly increasing in columns. The \textbf{content} of a tableau is the sequence $(b_1,b_2,...)$ where $i$ 
appears $b_i$ times in the filling for every $i\geq 1$.
\end{definition}

\begin{example}
Let $(\mu)=(2,2)$ and  $(\nu)=(5,4,1,1)$. The following fillings of the skew diagram $(\nu)/(\mu)$ give two tableaux
\vskip10pt
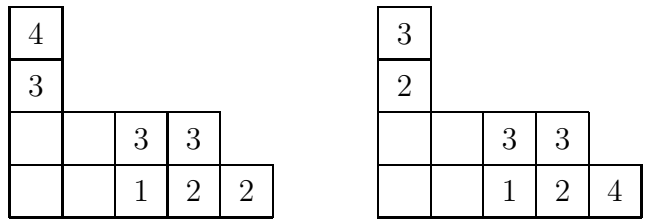
 \begin{xy}
<0.7cm,0cm>:
(-4.88,0); (0,0) **{ };
(0,0); (0,4) **@{-}, (5,0) **@{-};
(1,4); (1,0) **@{-}, (0,4) **@{-};
(5,1); (5,0) **@{-}, (0,1) **@{-};
(4,2); (4,0) **@{-}, (0,2) **@{-};
(0,3); (1,3) **@{-};
(2,2); (2,0) **@{-};
(3,2); (3,0) **@{-};
(2.5,0.5) *{1};
(3.5,0.5) *{2};
(4.5,0.5) *{2};
(3.5,1.5) *{3};
(2.5,1.5) *{3};
(0.5,2.5) *{3};
(0.5,3.5) *{4};
(7,0); (7,4) **@{-}, (12,0) **@{-};
(8,4); (8,0) **@{-}, (7,4) **@{-};
(12,1); (12,0) **@{-}, (7,1) **@{-};
(11,2); (11,0) **@{-}, (7,2) **@{-};
(7,3); (8,3) **@{-};
(10,2); (10,0) **@{-};
(9,2); (9,0) **@{-};
(9.5,0.5) *{1};
(10.5,0.5) *{2};
(11.5,0.5) *{4};
(10.5,1.5) *{3};
(9.5,1.5) *{3};
(7.5,3.5) *{3};
(7.5,2.5) *{2};
\end{xy}
\vskip10pt
These two tableaux have both content (1,2,3,1). 
\end{example}
Young diagrams are closely related to the algebra of symmetric polynomials, since a very important basis for this algebra is 
indexed by partitions. We describe the algebra of symmetric polynomials as well as its basis formed by Schur polynomials in
 the following paragraph.

The algebra of symmetric polynomials in $N$ variables is the algebra of polynomials \newline 
$f\in\mathbb{Q}[x_1,x_2,...,x_N]$ invariant under the action of the symmetric group $S_N$ that permutes the variables. This
algebra is denoted by $\mathbf{\Lambda _N=\mathbb{Q}[x_1,x_2,...,x_N]^{S_N}}$. 

For $0<m\leq N$ define the elementary symmetric polynomial
\[ e_m =\sum_{1\leq i_1<...<i_m\leq N}x_{i_1}... x_{i_m},\]
and $e_m=0$ for $m>N$. 

For $m>0$ we define the homogeneous symmetric function 
\[ h_m = \sum_{1\leq i_1 ...\leq i_m\leq N}x_{i_1}... x_{i_m}\ , \]
\noindent
and we define $e_m=0=h_m$ for $m<0$. We also define $e_0=1=h_0$.

A basis for $\Lambda_N$ is given by the Schur polynomials $S_{(\mu)}$ indexed by partitions $(\mu)=(\mu_1,\mu_2,...,\mu_N)$
 with $l(\mu)\leq N$, defined by
\begin{equation} \label{E:schur}
 S_{(\mu)} =det(h_{\mu _i -i+j})=  \begin{vmatrix}
    h_{\mu_1} & h_{\mu_1 +1} & \dots & h_{\mu_1 +N-1} \\
    h_{\mu_2 -1} & h_{\mu_2} & \dots & h_{\mu_2 +N-2} \\
    \vdots & \vdots & \ddots & \vdots \\
    h_{\mu_N -N+1} & h_{\mu_N -N+2} & \dots & h_{\mu_N }
 \end{vmatrix} 
\end{equation}
if $1\leq l(\mu)\leq N$. Note that $S_{(\mu)}=h_m$ if $(\mu) =(m)$ consists of a single part.

The elementary symmetric polynomials are also generators of the algebra of $\Lambda _N$, since Schur polynomials can also be
expressed as a determinant of them. If $(\mu)=(\mu_1,\mu_2,...,\mu_N)$ and its conjugate partition is 
$(\tilde\mu)=(\tilde\mu_1,\tilde\mu_2,...,\tilde\mu_t)$ where $t=\mu_1$, then the Schur polynomial $S_{(\mu)}$ is also given 
by
\begin{equation} \label{E:schur1}
 S_{(\mu)} =det(e_{\tilde\mu _i -i+j})=  \begin{vmatrix}
    e_{\tilde\mu_1} & e_{\tilde\mu_1 +1} & \dots & e_{\tilde\mu_1 +t-1} \\
    e_{\tilde\mu_2 -1} & e_{\tilde\mu_2} & \dots & e_{\tilde\mu_2 +t-2} \\
    \vdots & \vdots & \ddots & \vdots \\
    e_{\tilde\mu_t -t+1} & e_{\tilde\mu_t -t+2} & \dots & e_{\tilde\mu_t }
 \end{vmatrix}. 
\end{equation}
Note that $S_{(\mu)}=e_m$ if $(\mu) =(1^m)$.

Equations \eqref{E:schur} and \eqref{E:schur1} are known in the literature as the Jacobi-Trudi determinants.

It was proved by Goodman and Wenzl \cite{GW} that the fusion algebra $\mathcal{F} (A_{N-1},k)$ associated to 
$\widehat{sl}_N$ on 
level $k$ is isomorphic to the quotient algebra of symmetric polynomials $\Lambda_N/I^{(N,k)}$, where 
\[ I^{(N,k)}=\langle S_{(1^N)} -1,S_{(\mu)} \mid \mu_1 -\mu_N =k+1 \rangle. \]

This condition on partitions in $I^{(N,k)}$ limits the shape of Young diagrams that are relevant to the quotient algebra
 $\Lambda _N/I^{(N,k)}$.
The following definition exactly describes these diagrams.
\begin{definition}
Let $(\mu)$ be a partition so that $\mu_1 -\mu_N\leq k$, then we say that $(\mu)$ is $\mathbf{(N,k)}$-\textbf{restricted}. 
The set of $(N,k)$-restricted partitions will be denoted by $\mathbf{\Pi^{(N,k)}}$.
\end{definition}

\begin{example}
The Young diagrams for the partitions in the definition are exactly the ones whose distance between the first column of 
height $N-1$ and the last column, is less than or equal to $k$, as the next diagram shows.
\vskip10pt
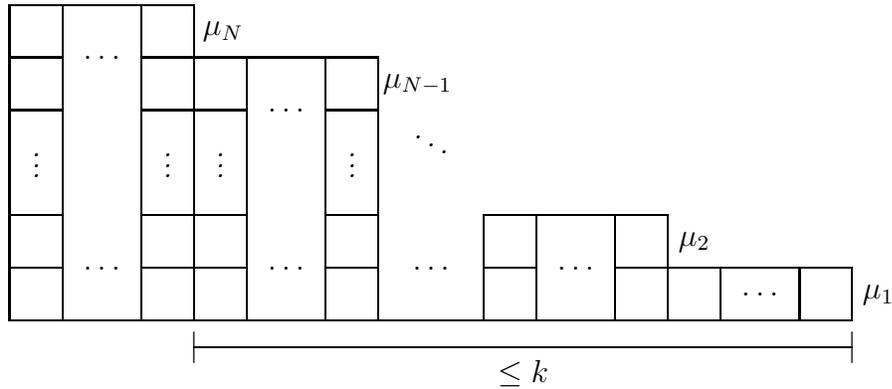
 \begin{xy}
<0.7cm,0cm>:
(-2.28,0); (0,0) **{ };
(0,0); (0,6) **@{-}, (1,0) **@{-};
(1,6); (0,6) **@{-}, (1,0) **@{-};
(0,1); (1,1) **@{-},
(0,2); (1,2) **@{-},
(0,4); (1,4) **@{-},
(0,5); (1,5) **@{-},
(1,0); (2.5,0) **@{-}, 
(1,6); (2.5,6) **@{-}, 
(1.75,5) *{\dots};
(1.75,1) *{\dots};
(0.5,3.1) *{\vdots};
(3,3.1) *{\vdots};
(4,3.1) *{\vdots};
(6.5,3.1) *{\vdots};
(4.05,5.5) *{\mu_N};
(7.75,4.5) *{\mu_{N-1}};
(13,1.5) *{\mu_2};
(16.5,0.5) *{\mu_1};
(2.5,0); (2.5,6) **@{-}, (4.5,0) **@{-};
(3.5,6); (2.5,6) **@{-}, (3.5,0) **@{-};
(2.5,1); (4.5,1) **@{-},
(2.5,2); (4.5,2) **@{-},
(2.5,4); (4.5,4) **@{-},
(2.5,5); (4.5,5) **@{-},
(4.5,0); (4.5,5) **@{-},
(5.25,1) *{\dots};
(5.25,4) *{\dots};
(4.5,0); (6,0) **@{-},
(4.5,5); (6,5) **@{-},
(6,0); (6,5) **@{-}, (7,0) **@{-};
(7,5); (6,5) **@{-}, (7,0) **@{-};
(6,1); (7,1) **@{-},
(6,2); (7,2) **@{-},
(6,4); (7,4) **@{-},
(8,1) *{\dots}; 
(7,0); (9,0) **@{-},
(8,3.5) *{\ddots};
(9,0); (9,2) **@{-}, (10,0) **@{-};
(10,2); (9,2) **@{-}, (10,0) **@{-};
(9,1); (10,1) **@{-},
(10,0); (11.5,0) **@{-},
(10,2); (11.5,2) **@{-},
(10.75,1) *{\dots};
(11.5,0); (11.5,2) **@{-}, (13.5,0) **@{-};
(12.5,2); (11.5,2) **@{-}, (12.5,0) **@{-};
(11.5,1); (13.5,1) **@{-},
(13.5,0); (13.5,1) **@{-},
(13.5,0); (15,0) **@{-},
(13.5,1); (15,1) **@{-},
(14.25,0.5) *{\dots};
(15,0); (15,1) **@{-}, (16,0) **@{-};
(16,1); (15,1) **@{-}, (16,0) **@{-};
(3.5,-.5); (16,-.5) **@{-};
(3.5,-.5) *{|};
(16,-.5) *{|};
(9.75,-1) *{\leq k};
\end{xy}
\vskip10pt
\end{example}

After giving some definitions we present below some results of Goodman and Wenzl, (see \cite{GW} or \cite{Tu}) which 
provide the multiplication in the fusion algebra $\mathcal{F}(A_{N-1},k)$.

\begin{definition} 
 We call a skew diagram $(\nu) /(\mu)$ a $\mathbf{(N,k)}$-\textbf{cylindric m-row strip} if it is an $m$-row strip for some 
$m$ and $\nu_1-\mu_N\leq k$.
\end{definition}

\begin{example}
Let $\mu=(3,2,1)$ and $\nu=(4,3,2)$ be partitions, $N=3$ and $k=3$. The skew partition $\nu/\mu$ is a (3,3)-cylindric 3-row 
strip, since $\nu_1-\mu_3=3$. This can be visualized from the diagram
\vskip20pt
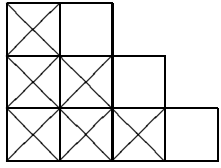
 \begin{xy}
<0.7cm,0cm>:
(-8.88,0); (0,0) **@{ };
(0,0); (0,3) **@{-}, (4,0) **@{-}, (2,2) **@{-};
(2,3); (0,3) **@{-}, (2,0) **@{-};
(1,3); (1,0) **@{-}, (0,2) **@{-};
(3,2); (0,2) **@{-}, (3,0) **@{-};
(4,1); (0,1) **@{-}, (4,0) **@{-};
(0,1); (1,2) **@{-};
(1,0); (2,1) **@{-}, (0,1)**@{-};
(2,0); (3,1) **@{-}, (0,2)**@{-};
(0,3); (3,0) **@{-};
\end{xy}
\vskip10pt
In fact $\nu/\mu$ is a $(3,k)$-cylindric 3-row strip for any $k\geq 3$ and it is not a $(3,2)$-cylindric 3-row strip.
\end{example}

The isomorphism $\Lambda_N/I^{(N,k)}\cong \mathcal{F}(A_{N-1},k)$ gives the following dictionary

\begin{align} \label{E:dict}
 \Lambda_N^{S_N}/I^{(N,k)} & \longleftrightarrow \mathcal{F}(A_{N-1},k) \\
 S_{(\mu)} & \longleftrightarrow V^{\mu} \notag\\
S_{(\mu)}S_{(\lambda)}=\sum_{(\nu)}N_{(\mu),(\lambda)}^{(k)(\nu)}S_{(\nu)} & \longleftrightarrow 
V^{\mu}\otimes_k V^{\lambda}=\bigoplus_{\nu}N_{\mu,\lambda}^{(k)\nu}V^{\nu}, \notag
\end{align}
where the left hand side is indexed by partitions $(\mu)=(\mu_1,...,\mu_N)\in\Pi^{(N,k)}$ and the right hand side is indexed 
by weights $\mu=\sum_{j=1}^{N-1}(\mu_j-\mu_{j+1})\lambda_j\in P_k^+$. The correspondence \eqref{E:wp} allows us to move from
left to right in the dictionary and the right inverse \eqref{E:pw} of \eqref{E:wp} allows us to move in the opposite 
direction in the dictionary. We also get an equality of structure constants in both rings, i.e., 
$N_{(\mu),(\lambda)}^{(k)(\nu)}=N_{\mu,\lambda}^{(k)\nu}$.

\section{Fusion Pieri rules }

From the equivalence relation $\thicksim$ defined in \eqref{E:er} we know that given a weight 
$\mu=\sum_{j=1}^{N-1}a_j\lambda_j\in P_k^+$ there are infinitely many partitions in $\Pi^{(N,k)}$ that correspond to $\mu$, 
since the correspondence \eqref{E:pw} is onto but not one-to-one. The next lemma says that the Schur polynomials 
corresponding to partitions that are equivalent under $\thicksim$ are equal in the quotient $\Lambda_N/I^{(N,k)}$. 
Just recall that two partitions $(\mu)=(\mu_1,...,\mu_N)$ and $(\nu)=(\nu_1,...,\nu_N)$ are equivalent under $\thicksim$ if 
and only if there exists a positive integer $c$ such that $\mu_i=\nu_i+c$ for $i=1,...,N$. Then it would be sufficient to 
prove the following.

\begin{lemma} \label{L:eq}
Let $(\mu)=(\mu_1,...,\mu_N)\in\Pi^{(N,k)}$, then in $\Lambda_N/I^{(N,k)}$ we have
\[ S_{(\mu)}=S_{(\mu_1-\mu_N,...,\mu_{N-1}-\mu_N)}. \]
\end{lemma}

To prove this lemma we need the following result due to Goodman-Wenzl (see \cite{GW} or \cite{Tu}). 

\begin{theorem} \label{T:em} \cite{GW} (\textbf{Fusion Pieri rule} for multiplication by $e_m$) 

Let $(\mu)\in\Pi^{(N,k)}$ and $m\leq N$. Then in $\Lambda_N/I^{(N,k)}$ we have
\begin{equation} \label{E:em}
 S_{(\mu)}e_m=\sum_{\substack{(\nu)\in\Pi^{(N,k)},\\ (\nu)/(\mu) \text{ is an m-column strip}}}S_{(\nu)}. 
\end{equation}
\end{theorem}

\begin{corollary} 
For $(\mu)=(\mu_1,...,\mu_N)\in\Pi^{(N,k)}$ we have the equality in $\Lambda_N/I^{(N,k)}$
\[ S_{(\mu)}=S_{(\mu_1+1,...,\mu_N+1)}. \]
\end{corollary}
\begin{proof}
From Theorem~\ref{T:em} we have
\[ S_{(\mu)}e_N=\sum_{\substack{(\nu)\in\Pi^{(N,k)},\\ (\nu) /(\mu) \text{ is an N-column strip}}}S_{(\nu)}. \]
The partition $(\mu)$ has length at most $N$, so its Young diagram has at most $N$ rows, therefore the only 
$(\nu)\in\Pi^{(N,k)}$ such that $(\nu) /(\mu)$ is an N-column strip is the one obtained by adding one box on each row, so we
 get
\[ S_{(\mu)}e_N=S_{(\mu_1+1,...,\mu_N+1)}. \]
Now, since $e_N-1$ is in $I^{(N,k)}$, we get the desired equality in $\Lambda_N/I^{(N,k)}$.
\end{proof}

Lemma~\ref{L:eq} is a direct consequence of this corollary. It can be also proved directly using \eqref{E:schur1}.

\begin{example}
Let $N=4$, $k=3$ and $(\mu)=(5,4,4,3)$. From Lemma~\ref{L:eq} we have that in $\Lambda_4/I^{(4,3)}$
\[ S_{(\mu)}= S_{(2,1,1)}. \]
We can actually see this equality by using the definition of the Schur polynomial in terms of the elementary symmetric 
polynomials and the fact that $e_4=1$ in the quotient algebra $\Lambda_4/I^{(4,3)}$. The conjugate of $(\mu)=(5,4,4,3)$ is 
the partition $(\tilde\mu)=(4,4,4,3,1)$ so from the Jacobi-Trudi determinant given by \eqref{E:schur1} we get 
\[ S_{(\mu)} =det(e_{\tilde\mu _i -i+j})=  \begin{vmatrix}
    e_4 & e_5 & e_6 & e_7 & e_8\\
    e_3 & e_4 & e_5 & e_6 & e_7 \\
    e_2 & e_3 & e_4 & e_5 & e_6\\
     e_0 & e_1 & e_2 & e_3 & e_4\\
e_{-3} & e_{-2} & e_{-1} & 1 & e_1 
 \end{vmatrix} .\]
Using the relations $e_{-1}=e_{-2}=e_{-3}=e_5=e_6=e_7=e_8=0$ and $e_0=1$ in $\Lambda_4$ and $e_4=1$ in $\Lambda_4/I^{(4,3)}$,
this determinant simplifies to
\[ S_{(\mu)} =det(e_{\tilde\mu _i -i+j})=  \begin{vmatrix}
    1 & 0 & 0 & 0 & 0\\
    e_3 & 1 & 0 & 0 & 0 \\
    e_2 & e_3 & 1 & 0 & 0\\
     1 & e_1 & e_2 & e_3 & 1\\
0 & 0 & 0 & 1 & e_1 
 \end{vmatrix} =\begin{vmatrix} e_3 & 1\\1 & e_1 \end{vmatrix}.\]
The conjugate of the partition $(2,1,1)$ is the partition $(3,1)$, and from \eqref{E:schur1} we get 
\[ S_{(2,1,1)} =det(e_{\tilde\mu _i -i+j})=  \begin{vmatrix} e_3 & 1\\1 & e_1 \end{vmatrix}. \]
This example shows how Lemma~\ref{L:eq} yields the equality in a much simpler way. 
\end{example}

Lemma~\ref{L:eq} says that from a Young diagram in $\Pi^{(N,k)}$ we can remove all of the columns of height N, and the Schur 
polynomials corresponding to the original diagram and the new diagram are equal in $\Lambda_N/I^{(N,k)}$.

The next theorem illustrates another result of \cite{GW}, connected to the orbits product from 
Theorem \ref{T:orbitproduct}. We will see some examples of how the connection works in this chapter and we will give the 
proof in the next chapter.

\begin{theorem} \label{T:hm} \cite{GW} (\textbf{Fusion Pieri rule} for multiplication by $h_m$) 

Let $(\mu)\in\Pi^{(N,k)}$ and $m\leq k$. Then in $\Lambda_N/I^{(N,k)}$ we have
\begin{equation} \label{E:Pieri}
 S_{(\mu)}h_m=\sum_{\substack{(\nu)\in\Pi^{(N,k)},\\ (\nu) /(\mu)\text{ is an }(N,k)\text{-cylindric m-row strip}}}S_{(\nu)}. 
\end{equation}
\end{theorem}

\begin{example} \label{Ex:9}
Let $N=4$, $k=3$ and $(\mu)=(3,3,1,1)$. Using the previous theorem we get that $S_{(\mu)}h_1$ equals the sum of 
$S_{(\nu)}$'s so that $(\nu)$ is a $(4,3)$-cylindric 1-row strip. The diagram below shows all possible $(\nu)$'s.
\vskip20pt
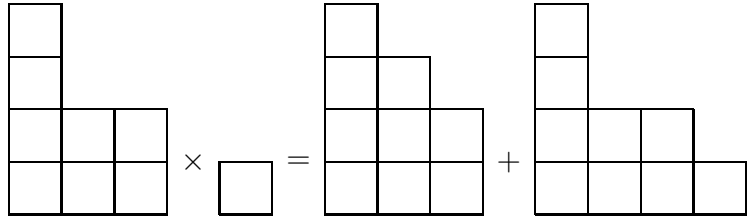
\begin{xy}
<0.7cm,0cm>:
(-3.88,0); (0,0) **@{ };
(0,0); (0,4) **@{-}, (3,0) **@{-};
(1,4); (0,4) **@{-}, (1,0) **@{-};
(0,3); (1,3) **@{-};
(3,2); (3,0) **@{-}, (0,2) **@{-};
(3,1); (0,1) **@{-};
(2,0); (2,2) **@{-};
(3.5,1), *{\times};
(4,0); (4,1) **@{-}, (5,0) **@{-};
(5,1); (4,1) **@{-}, (5,0) **@{-};
(5.5,1), *{=};
(6,0); (6,4) **@{-}, (9,0) **@{-};
(7,4); (6,4) **@{-}, (7,0) **@{-};
(8,3); (6,3) **@{-}, (8,0) **@{-};
(9,2); (6,2) **@{-}, (9,0) **@{-};
(9,1); (6,1) **@{-};
(9.5,1), *{+};
(10,0); (10,4) **@{-}, (14,0) **@{-};
(11,4); (10,4) **@{-}, (11,0) **@{-};
(14,1); (10,1) **@{-}, (14,0) **@{-};
(13,2); (10,2) **@{-}, (13,0) **@{-};
(10,3); (11,3) **@{-};
(12,0); (12,2) **@{-}
\end{xy}
\vskip20pt
Therefore in the algebra $\Lambda_4/I^{(4,3)}=\mathbb{Q}[x_1,x_2,x_3,x_4]^{S_4}/I^{(4,3)}$ we get that 
\[ S_{(3,3,1,1)}h_1=S_{(3,3,2,1)} + S_{(4,3,1,1)}. \]
From Lemma~\ref{L:eq} we know that we can remove the column of length 4 from each of the Young diagrams, so we get
\begin{equation} \label{E:9}
S_{(2,2)}h_1=S_{(2,2,1)} + S_{(3,2)}.
\end{equation}
Now using the Dictionary \eqref{E:dict}, this product translates into the fusion product in $\mathcal{F}(A_3,3)$
\[ V^{2\lambda_2}\otimes_3 V^{\lambda_1} = V^{\lambda_2 +\lambda_3}\oplus V^{\lambda_1 +2\lambda_2}. \]
\end{example}

\begin{example} \label{Ex:10}
Let $N=4$, $k=3$ and $(\mu)=(3,2,2,1)$. To find the product $S_{(\mu)}h_1$, we use the diagram below, which shows all 
diagrams $(\nu)\in\Pi^{(4,3)}$ where $(\nu)/(\mu)$ is a (4,3)-cylindric 1-row strip 
\vskip20pt
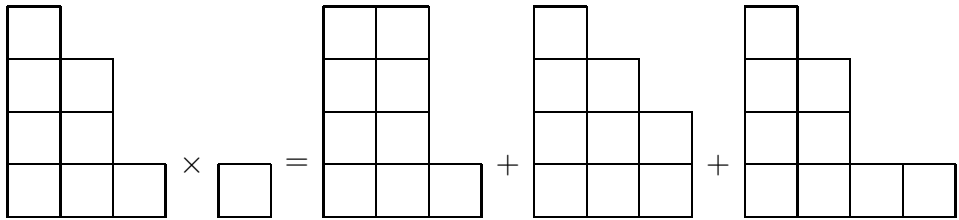
\begin{xy}
<0.7cm,0cm>:
(-1.88,0); (0,0) **@{ };
(0,0); (0,4) **@{-}, (3,0) **@{-};
(1,4); (0,4) **@{-}, (1,0) **@{-};
(2,3); (0,3) **@{-}, (2,0) **@{-};
(3,1); (0,1) **@{-}, (3,0) **@{-};
(2,2); (0,2) **@{-};
(3.5,1), *{\times};
(4,0); (4,1) **@{-}, (5,0) **@{-};
(5,1); (4,1) **@{-}, (5,0) **@{-};
(5.5,1), *{=};
(6,0); (6,4) **@{-}, (9,0) **@{-};
(8,4); (6,4) **@{-}, (8,0) **@{-};
(7,4); (7,0) **@{-};
(9,1); (6,1) **@{-}, (9,0) **@{-};
(8,2); (6,2) **@{-};
(8,3); (6,3) **@{-};
(9.5,1), *{+};
(10,0); (10,4) **@{-}, (13,0) **@{-};
(11,4); (10,4) **@{-}, (11,0) **@{-};
(12,3); (12,0) **@{-}, (10,3) **@{-};
(13,2); (10,2) **@{-}, (13,0) **@{-};
(13,1); (10,1) **@{-};
(13.5,1), *{+};
(14,0); (14,4) **@{-}, (18,0) **@{-};
(15,4); (14,4) **@{-}, (15,0) **@{-};
(16,3); (16,0) **@{-}, (14,3) **@{-};
(18,1); (14,1) **@{-}, (18,0) **@{-};
(17,1); (17,0) **@{-};
(14,2); (16,2) **@{-};
\end{xy}
\vskip20pt
Therefore, in the algebra $\Lambda_4/I^{(4,3)}$, we get that 
\[ S_{(3,2,2,1)}h_1=S_{(3,2,2,2)} + S_{(3,3,2,1)} + S_{(4,2,2,1)}. \]
By removing all columns of height 4 in each of the diagrams on both sides we get
\begin{equation} \label{E:10}
S_{(2,1,1)}h_1=S_{(1)} + S_{(2,2,1)} + S_{(3,1,1)}= h_1 + S_{(2,2,1)} + S_{(3,1,1)},
\end{equation}
which, by the Dictionary \eqref{E:dict}, translates into the fusion product in $\mathcal{F}(A_3,3)$ 
\[ V^{\lambda_1+\lambda_3}\otimes_3 V^{\lambda_1} = V^{\lambda_1}\oplus V^{\lambda_2 +\lambda_3}\oplus 
V^{2\lambda_1 +\lambda_3} \]
\end{example} 

The previous theorem can be simplified. Note that since $\mu$ is a weight of $A_{N-1}$, we can take the Young diagram for 
$(\mu)$ to be inside the rectangle $(N-1)\times k$, without loss of generality. Then  $(\nu)$ in the right hand side of 
\eqref{E:Pieri} has to be inside the rectangle $N\times k$. We prove this statement in the following lemma  after setting up 
some notation.

\begin{notation}
For a partition $(\mu)\in \Pi^{(N,k)}$, we write $(\mu)\subseteq r\times k$ when the Young diagram of $(\mu)$ is contained 
in the rectangle $r\times k$. 
\end{notation}

\begin{lemma} \label{L:nu1}
Let $(\mu)\subseteq (N-1)\times k$ and let $(\nu)\in\Pi^{(N,k)}$, where $(\nu) /(\mu)$ is a $(N,k)$-cylindric m-row
 strip for $m\leq k$,  then $(\nu)\subseteq N\times k$.
\end{lemma}
\begin{proof}
Let $(\mu)\subseteq (N-1)\times k$ and assume that $(\nu)\in\Pi^{(N,k)}$ where $(\nu) /(\mu)$ is a $(N,k)$-cylindric m-row 
strip. Thus $(\mu)$ is a partition of length at most $N-1$, say $(\mu)=(\mu_1,...,\mu_{N-1})$. Note that since 
$(\mu)\subseteq (N-1)\times k$ and $(\nu) /(\mu)$ is an m-row strip, then the height of any column of $(\nu)$ is at most $N$.
That is, $(\nu)$ is a partition of length at most $N$, say $(\nu)=(\nu_1,...,\nu_N)$.

We claim that $(\nu)$ cannot have rows of length $>k$, because otherwise $\nu_1 -\mu_N=\nu_1>k$ which means $(\nu)/(\mu)$ is 
not $(N,k)$-cylindric. Therefore $(\nu)$ is inside the rectangle $N\times k$.

\end{proof}

\begin{remark} \label{R:nu2}
The condition $(\nu) /(\mu)$ is $(N,k)$-cylindric is implied from the assumption that $(\nu)\subseteq N\times k$.
This is clear because $\nu_1\leq k$ implies $\nu_1-\mu_N\leq k$.
\end{remark}

Lemma~\ref{L:nu1} and Remark~\ref{R:nu2} allow us to rewrite Theorem~\ref{T:hm} and equation~\eqref{E:Pieri} as follows:

\begin{theorem} \label{T:hm2} (\textbf{Fusion Pieri rule} for multiplication by $h_m$) 

Let $(\mu)\subseteq(N-1)\times k$ and let $m\leq k$. Then in $\Lambda_N/I^{(N,k)}$ we have
\begin{equation} \label{E:p2}
 S_{(\mu)}h_m=\sum_{\substack{(\nu)\subseteq{N\times k},\\ (\nu) /(\mu)\text{ is an m-row strip } }}S_{(\nu)}.
\end{equation}
\end{theorem}

\begin{remark} \label{R:red}
Theorem~\ref{T:hm} seems a little more general than this last theorem, since it applies to partitions $(\mu)$ in 
$\Pi^{(N,k)}$, but in general if $(\mu)=(\mu_1,\mu_2,...,\mu_N)\in\Pi^{(N,k)}$ we can apply Theorem~\ref{T:hm2} to the 
partition $(\mu_1-\mu_N,\mu_2-\mu_N,...,\mu_{N-1}-\mu_N)\subseteq(N-1)\times k$ and we get the same result.
\end{remark}

Theorem~\ref{T:hm2} will be useful to prove our main result (Theorem \ref{T:main}). Here we show some examples of how this 
theorem applies.

\begin{example}
Let $(\mu)=(2,1)$, $N=3$ and $k=3$. Using Theorem~\ref{T:hm2} to compute $S_{(\mu)}h_2$ we get
\vskip20pt
\begin{xy}
<0.7cm,0cm>:
(-2.88,0); (0,0) **{ };
(0,0); (0,2) **@{-}, (2,0) **@{-};
(2,1); (0,1) **@{-}, (2,0) **@{-},
(1,2); (1,0) **@{-}, (0,2) **@{-};
(2.5,1) *{\times};
(3,0); (3,1) **@{-}, (5,0) **@{-};
(5,1); (3,1) **@{-}, (5,0) **@{-};
(4,1); (4,0) **@{-};
(5.5,1) *{=}; (8.5,1) *{+}; (12.5,1) *{+};
(6,0); (6,3) **@{-}, (8,0) **@{-};
(8,2); (6,2) **@{-}, (8,0) **@{-};
(7,3); (6,3) **@{-}, (7,0) **@{-};
(8,1); (6,1) **@{-};
(9,0); (9,3) **@{-}, (12,0) **@{-};
(12,1); (9,1) **@{-}, (12,0) **@{-};
(10,3); (10,0) **@{-}, (9,3) **@{-};
(11,1); (11,0) **@{-};
(10,2); (9,2) **@{-};
(13,0); (13,2) **@{-}, (16,0) **@{-};
(16,1); (13,1) **@{-}, (16,0) **@{-};
(15,2); (15,0) **@{-}, (13,2) **@{-};
(14,2); (14,0) **@{-};
\end{xy}
\vskip10pt
\noindent
that is, in the algebra $\Lambda_3/I^{(3,3)}$, we have that:
\[ S_{(2,1)}h_2=S_{(2,2,1)}+S_{(3,1,1)}+S_{(3,2)}. \]
Now, the weight associated to the partition $(\mu)=(2,1)$ is $\mu =\lambda_1+\lambda_2$ and the weight associated to
$h_2=S_{(2,0,0)}$ is $2\lambda_1$. Then the above equation translates into the fusion product in $\mathcal{F}(A_2,3)$
\[ V^{\lambda_1+\lambda_2}\otimes_3 V^{2\lambda_1}= V^{\lambda_2}\oplus V^{2\lambda_1}\oplus V^{\lambda_1+2\lambda_2}. \]
Under the correspondence \eqref{E:ow} we have that the orbits corresponding to the weights $\lambda_1+\lambda_2$, 
$2\lambda_1$, $\lambda_2$ and $\lambda_1+2\lambda_2$ are respectively [(2,1,0)], [(1,1,0)], [(2,0,0)] and [(2,2,1)]. Looking 
at \eqref{E:orbpieri} in Example~\ref{E:or1} we see that this fusion product agrees with the orbit calculation from that 
example. This is not a coincidence, it follows from \eqref{E:FW} proved by Feingold and Weiner in \cite{FW}.
\end{example}

In the next example we see another fusion product done via Pieri rules whose outcome agrees with an orbit calculation for 
$N=4$. This does not follow from \eqref{E:FW}, but it is not a coincidence either as we will see in the next chapter.
\vskip10pt
\begin{example}
Let $(\mu)=(3,2,1)$, $N=4$ and $k=3$. Using the Theorem~\ref{T:hm2} we get
\vskip20pt
 \begin{xy}
<0.7cm,0cm>:
(-1.88,0); (0,0) **{ };
(0,0); (3,0) **@{-}, (0,3) **@{-};
(1,3); (1,0) **@{-}, (0,3) **@{-};
(2,2); (2,0) **@{-}, (0,2) **@{-};
(3,1); (3,0) **@{-}, (0,1) **@{-};
(3.5,1) *{\times};
(4,0); (6,0) **@{-}, (4,1) **@{-};
(6,1); (6,0) **@{-}, (4,1) **@{-};
(5,0); (5,1) **@{-};
(6.5,1) *{=};
(7,0); (10,0) **@{-}, (7,3) **@{-};
(9,3); (9,0) **@{-}, (7,3) **@{-};
(10,2); (10,0) **@{-}, (7,2) **@{-};
(10,1); (7,1) **@{-};
(8,3); (8,0) **@{-};
(10.5,1) *{+};
(11,0); (14,0) **@{-}, (11,4) **@{-};
(12,4); (12,0) **@{-}, (11,4) **@{-};
(12,3); (11,3) **@{-};
(14,2); (14,0) **@{-}, (11,2) **@{-};
(14,1); (11,1) **@{-};
(13,0); (13,2) **@{-};
(14.5,1) *{+};
(15,0); (18,0) **@{-}, (15,4) **@{-};
(16,4); (16,0) **@{-}, (15,4) **@{-};
(17,3); (17,0) **@{-}, (15,3) **@{-};
(17,2); (15,2) **@{-};
(18,1); (18,0) **@{-}, (15,1) **@{-};
\end{xy}
\vskip10pt
This shows that $S_{(3,2,1)}h_2=S_{(3,3,2)}+S_{(3,3,1,1)}+S_{(3,2,2,1)}$ where the computation is done in the algebra 
$\Lambda_4/I^{(4,3)}$. By taking away columns of length 4 we get
\[ S_{(3,2,1)}h_2=S_{(3,3,2)}+S_{(2,2)}+S_{(2,1,1)}. \]
By the Dictionary \eqref{E:dict} we get the fusion product
\[ V^{\lambda_1+\lambda_2+\lambda_3}\otimes_3 V^{2\lambda_1}= V^{\lambda_2+2\lambda_3}\oplus V^{2\lambda_2}\oplus 
V^{\lambda_1+\lambda_3}. \]
Under the correspondence \eqref{E:ow} we have that the orbits corresponding to the weights $\lambda_1+\lambda_2+\lambda_3$, 
$2\lambda_1$, $\lambda_2+2\lambda_3$, $2\lambda_2$ and $\lambda_1+\lambda_3$ are respectively [(3,2,1)], [(1,1,0)], 
[(3,3,2)], [(2,2,0)] and [(3,1,0)]. Looking at \eqref{E:example2} from Example~\ref{E:or2}, we see again that the 
orbit computation gives exactly the same answer as the Pieri rule computation above.
\end{example}
\vskip10pt

\begin{example} \label{E:iteration}
Let $(\mu)=(3,2,1)$, $N=4$ and $k=3$. Let us compute the product $S_{(3,2,1)}h_2 h_1$ by iteration. In the previous example 
we found that $S_{(3,2,1)}h_2=S_{(3,3,2)}+S_{(2,2)}+S_{(2,1,1)}$, therefore we have
 \[ S_{(3,2,1)}h_2 h_1=S_{(3,3,2)}h_1+S_{(2,2)}h_1+S_{(2,1,1)}h_1.\]
The products $S_{(2,2)}h_1$ and $S_{(2,1,1)}h_1$ were already computed in \eqref{E:9} and \eqref{E:10} from 
Examples~\ref{Ex:9} and ~\ref{Ex:10}, so we only need to compute the product $S_{(3,3,2)}h_1$. The diagram below shows all 
diagrams $(\nu)\subseteq N\times k$ such that $(\nu)/(\mu)$ is a 1-row strip, giving the answer according to 
Theorem~\ref{T:hm2}
\vskip10pt
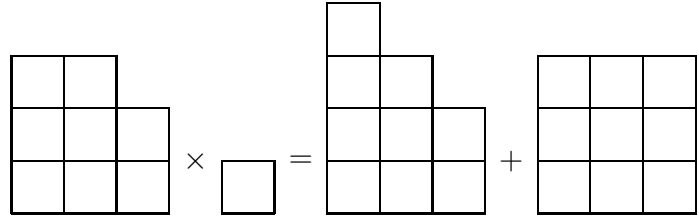
\begin{xy}
<0.7cm,0cm>:
(-4.38,0); (0,0) **@{ };
(0,0); (0,3) **@{-}, (3,0) **@{-};
(2,3); (0,3) **@{-}, (2,0) **@{-};
(3,2); (0,2) **@{-}, (3,0) **@{-};
(3,1); (0,1) **@{-};
(1,0); (1,3) **@{-};
(3.5,1), *{\times};
(4,0); (4,1) **@{-}, (5,0) **@{-};
(5,1); (4,1) **@{-}, (5,0) **@{-};
(5.5,1), *{=};
(6,0); (6,4) **@{-}, (9,0) **@{-};
(7,4); (6,4) **@{-}, (7,0) **@{-};
(8,3); (6,3) **@{-}, (8,0) **@{-}; 
(9,2); (6,2) **@{-}, (9,0) **@{-};
(9,1); (6,1) **@{-};
(9.5,1), *{+};
(10,0); (10,3) **@{-}, (13,0) **@{-};
(13,3); (10,3) **@{-}, (13,0) **@{-};
(13,2); (10,2) **@{-};
(13,1); (10,1) **@{-};
(11,0); (11,3) **@{-};
(12,3); (12,0) **@{-};
\end{xy}
\vskip20pt
Therefore $S_{(2,1,1)}h_1=S_{(3,3,2,1)}+S_{(3,3,3)}=S_{(2,2,1)}+S_{(3,3,3)}$. Now putting this together with \eqref{E:9} and 
\eqref{E:10}, we get 
\begin{equation} \label{E:13}
S_{(3,2,1)}h_2 h_1= 3 S_{(2,2,1)} + S_{(3,3,3)} + S_{(3,2)} + S_{(3,1,1)} + h_1.
\end{equation}

\end{example}

\section{Iterative fusion Pieri rules }

The next theorem illustrates how the iterative product of $S_\mu$ by several $h_m$'s, as in the previous example, can be done
without iteration if instead of using diagrams we use tableaux. We need a definition before stating the theorem.

\begin{definition}
Let $(\mu)\subseteq (N-1)\times k$, $(\nu)\in \Pi^{(N,k)}$ and $\epsilon =(\epsilon_1,\epsilon_2,...,\epsilon_r)$ be a 
sequence of nonnegative integers. A tableau of shape $(\nu)/(\mu)$ is called an $\mathbf{(N,k)}$-\textbf{cylindric tableau} 
if for each $1\leq p\leq \nu_N$ the entry in the top row and column $p$ of $(\nu)/(\mu)$ is strictly less than the entry in 
the bottom row and column $k+p$. We denote by $K_{(\nu) /(\mu),\epsilon}^{(N,k)}$ the number of $(N,k)$-cylindric tableaux of
shape $(\nu) /(\mu)$ and content $\epsilon$. The  number $K_{(\nu) /(\mu),\epsilon}^{(N,k)}$ is called the 
$\mathbf{(N,k)}$-\textbf{fusion skew Kostka number} and for $(\mu) =0$, the number $K_{(\nu) ,\epsilon}^{(N,k)}$ is called 
the $\mathbf{(N,k)}$-\textbf{fusion Kostka number}.
\end{definition}

\begin{example} \label{E:1}
Let $(\mu)=(3,2,1)$, $(\nu)=(4,2,2,1)$, $N=4$, $k=3$ and $\epsilon =(2,1)$. Below we draw all tableaux of shape 
$(\nu) /(\mu)$ and content $\epsilon$ and indicate whether each tableau is $(4,3)$-cylindric or not
\vskip10pt
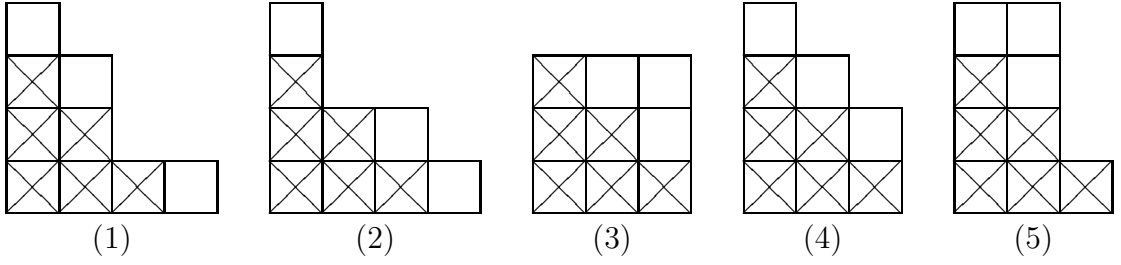
\begin{xy}
<0.7cm,0cm>:
(-2.63,0); (0,0) **@{ };
(0,0); (0,4) **@{-}, (4,0) **@{-};
(1,4); (0,4) **@{-}, (1,0) **@{-};
(4,1); (4,0) **@{-}, (0,1) **@{-};
(2,3); (0,3) **@{-}, (2,0) **@{-};
(3,1); (3,0) **@{-};
(2,2); (0,2) **@{-};
(6,0); (6,4) **@{-}, (10,0) **@{-};
(7,4); (6,4) **@{-}, (7,0) **@{-};
(10,1); (10,0) **@{-}, (6,1) **@{-};
(8,3); (6,3) **@{-}, (8,0) **@{-};
(9,1); (9,0) **@{-};
(8,2); (6,2) **@{-};
(12,0); (12,4) **@{-}, (16,0) **@{-};
(13,4); (12,4) **@{-}, (13,0) **@{-};
(16,1); (16,0) **@{-}, (12,1) **@{-};
(14,3); (12,3) **@{-}, (14,0) **@{-};
(15,1); (15,0) **@{-};
(14,2); (12,2) **@{-};
(0.5,3.5) *{1};
(1.5,2.5) *{1};
(3.5,0.5) *{2};
(6.5,3.5) *{1};
(7.5,2.5) *{2};
(9.5,0.5) *{1};
(12.5,3.5) *{2};
(13.5,2.5) *{1};
(15.5,0.5) *{1};
(2,-0.5) *{(3,4)-\text{cylindric}};
(8,-0.5) *{\text{not }(3,4)-\text{cylindric}};
(14,-0.5) *{\text{not }(3,4)-\text{cylindric}};
\end{xy}
\vskip10pt
\noindent
This shows that $K_{(4,2,2,1) /(3,2,1),(2,1)}^{(4,3)}=1$.
\end{example}

\begin{example}
This example shows the dependence on $k$ of the fusion skew Kostka number. Let $(\mu)$, $(\nu)$, $N$ and $\epsilon$ be as in 
the previous example and let $k=4$. One can see that all diagrams above are $(4,4)$-cylindric, therefore we have that
\newline $K_{(4,2,2,1) /(3,2,1),(2,1)}^{(4,4)}=3$.
\end{example}

\begin{example} \label{E:4}
Let $(\mu)=(3,2,1)$, $(\nu)=(3,3,2,1)$, $N=4$, $k=3$ and $\epsilon=(2,1)$. Below we show all tableaux of shape $(\nu)/(\mu)$ 
and content $\epsilon$, all of them are $(4,3)$-cylindric 
\vskip10pt
\begin{xy}
<0.7cm,0cm>:
(-5.38,0); (0,0) **@{ };
(0,0); (0,4) **@{-}, (3,0) **@{-};
(1,4); (0,4) **@{-}, (1,0) **@{-};
(3,2); (0,2) **@{-}, (3,0) **@{-};
(2,3); (0,3) **@{-}, (2,0) **@{-};
(3,1); (0,1) **@{-};
(4,0); (4,4) **@{-}, (7,0) **@{-};
(5,4); (4,4) **@{-}, (5,0) **@{-};
(7,2); (4,2) **@{-}, (7,0) **@{-};
(6,3); (4,3) **@{-}, (6,0) **@{-};
(7,1); (4,1) **@{-};
(8,0); (8,4) **@{-}, (11,0) **@{-};
(9,4); (8,4) **@{-}, (9,0) **@{-};
(11,2); (8,2) **@{-}, (11,0) **@{-};
(10,3); (8,3) **@{-}, (10,0) **@{-};
(11,1); (8,1) **@{-};
(0.5,3.5) *{1};
(1.5,2.5) *{1};
(2.5,1.5) *{2};
(4.5,3.5) *{1};
(5.5,2.5) *{2};
(6.5,1.5) *{1};
(8.5,3.5) *{2};
(9.5,2.5) *{1};
(10.5,1.5) *{1};
\end{xy}
\vskip10pt
\noindent
Therefore we have $K_{(3,3,2,1) /(3,2,1),(2,1)}^{(4,3)}=3$.
\end{example}

\begin{theorem} \label{T:prules}
 Let $\mu\in P_+^k$, $\epsilon=(\epsilon_1,\epsilon_2,\dots,\epsilon_r)$ a sequence of nonnegative integers and
$h_{\epsilon}=h_{\epsilon_1}\cdots h_{\epsilon_r}$. Then we have
\[ S_{(\mu)}h_{\epsilon}=\sum_{(\nu)\in\Pi^{(N,k)}}K_{\nu /\mu,\epsilon}^{(N,k)}S_{(\nu)} .\]
\end{theorem}

\begin{example}
Let us use this theorem to compute the product $S_{(3,2,1)}h_2 h_1$, with $N=4$ and $k=3$. We want to verify that our 
result matches the answer in Example \ref{E:iteration}. According to last theorem, we need to list all diagrams 
$(\nu)\in \Pi^{(4,3)}$ such that $(\nu)/(\mu)$ is a $\mathbf{(4,3)}$-cylindric tableau with content (2,1). Then it is clear 
that $(\nu)/(\mu)$ consists of 3 boxes. We show below all diagrams $(\nu)\in \Pi^{(4,3)}$ such that $(\nu)/(\mu)$ consists 
of 3 boxes:
\vskip10pt
\begin{xy}
<0.7cm,0cm>:
(-.38,0); (0,0) **@{ };
(0,0); (0,4) **@{-}, (4,0) **@{-},  (2,2) **@{-};
(1,4); (0,4) **@{-}, (1,0) **@{-};
(2,3); (0,3) **@{-}, (2,0) **@{-};
(4,1); (0,1) **@{-}, (4,0) **@{-};
(0,2); (2,2) **@{-}, (1,3) **@{-}, (2,0) **@{-};
(3,1); (3,0) **@{-}, (2,0) **@{-};
(0,1); (1,2) **@{-}, (1,0) **@{-};
(1,0); (2,1) **@{-};
(3,0); (0,3) **@{-};
(5,0); (9,0) **@{-}, (5,4) **@{-}, (7,2) **@{-};
(6,4); (5,4) **@{-}, (6,0) **@{-};
(8,2); (5,2) **@{-}, (8,0) **@{-};
(9,1); (5,1) **@{-}, (9,0) **@{-};
(7,0); (7,2) **@{-}, (8,1) **@{-};
(6,3); (5,3) **@{-};
(6,0); (7,1) **@{-};
(5,3); (8,0) **@{-};
(5,2); (6,3) **@{-}, (7,0) **@{-};
(5,1); (6,2) **@{-}, (6,0) **@{-};
(10,0); (10,3) **@{-}, (13,0) **@{-}, (12,2) **@{-};
(13,3); (10,3) **@{-}, (13,0) **@{-};
(10,2); (13,2) **@{-}, (12,0) **@{-}, (11,3) **@{-};
(10,1); (13,1) **@{-}, (11,2) **@{-}, (11,0) **@{-};
(12,0); (12,3) **@{-}, (13,1) **@{-};
(11,0); (11,3) **@{-}, (12,1) **@{-};
(14,0); (14,4) **@{-}, (17,0) **@{-};
(10,3); (13,0) **@{-};
(15,4); (14,4) **@{-}, (15,0) **@{-};
(16,3); (14,3) **@{-}, (16,0) **@{-};
(17,2); (14,2) **@{-}, (17,0) **@{-};
(17,1); (14,1) **@{-};
(14,3); (17,0) **@{-};
(14,2); (15,3) **@{-}, (16,0) **@{-};
(14,1); (15,2) **@{-}, (15,0) **@{-};
(14,0); (16,2) **@{-};
(15,0); (16,1) **@{-};
(16,0); (17,1) **@{-};
(20,0); (21,1) **@{-};
(18,3); (21,0) **@{-};
(18,2); (19,3) **@{-}, (20,0) **@{-};
(18,1); (19,2) **@{-}, (19,0) **@{-};
(19,0); (20,1) **@{-};
(18,0); (18,4) **@{-}, (20,2) **@{-}, (21,0) **@{-};
(20,4); (20,0) **@{-}, (18,4) **@{-};
(19,4); (19,0) **@{-};
(20,2); (18,2) **@{-};
(20,3); (18,3) **@{-};
(21,1); (18,1) **@{-}, (21,0) **@{-};
(2,-.5) *{(1)};
(7,-.5) *{(2)};
(11.5,-.5) *{(3)};
(15.5,-.5) *{(4)};
(19.5,-.5) *{(5)};
\end{xy}
\vskip20pt
Now we have to look at all possible fillings of these skew diagrams that make them into a $\mathbf{(4,3)}$-cylindric tableau.
Examples \ref{E:1} and \ref{E:4} above show the number of $(4,3)$-cylindric tableaux for the shapes (1) and (4), and it is 
easily seen that for the rest of the diagrams the Kostka skew number is equal to 1, so we have:
\begin{align}
 S_{(\mu)} h_2 h_1 & = S_{(4,2,2,1)}+ S_{(4,3,1,1)} +S_{(3,3,3)}+3S_{(3,3,2,1)}+S_{(3,2,2,2)} \notag\\
& =S_{(3,1,1)}+S_{(3,2)}+S_{(3,3,3)}+3S_{(2,2,1)}+S_{(1)} \notag\\
& =S_{(3,1,1)}+S_{(3,2)}+S_{(3,3,3)}+3S_{(2,2,1)}+h_1 .\notag
\end{align}
Just note that this matches \eqref{E:13} from Example \ref{E:iteration}.

\end{example}

\section{Rank-level duality } \label{S:ranklevel}

Several authors have stated the famous type A rank-level duality, see for instance \cite{Fu},  and many have tried to 
prove that $\mathcal{F}(A_{N-1},k)\cong \mathcal{F}(A_{k-1},N)$, which is actually not true (for instance 
$\mathcal{F}(A_1,3)$ is not isomorphic to $\mathcal{F}(A_2,2)$, in fact they do not even have the same dimensions which are 4 
and 6, respectively. However see \cite{Tu} or \cite{GW} for a proof). The isomorphism we are going to prove in this section 
was stated without proof in \cite{Fu}.

Before we start the proof of the isomorphism let's introduce some notation. We have 
\[ \mathcal{F}(A_{N-1},k)=\Lambda_N/I^{(N,k)}=\mathbb{Q}[x_1,x_2,...,x_N]^{S_N}/I^{(N,k)} \]
 is the algebra defined in Section~\ref{S:sympol} and 
\[ \mathcal{F}(A_{k-1},N)=\Lambda_k/I^{(k,N)}=\mathbb{Q}[y_1,y_2,...,y_k]^{S_k}/I^{(k,N)}. \]
  We keep the notation $S_{(\mu)}$ for Schur polynomials in $\Lambda_N/I^{(N,k)}$. To avoid confusion the Schur polynomials 
in $\Lambda_k/I^{(k,N)}$ will be denoted by $S'_{(\tilde\mu)}$. We have to use $ (\tilde\mu)$ to index Schur polynomials in 
$\Lambda_k$, because $\Lambda_N/I^{(N,k)}$ is indexed by partitions inside the rectangle $N\times k$ while 
$\Lambda_k/I^{(k,N)}$ is indexed by partitions inside the rectangle $k\times N$. We also keep the notation 
$h_m=h_m(x_1,x_2,...,x_N)$  and $e_s=e_s(x_1,x_2,...,x_N)$ for $h_m$ and $e_s$ in $\Lambda_N/I^{(N,k)}$, but we will use 
$h'_s$ to denote $h_s(y_1,y_2,...,y_k)$  and $e'_m$ for $e_m(y_1,y_2,...,y_k)$ in 
$\Lambda_k/I^{(k,N)}$, where $1\leq s\leq N$ and $1\leq m\leq k$.

Let $(\mu)\subseteq(N-1)\times k$ be a partition. Equation \eqref{E:em} from Theorem \ref{T:em} says that in
$\Lambda_N/I^{(N,k)}$ we have
\[ S_{(\mu)}e_m=\sum_{\substack{(\nu)\in\Pi^{(N,k)},\\(\nu) /(\mu) \text{ is an m-column strip}}}S_{(\nu)}. \]
The following remark will be the key in the proof of the rank-level duality isomorphism.

\begin{remark}
If we assume that $\mu_1<k$, that is, if $(\mu)\subseteq(N-1)\times (k-1)$, then the partitions $(\nu)$ on the right hand 
side of the equation above satisfy $(\nu)\subseteq N\times k$. 
\end{remark}
\begin{proof}
Assume that $\nu_1\geq k+1$, then since $\mu_1\leq k-1$ we have that $\nu_1-\mu_1\geq 2$. That means the bottom row of 
the Young diagram of $(\nu)$ has at least two more boxes than the bottom row of the Young diagram of $(\mu)$. This 
contradicts the fact that $(\nu)/(\mu)$ is an m-column strip.
\end{proof}

From this remark we get the following 
 
\begin{lemma} \label{L:nuevopieri}
If $(\mu)\subseteq(N-1)\times (k-1)$ then 
\begin{equation} \label{E:modifiedpieri}
S_{(\mu)}e_m=\sum_{\substack{(\nu)\subseteq N\times k,\\ (\nu) /(\mu) \text{ is an m-column strip}}}  S_{(\nu)}.
\end{equation}
\end{lemma}

Consider the partition $(1^m,0^{N-m})$. Since $e_m=S_{(1^m)}$ and the conjugate of $(1^m)$ is the partition 
$(\widetilde{1^m})=(m,0,0,...)$, applying conjugate to the equation \eqref{E:modifiedpieri} we have 
\begin{equation} \label{E:modifiedpieri2}
S'_{(\tilde\mu)}h'_m=\sum_{\substack{(\tilde\nu)\subseteq k\times N,\\ (\tilde\nu) /(\tilde\mu) \text{ is an m-row strip}}}
S'_{(\tilde\nu)}. 
\end{equation}
 This equation is exactly the result obtained in Theorem \ref{T:hm2}. This suggests that there might be a an algebra map from 
$\Lambda_N/I^{(N,k)}$ to $\Lambda_k/I^{(k,N)}$ defined by ``conjugation" 
\begin{equation}
 S_{(\mu)}  \longrightarrow S'_{(\tilde\mu)}. \notag
\end{equation}
This map fails to be an algebra homomorphism since $e_N=S_{(1^N)}=1\in\Lambda_N/I^{(N,k)}$  is mapped to 
$S'_{(\widetilde{1^N})}=S'_{(N)}=h'_N$ but $h'_N\ne 1$ in  $\Lambda_k/I^{(k,N)}$. But if we modify the ideal $I^{(N,k)}$ by 
adding the relation $h_k=1$ and modify the ideal $I^{(k,N)}$  by adding the relation $h'_N=1$, we will see how the above map 
becomes an isomorphism between the two new quotient algebras.

Consider the ideal $\bar{I}^{(N,k)}$ of $\Lambda_N=\mathbb{Q}[x_1,x_2,...,x_N]^{S_N}$ defined by
\[ \bar{I}^{(N,k)}=\langle S_{(1^N)} -1,S_{(\mu)}, h_k-1 \mid \mu_1 -\mu_N =k+1 \rangle \]
and the ideal $\bar{I}^{(k,N)}$ of $\Lambda_k=\mathbb{Q}[y_1,y_2,...,y_k]^{S_k}$  defined by
\[ \bar{I}^{(k,N)}=\langle S'_{(\tilde{1^k})} -1,S'_{(\tilde\mu)}, h'_N-1 \mid \tilde{\mu}_1 -\tilde{\mu}_k =N+1 \rangle. \]

Then we get a homomorphism 
 \begin{align} \label{E:ranklevelduality}
 \Lambda_N/\bar{I}^{(N,k)} &\longmapsto \Lambda_k/\bar{I}^{(k,N)}  \\
S_{(\mu)} & \longmapsto S'_{(\tilde\mu)}. \notag
\end{align}

The map sends $h_m=S_{(m)}$ to $S'_{(\tilde m)}=S'_{(1^m)}=e'_m$ for $0\leq m\leq k-1$ and $e_m=S_{(1^m)}$ to $S'_{(m)}=h'_m$
for $0\leq m\leq N-1$. Moreover, by Theorem \ref{T:hm2} and Lemma \ref{L:nuevopieri} we get that the product $S_{(\mu)}h_m$
in $\Lambda_N/\bar{I}^{(N,k)}$ matches the product $S'_{(\tilde\mu)}e'_m$ in $\Lambda_k/\bar{I}^{(k,N)}$ 
(See \eqref{E:modifiedpieri} and \eqref{E:modifiedpieri2}). We also get that the product $S_{(\tilde\mu)}e_m$ in 
$\Lambda_N/\bar{I}^{(N,k)}$ matches the product $S'_{(\mu)}h'_m$ in $\Lambda_k/\bar{I}^{(k,N)}$. Therefore, the homomorphism 
above sends generators to generators and the corresponding multiplication by each generator in both rings agrees, therefore 
we have an isomorphism between both algebras. This isomorphism is known as the Type A rank-level duality.

Using the notation for fusion algebras, the Dictionary \eqref{E:dict} relates $h_k=S_{(k,0\dots)}$ with $V^{k\lambda_1}$ in
$\mathcal{F}(A_{N-1},k)$ and $h'_N=S_{(N,0\dots)}$ with $V^{N\lambda_1}$ in $\mathcal{F}(A_{k-1},N)$, so we get an 
isomorphism between the algebras 
\[ \mathcal{F}'(A_{N-1},k)\cong \mathcal{F}'(A_{k-1},N), \]
where
\[ \mathcal{F}'(A_{N-1},k)=\mathcal{F}(A_{N-1},k)/\langle V^{k\lambda_1}-V^0\rangle\]
and
\[ \mathcal{F}'(A_{k-1},N)=\mathcal{F}(A_{k-1},N)/\langle V^{N\lambda_1}-V^0\rangle. \]

\section{Simple currents }

In the particular case when $m=k$ of Equation \eqref{E:p2} in Theorem \ref{T:hm2} we get that for 
$(\mu)\subseteq (N-1)\times k$
\[ S_{(\mu)}h_k=\sum_{\substack{(\nu)\subseteq{N\times k},\\ (\nu) /(\mu)\text{ is a $k$-row strip}}}S_{(\nu)}. \]
Looking at the right hand side of this equation, we get that there is a unique partition $(\nu)\subseteq N\times k$ so that
$(\nu) /(\mu)$ is a $k$-row strip, namely if $(\mu)=(\mu_1,\mu_2,...,\mu_{N-1})$ then $(\nu)=(k,\mu_1,\mu_2,...,\mu_{N-1})$.
Therefore the equation above becomes
\begin{equation} \label{E:simplecurrent}
 S_{(\mu)}h_k=S_{(k,\mu_1,\mu_2,...,\mu_{N-1})}= S_{(k-\mu_{N-1},\mu_1-\mu_{N-1},\mu_2-\mu_{N-1},...,
\mu_{N-2}-\mu_{N-1})}.
\end{equation}

The basis elements with this property are called simple currents. We establish the precise definition of a 
simple current next.
\begin{definition}
Let $B=\{x_a\mid a\in A\}$ be a basis for a fusion algebra $\mathcal{F}$. We say that $x_a\in B$ is a simple current if
for any element $x_b\in B$ the product $x_ax_b=x_c\in B$ for some $c\in A$.
\end{definition}

From equation \eqref{E:simplecurrent} we get that $h_k$ and therefore $h_k^t$ for $t\geq 1$ are simple currents in 
$\Lambda_N/I^{(N,k)}$. In fact, multiplying $S_{(\mu)}$ iteratively by $h_k$ we get 
\begin{equation}
S_{(\mu)}h_k^t=S_{(\mu')}
\end{equation}
where 
\begin{align} \label{E:productbysimplecurrent}
(\mu')=&(k-\mu_{N-t}+\mu_{N-(t-1)},k-\mu_{N-t}+\mu_{N-(t-2)},...,k-\mu_{N-t}+\mu_{N-1},k-\mu_{N-t},
\notag\\&\mu_1-\mu_{N-t},...,\mu_{N-(t+1)}-\mu_{N-t}).
\end{align}
In particular, if $(\mu)=(k)$, since $S_{(\mu)}=h_k$ we get $h_k^{t+1}=S_{(k^{(t+1)})}$ and therefore $h_k^N=S_{(k^N)}=1$ 
in $\Lambda_N/I^{(N,k)}$, since the partition $(k^N)$ consists of $k$ columns of height $N$.

 It is also known that the set of simple currents of the basis for $\Lambda_N/I^{(N,k)}$ indexed  by partitions whose 
Young diagram are inside the $(N-1)\times k$ rectangle, is the set $\{h_k^t\mid t\geq  0\}$ and that this set actually spans 
a sub-algebra of $\Lambda_N/I^{(N,k)}$. Moreover, since $h_k^N=1$ in $\Lambda_N/I^{(N,k)}$, we get that the set of simple 
currents forms a cyclic group under multiplication isomorphic to $\mathbb{Z}_N$ and the sub-algebra they span is the group 
algebra $\mathbb{Q}[\mathbb{Z}_N]$.

So we can see that the algebra $\Lambda_N/\bar{I}^{(N,k)}$ can be obtained from the fusion algebra $\Lambda_N/I^{(N,k)}$ by 
identifying simple currents with $1$, in other words by making the quotient by the ideal generated by $h_k-1$,
 $\langle h_k-1 \rangle$.


We conclude the chapter by recalling a known characterization of the fusion coefficients $N_{\lambda,\mu}^{(k)\nu}$
when $\lambda$ is a multiple of the first fundamental weight $\lambda_1$ of type $A_{N-1}$, and an extension of it using the
theory of simple currents. The precise characterization is as follows (see \cite{Tu}):

\begin{equation} \label{E:tudose}
N_{\mu,m\lambda_1}^{(k)\nu}=
\begin{cases}
1, &\text{ if there exists $(\bar\nu)\subseteq N\times k$ such that $(\bar\nu)/(\mu)$ is an m-row strip}\\
& \text{ and $(\bar\nu)\thicksim(\nu)$,}\\
0, &\text{ otherwise, }
\end{cases}
\end{equation}
where $m\leq k$ and $\thicksim$ is the equivalence relation defined in Equation \eqref{E:er} and $(\mu)$ and $(\nu)$ are the 
partitions corresponding to $\mu$ and $\nu$.

Let $(\mu')$ be the partition from Equation \eqref{E:productbysimplecurrent} and let $(\lambda)$ be the partition such that
$S_{(\lambda)}=h_k^t h_m $ for $m\leq k$ and $0\leq t\leq N-1$. Then we get
\begin{align}
 S_{(\mu)}S_{(\lambda)}= &  S_{(\mu)} h_k^t h_m \notag\\
= & S_{(\mu')} h_m \notag
\end{align}
Therefore we get
\begin{equation} \label{E:graciasadios}
N_{\mu,\lambda}^{(k)\nu}=N_{\mu',m\lambda_1}^{(k)\nu}.
\end{equation}


\chapter{Connection between Young diagrams and $S_k$-orbits of $\mathbb{Z}_N^k$} \label{S:ydo}

\section{Main theorem and results }

Let $\lambda=a_1\lambda_1+a_2\lambda_2+...+a_{N-1}\lambda_{N-1}$ be a dominant integral weight for $\mathfrak{g}=sl_N$ of 
type $A_{N-1}$. The orbit in standard form associated to $\lambda$ is given by \eqref{E:orb} and the partition is given by 
\eqref{E:tab}. These two points of view are related as the following lemma shows. 
\begin{lemma} \label{L:conj}
 The orbit in standard form \eqref{E:orb} is a partition whose conjugate is given by \eqref{E:tab}.
\end{lemma}
\begin{proof}
Consider the diagram

\vskip10pt

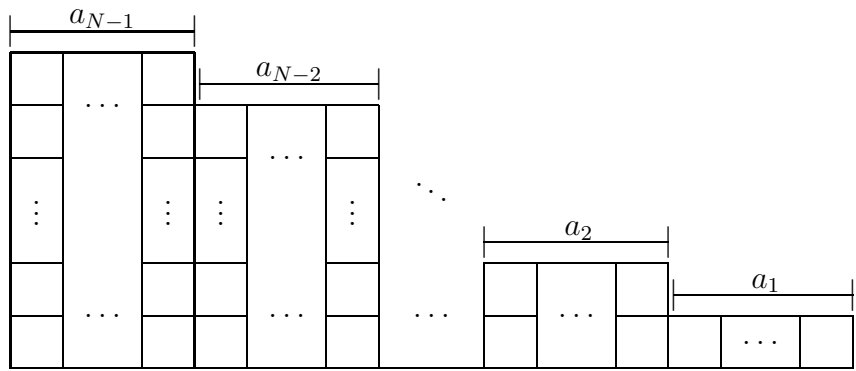
 \begin{xy}
<0.7cm,0cm>:
(-2.88,0); (0,0) **{ };
(0,0); (0,6) **@{-}, (1,0) **@{-};
(1,6); (0,6) **@{-}, (1,0) **@{-};
(0,1); (1,1) **@{-},
(0,2); (1,2) **@{-},
(0,4); (1,4) **@{-},
(0,5); (1,5) **@{-},
(1,0); (2.5,0) **@{-}, 
(1,6); (2.5,6) **@{-}, 
(1.75,5) *{\dots};
(1.75,1) *{\dots};
(0.5,3.1) *{\vdots};
(3,3.1) *{\vdots};
(4,3.1) *{\vdots};
(6.5,3.1) *{\vdots};
(2.5,0); (2.5,6) **@{-}, (4.5,0) **@{-};
(3.5,6); (2.5,6) **@{-}, (3.5,0) **@{-};
(2.5,1); (4.5,1) **@{-},
(2.5,2); (4.5,2) **@{-},
(2.5,4); (4.5,4) **@{-},
(2.5,5); (4.5,5) **@{-},
(4.5,0); (4.5,5) **@{-},
(5.25,1) *{\dots};
(5.25,4) *{\dots};
(4.5,0); (6,0) **@{-},
(4.5,5); (6,5) **@{-},
(6,0); (6,5) **@{-}, (7,0) **@{-};
(7,5); (6,5) **@{-}, (7,0) **@{-};
(6,1); (7,1) **@{-},
(6,2); (7,2) **@{-},
(6,4); (7,4) **@{-},
(8,1) *{\dots}; 
(7,0); (9,0) **@{-},
(8,3.5) *{\ddots};
(9,0); (9,2) **@{-}, (10,0) **@{-};
(10,2); (9,2) **@{-}, (10,0) **@{-};
(9,1); (10,1) **@{-},
(10,0); (11.5,0) **@{-},
(10,2); (11.5,2) **@{-},
(10.75,1) *{\dots};
(11.5,0); (11.5,2) **@{-}, (13.5,0) **@{-};
(12.5,2); (11.5,2) **@{-}, (12.5,0) **@{-};
(11.5,1); (13.5,1) **@{-},
(13.5,0); (13.5,1) **@{-},
(13.5,0); (15,0) **@{-},
(13.5,1); (15,1) **@{-},
(14.25,0.5) *{\dots};
(15,0); (15,1) **@{-}, (16,0) **@{-};
(16,1); (15,1) **@{-}, (16,0) **@{-};
(0,6.4) *{|}; 
(3.5,6.4) *{|};
(0,6.4); (3.5,6.4) **@{-};
(1.75,6.62) *{a_{N-1}};
(3.6,5.4) *{|}; 
(7,5.4) *{|};
(3.6,5.4); (7,5.4) **@{-};
(5.3,5.62) *{a_{N-2}};
(9,2.4) *{|}; 
(12.5,2.4) *{|};
(9,2.4); (12.5,2.4) **@{-};
(10.8,2.62) *{a_2};
(12.6,1.4) *{|}; 
(16,1.4) *{|};
(12.6,1.4); (16,1.4) **@{-};
(14.35,1.62) *{a_1};
\end{xy}
\vskip10pt
\noindent
where the label $a_i$ on top equals the number of columns in that block, and the subscript i, equals the height of each 
column in that block. This means that, there are $a_{N-1}$ columns of height $N-1$, $a_{N-2}$ columns of height $N-2$,...,
$a_2$ columns of height 2 and $a_1$ columns of height 1.
Reading the diagram from left to right we get the partition \eqref{E:orb}, and reading the diagram from bottom to top, we  
get the following: The first row has length $a_{N-1}+a_{N-2}+...+a_1=\sum_{j=1}^{N-1} a_j$. The second row has length 
$a_{N-1}+a_{N-2}+...+a_2=\sum_{j=2}^{N-1} a_j$, and the last row has length $a_{N-1}$, that is, the diagram above corresponds 
to the partition \eqref{E:tab}. Thus \eqref{E:orb} and \eqref{E:tab} are conjugates of each other.
\end{proof}

\begin{remark} \label{R:tuple}
This lemma says that from a Young diagram $(\mu)\subseteq (N-1)\times k$, we can associate with $(\mu)$ a $k$-tuple 
whose $j$th entry is given by the height of the $j$-th column of the diagram. This $k$-tuple is actually the representative 
in standard form of the orbit in $\mathbb{Z}_N^k$ whose corresponding partition is $(\mu)$. This $k$-tuple will be denoted 
by $\hat\mu$.
\end{remark}

\begin{remark} \label{R:nu'}
If $(\nu)=(\nu_1,...,\nu_N)\subseteq N\times k$ we can also associate with $(\nu)$ a $k$-tuple in $\mathbb{Z}_N^k$, by 
letting the $j$-th entry of the $k$-tuple be the height of the $j$-th column of $(\nu)$ modulo $N$. This $k$-tuple is a 
representative from the orbit $[\nu]$ associated to the partition $(\nu)$, but is not necessarily in standard form. 
\textbf{We denote this} $\mathbf{k}$-\textbf{tuple by} $\mathbf{\bar\nu}$. The representative in standard form for the orbit 
$[\nu]$ can be obtained directly as follows: Let $(\nu')=(\nu_1-\nu_N,...,\nu_{N-1}-\nu_N)$, then 
$(\nu')\subseteq (N-1)\times k$, so the $k$-tuple whose $j$-th entry is the $j$-th column of the Young diagram for
 $(\nu')$ is the representative in standard form of the orbit $[\nu]$ associated to the partition $(\nu)$. Therefore
$\nu'\in [\nu]$. We also get that $(\nu)\thicksim(\bar\nu)\thicksim(\nu')$ where $\thicksim$ is the equivalence relation 
defined in \eqref{E:er}.
\end{remark}

Seeing the standard form of an orbit in $\mathbb{Z}_N ^k$ as a partition will help us prove why, in some special cases, the 
orbit product \eqref{E:orbprod} matches the Pieri rules. The next theorem illustrates some of these cases.

\begin{theorem} \label{T:main}
Let $\lambda=m\lambda_1$ be a multiple of the first fundamental weight for $A_{N-1}$, $m\leq k$,
$\mu=a_1\lambda_1+...+a_{N-1}\lambda_{N-1}$ any other weight of level $k$ and $[\lambda]$ and $[\mu]$ their 
corresponding orbits in $\mathbb{Z}_N ^k$. Then $N_{\mu,\lambda}^{(k)\nu}=M_{[\mu],[\lambda]}^{(k)[\nu]}$ for any weight
$\nu$ of level $k$.
In other words, the product of the orbits $[\mu]\times [\lambda]$ matches Pieri rules for the multiplication $S_{(\mu)} h_m$,
 \eqref{E:p2}. 
\end{theorem}
\begin{proof}
The first step of the proof is to prove the following claim

\noindent
Claim: $N_{\mu,\lambda}^{(k)\nu}=1$ if and only if $M_{[\mu],[\lambda]}^{(k)[\nu]}=1$

The theorem will follow from this claim, Corollary ~\ref{C:char} and Equation \eqref{E:tudose}.

\noindent
Proof of the claim: Let's assume that $N_{\mu,\lambda}^{(k)\nu}=1$. From Equation \eqref{E:tudose}, we know that there is a 
partition $(\nu')\subseteq N\times k$ so that $(\nu')/(\mu)$ is an m-row strip, and $(\nu')\thicksim(\nu)$, then from 
Remark \ref{R:nu'} we get that $[\nu]=[\nu']$.
Since  $(\nu')\subseteq N\times k$ and $(\nu')/(\mu)$ is an m-row strip, $(\nu')$ can be obtained from $(\mu)$ by adding m 
boxes with no two in the same column, hence there exists $y\in [(1^m,0^{k-m})]$ so that 
\[ \hat{\nu}+y=\nu'. \]
Now since $y\in [(1^m,0^{k-m})]$ and $(\nu')\in[\nu']=[\nu]$ we have that $M_{[\mu],[\lambda]}^{(k)[\nu]}\ne 0$, then from
Corollary \ref{C:char} we get $M_{[\mu],[\lambda]}^{(k)[\nu]}=1$.

In the other direction, let's assume that $M_{[\mu],[\lambda]}^{(k)[\nu]}=1$. Then from Corollary ~\ref{C:char} we have that
$\hat\nu=\left((N-1)^{a_{N-1}-m_{N-1}+m_{N-2}},...,1^{a_1-m_1+m_0},0^{a_0-m_0+m_{N-1}}\right)$, for some integers $m_0$, 
$m_1$,...,$m_{N-1}$ such that $\sum_{i=0}^{N-1}m_i=m$ and $0\leq m_i\leq a_i$ for $0\leq i\leq N-1$.

Consider the $k$-tuple $\bar\nu\in\mathbb{Z}_N^k$ defined in Remark \ref{R:nu'}, 
\[ \bar\nu=\left(0^{m_{N-1}},(N-1)^{a_{N-1}-m_{N-1}+m_{N-2}},...,1^{a_1-m_1+m_0},0^{a_0-m_0}\right) \]
so $\bar\nu\in[\nu]$. Define the Young diagram $(\nu')$ whose first $m_{N-1}$ columns have height $N$ and for $j>m_{N-1}$, 
the $j^{th}$ column has height equal to the $j^{th}$ entry of $\bar\nu$. 

Now since the $k$-tuples $\hat\mu$ and $\nu'$ satisfy the equation
\[  \hat\mu +(1^{m_{N-1}},0^{a_{N-1}+{m_N-1}},...,1^{m_0},0^{a_0-m_0})=\nu' \quad\text{where}\quad \sum_{i=0}^{N-1}m_i=m, \]
then we have that $(\nu')/(\mu)$ is an m-row strip and since $(\nu')\thicksim(\bar\nu)\thicksim(\nu)$, then we get from
\eqref{E:tudose} that $N_{\mu,\lambda}^{(k)\nu}=1$.
\end{proof}

Note that in general multiplication of two orbits gives a linear combination of orbits. The orbit $[k\lambda_1]$ is special 
since every product of this orbit with any other gives a single orbit as an answer. We encode this result in the  following 
lemma.

\begin{lemma} \label{L:single}
 Let $[\mu]$ be any $S_k$-orbit of $\mathbb{Z}_N^k$. The product $[k\lambda_1]\times[\mu]$ equals a single orbit.
\end{lemma}
\begin{proof}
Since the $k$-tuple $(1^k)$ is fixed by the action of $S_k$, the orbit of $[k\lambda_1]=[(1^k)]$ consists of a single 
$k$-tuple. The result easily follows from this fact.
\end{proof}

Although the product defined on orbits is not associative, as shown in example~\ref{E:nonassos}, we can prove associativity
in the special case when one of the orbits involved is $[(t^k)]$, $1\leq t\leq N-1$. We state the result in the following 
lemma.

\begin{lemma} \label{L:ass}
Let $[a]$ and $[b]$ be $S_k$-orbits of $\mathbb{Z}_N^k$, and let $0\leq t\leq N-1$. Then we have:
\begin{equation} \label{E:assoc}
 \left( [a]\times [(t^k)]\right)\times [b] =[a]\times\left([(t^k)]\times[b]\right). 
\end{equation}
\end{lemma}
\begin{proof}
From lemma \ref{L:single}, the products $[a]\times [(t^k)]$ and $[(t^k)]\times[b]$ each consists of a single orbit. If 
$a=(a_1,...,a_k)$ is a representative from the orbit $[a]$, then a representative of the orbit $[a]\times [(t^k)]$ is given
by $(a_1+t,...,a_k+t)$ (where every entry is considered modulo $N$). Let us denote by $[a+t]$ the orbit of 
$[a]\times [(t^k)]$, and since the same applies to the orbit $[(t^k)]\times[b]$, we denote this orbit by $[b+t]$, that is, 
\begin{equation} \label{E:a+t}
[a]\times [(t^k)] = [a+t] \quad\text{and}\quad [(t^k)]\times[b] = [b+t].
\end{equation}
Then from \eqref{E:assoc} and \eqref{E:a+t} all we have to do is prove the following equation:
\[ [a+t]\times [b] = [a]\times [b+t]. \]
This is equivalent to proving that for any orbit $[c]$ of $\mathbb{Z}_N^k$ we have that:
\begin{equation} \label{E:M}
 M_{[a+t],[b]}^{(k)[c]}=M_{[a],[b+t]}^{(k)[c]} .
\end{equation}
To prove this equality, first observe that there is a one-to-one map given by: 
\begin{align} \label{E:T}
  T([a+t],[b],[c]) & \longrightarrow T([a],[b+t],[c]) \\
  (x,y,z) & \longmapsto (x-t,y+t,c) \notag
\end{align}
where $x-t$ is the $k$-tuple  obtained from $x$ by subtracting $t$ from every entry of $x$ and $y+t$ is defined by adding $t$
to every entry of $y$, that is, $x-t=x-(t^k)$ and $y+t=y+(t^k)$. Let $\sigma\in S_k$. Since the action of $S_k$ on $k$-tuples
is linear, we have that  $\sigma(x-t)=\sigma (x-(t^k))=\sigma x-\sigma(t^k)$ and 
$\sigma(y+t)=\sigma (y+(t^k))=\sigma y+\sigma(t^k)$. Therefore the map \eqref{E:T} is invariant under the action of $S_k$. 
Thus the equality \eqref{E:M} follows. 
\end{proof}

These lemmas are useful to find other orbits whose fusion product can be computed by the orbits method, and those are given 
below.

\begin{theorem} \label{T:extend}
For $m\leq k$ let $[m\lambda_1]$ and $[k\lambda_1]$ be the orbits in $\mathbb{Z}_N^k$ corresponding to the weights 
$m\lambda_1$ and $k\lambda_1$, respectively. For $0\leq t\leq N-1$ let $\lambda$ be the weight whose corresponding orbit in 
$\mathbb{Z}_N^k$ is given by $[\lambda]=[m\lambda_1]\times[k\lambda_1]^t$, and let $\mu$ be any other weight on level k. Then
for every weight $\nu$ on level $k$ we have
\[ N_{\lambda,\mu}^{(k)\nu}=M_{[\lambda],[\mu]}^{(k)[\nu]}. \]
\end{theorem}
\begin{proof}
Since $[k\lambda_1]^t=[(t^k)]$, for an orbit $[\mu]$, we have:
\begin{equation} \label{E:associativity}
 [\lambda]\times[\mu]=([m\lambda_1]\times[k\lambda_1]^t)\times[ \mu]=[m\lambda_1]\times([k\lambda_1]^t\times[ \mu]). 
\end{equation}
From Lemma~\ref{L:single} the orbits product $[k\lambda_1]^t\times[ \mu]=[\mu+t]$ gives a single orbit. Let $\mu'$ be the 
weight so that $[\mu']=[\mu+t]$. Then from \eqref{E:associativity} we get
\[ [\lambda]\times[\mu]=[m\lambda_1]\times[\mu']. \]
Therefore $M_{[\lambda],[\mu]}^{(k)[\nu]}=M_{[m\lambda_1],[\mu']}^{(k)[\nu]}$, for all weights $\nu$.

From Theorem~\ref{T:main}, we have $M_{[m\lambda_1],[\mu']}^{(k)[\nu]}=N_{m\lambda_1,\mu'}^{(k)\nu}$ for all weights $\nu$. 
This combined with \eqref{E:graciasadios} gives us 
\[ N_{\lambda,\mu}^{(k)\nu}=M_{[m\lambda_1],[\mu']}^{(k)[\nu]}.\]
Therefore we get
\[ N_{\lambda,\mu}^{(k)\nu}=M_{[\lambda],[\mu]}^{(k)[\nu]}. \]
\end{proof}

We get the following characterization of these special level $k$ fusion coefficients.

\begin{theorem}
Let $\mu$ and $\nu$ be arbitrary weights on level $k$, for $0\leq t\leq N-1$ and $m\leq k$ let $\lambda$ be the weight 
associated to the orbit $[m\lambda_1]\times[k\lambda_1]^t$, let $\mu'$ be the weight whose corresponding orbit is 
$[\mu']=[\mu]\times[k\lambda_1]^t$, and assume that $\hat\mu= ((N-1)^{a_{N-1}},\dots,1^{a_1},0^{a_0})$. Then we have 
\[ M_{[\lambda],[\mu]}^{(k)[\nu]}=
\begin{cases} \label{E:characterization}
1,  & \text{if } \tilde\nu=((N-1)^{b_{N-1}-m_{N-1}+m_{N-2}},\dots,1^{b_1-m_1+m_0},0^{b_0-m_0+m_{N-1}}),\\
  & \text{where }  m_0,\dots,m_{N-1} \text{ are integers such that } \displaystyle{\sum_{i=0}^{N-1}} m_i=m, \\
 &   0\leq m_i\leq b_i \text{ and } b_i=\begin{cases} a_{i-t} &\text{ if } i\geq t \\ a_{N-(t-i)} &\text{ if } i<t 
\end{cases}\\
0, & \text{otherwise}
\end{cases} \]
\end{theorem}

As a summary of the results in this section, more precisely, from Theorem~\ref{T:extend}, we get an algorithm to compute 
fusion coefficients via orbits for Young diagrams of the form

\vskip10pt

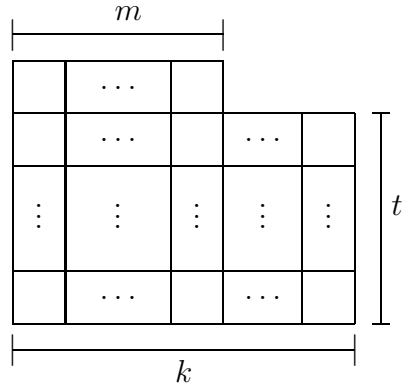
\begin{xy}
<0.7cm,0cm>:
(-7.18,0); (0,0) **{ };
(0,0); (6.5,0) **@{-}, (0,5) **@{-};
(6.5,4); (6.5,0) **@{-}, (0,4) **@{-};
(4,5); (4,0) **@{-}, (0,5) **@{-};
(1,0); (1,5) **@{-};
(3,0); (3,5) **@{-};
(5.5,0); (5.5,4) **@{-};
(0,1); (6.5,1) **@{-};
(0,3); (6.5,3) **@{-};
(2,.5) *{\dots};
(4.75,.5) *{\dots};
(2,3.5) *{\dots};
(2,4.5) *{\dots};
(4.75,3.5) *{\dots};
(.5,2.2) *{\vdots};
(2,2.2) *{\vdots};
(4.75,2.2) *{\vdots};
(6,2.2) *{\vdots};
(3.5,2.2) *{\vdots};
(0,-0.5) *{|};
(6.5,-0.5) *{|};
(0,-0.5); (6.5,-0.5) **@{-};
(3.25,-0.9) *{k};
(0,5.5) *{|};
(4,5.5) *{|};
(0,5.5); (4,5.5) **@{-};
(2.2,5.9) *{m};
(7,0); (7,4) **@{-};
(7.05,0) *{\_};
(7.05,4) *{\_};
(6.9,0) *{\_};
(6.9,4) *{\_};
(7.3,2.25) *{t};
\end{xy}
\vskip10pt
\noindent
whose corresponding weights are
\begin{equation} \label{E:weight}
 \lambda=
\begin{cases}
m\lambda_1,                       & t=0, \\
(k-m)\lambda_t + m\lambda_{t+1},  & 1\leq t\leq N-2, \\
(k-m)\lambda_{N-1},               & t=N-1 ,
\end{cases} 
\end{equation}
for $0\leq m\leq k$. The corresponding orbits are
\begin{equation} \label{E:orbits}
 [\lambda]=
\begin{cases}
[((t+1)^m,t^{k-m})],              &  0\leq t\leq N-2, \\
[((N-1)^{k-m},0^m)],  & t=N-1.
\end{cases} 
\end{equation}

\begin{example}
For $A_3$ level 3 ($N=4$ and $k=3$), let $\lambda=2\lambda_3+\lambda_2$ and let $\mu=2\lambda_1+\lambda_2$. The weight 
$\lambda$ has the form \eqref{E:weight} with $t=2$ and $m=2$. Theorem~\ref{T:extend} says that we can compute the fusion 
product $V^{\lambda}\otimes_3 V^{\mu}$ via orbits of $\mathbb{Z}_4^3$. The orbits associated to $\lambda$ and $\mu$ are
$[\lambda]=[(3,3,2)]$ and $[\mu]=[(2,1,1)]$. To compute their product, we have to remove the redundant equation from the list
\begin{align}
(3,3,2)+& (2,1,1)= (1,0,3), \notag\\
(3,3,2)+& (1,2,1)= (0,1,3), \notag\\
(3,3,2)+& (1,1,2)= (0,0,0). \notag
\end{align}
The second equation on this list is redundant, therefore we get 
\[ [\lambda]\times[\mu]= [(3,1,0)] +[(0,0,0)]. \]
The weights associated to $[(3,1,0)]$ and $[(0,0,0)]$ are $\nu_1=\lambda_1+\lambda_3$ and $\nu_2=0$, so we get
\[ V^{\lambda}\otimes_3 V^{\mu}=V^{\lambda_1+\lambda_3}\oplus V^{0}. \]
We can verify this result using the characterization \eqref{E:characterization}. Note that 
\[ a_3=0,\quad a_2=1,\quad a_1=2,\quad a_0=0, \]
so
\[ b_3=a_1=2,\quad b_2=a_0=0,\quad b_1=a_3=0,\quad b_0=a_2=1. \]
There are two different sets of integers $\{m_0,m_1,m_2,m_3\} $ satisfying the conditions of Theorem 
\eqref{E:characterization}, namely:
\[ \{m_0=1,m_1=0,m_2=0,m_3=1\}\quad\text{and}\quad \{m'_0=0,m'_1=0,m'_2=0,m'_3=2\}. \]
For the set on the left, $\{m_0=1,m_1=0,m_2=0,m_3=1\}$, we get
\[ \hat\nu_1=(3^{b_3-m_3+m_2},2^{b_2-m_2+m_1},1^{b_1-m_1+m_0},0^{b_0-m_0+m_3})=(3,1,0),\]
and for the set $\{m'_0=0,m'_1=0,m'_2=0,m'_3=2\}$, we get
\[ \hat\nu_2=(3^{b_3-m'_3+m'_2},2^{b_2-m'_2+m'_1},1^{b_1-m'_1+m'_0},0^{b_0-m'_0+m'_3})=(0,0,0),\]
which are the representatives in standard form of the orbits we got, and Theorem \eqref{E:characterization} predicts 
coefficient zero for the rest of the orbits, which matches the result found in this example.
\end{example}

\section{What goes wrong with orbits }

The orbits product $[a]\times [b]=\sum_{[c]\in\mathcal{O}}M_{[a],[b]}^{(k)[c]}[c]$ does not match Pieri rules for all cases,
as we show in this section.
\begin{example}
Consider the weight $\lambda=\lambda_1+\lambda_2$ of $A_2$ on level 3. The partition associated to $\lambda$ is 
$(\lambda)=(2,1,0)$. By applying iterative Pieri rules (Theorem~\ref{T:prules}) to compute 
$S_{(\lambda )} S_{(\lambda)}$ and the fact that $S_{(\lambda)}=h_1h_2-h_3$, we get:
\begin{align} 
S_{(\lambda )} S_{(\lambda)} & =S_{3}+ S_{(3,3)}+2S_{(\lambda)}+ S_{(0)} \notag\\
& =h_3+ S_{(3,3)}+ 2S_{(\lambda)}+1.\notag 
\end{align}
From the Dictionary \eqref{E:dict}, this translates into the fusion product
\begin{equation} \label{E:exa} 
V^{\lambda }\otimes_3 V^{\lambda}=V^{3\lambda_1}\oplus V^{3\lambda_2}\oplus2V^{\lambda_1+\lambda_2}\oplus V^{0}. 
\end{equation}
This result does not match the result gotten from the orbits point of view. To see this, first note that the orbit 
corresponding to $\lambda=\lambda_1+\lambda_2$ in $\mathbb{Z}_3^3$ is the orbit of $[(2,1,0)]$, and by a simple calculation 
of orbits we get that:
\[ [(2,1,0)]\times [(2,1,0)]=[(1,1,1)]+ [(2,2,2)]+ 3[(2,1,0)]+ [(0,0,0)]. \]
The orbits in this product represent respectively the weights $3\lambda_1$, $3\lambda_2$, $\lambda_1+\lambda_2$ and $0$, 
which are exactly the same weights we got in \eqref{E:exa}, however multiplicities do not all match.
\end{example}

\begin{example}
Consider the weights $\lambda=\lambda_1+\lambda_2$ and $\mu=2\lambda_2$ of $A_3$ on level 3. The partitions associated to 
$\lambda$ and $\mu$ are $(\lambda)=(2,1,0)$ and $(\mu)=(2,2)$. Using theorem~\ref{T:prules} to compute the fusion product 
between the Schur polynomials $S_{(\lambda)}$ and $S_{(\mu)}$ and the fact that $S_{(\mu)}=h_2^2-h_1h_3$, we get:
\[ S_{(\lambda )} S_{(\mu)}=S_{(3,2,2)}+S_{(3,3,1)}+S_{(2,1)}+ S_{(1,1,1)}. \]
Using the correspondence between orbits and partitions defined in the last section of Chapter \ref{S:pre} and given more
explicitly below in \eqref{E:opart}  we get
\[ V^{\lambda}\otimes_3 V^{\mu}=V^{\lambda_1+2\lambda_3}\oplus V^{2\lambda_2+\lambda_3}\oplus V^{\lambda_1+\lambda_2}\oplus 
  V^{\lambda_3}. \]
Now performing the orbit computation $[\lambda ]=[(2,1,0)]$ times $[\mu ]=[(2,2,0)]$ in $\mathbb{Z}_4^3$, we get:
\[ [(2,1,0)]\times [(2,1,0)]=[(3,2,2)]\oplus [(2,1,0)]\oplus [(3,0,0)]. \]
The orbits in this product represent, respectively, the weights $2\lambda_2+\lambda_3$, $\lambda_1+\lambda_2$ and $\lambda_3$. 
Observe that the orbit corresponding to the weight $\lambda_1+2\lambda_3$ is missing in the orbit computation.
\end{example}

\section{How the orbits method can be fixed }

Note that the Schur polynomial corresponding to weight $m\lambda_1$ is $S_{(m,0,...,0)}=h_m$, that is, weights that are 
multiples of the first fundamental weight $\lambda_1$ correspond to the homogeneous symmetric polynomials. According to the 
Jacobi-Trudi determinant \eqref{E:schur}, these polynomials, generate the quotient ring 
$\Lambda_N/I^{(N,k)}=\mathbb{Q}[x_1,x_2,...,x_N]^{S_N}/I^{(N,k)}$. Therefore for the rest of the weights of level $k$ whose 
Young diagram is not of the form \eqref{E:weight}, we can compute the fusion products via orbits by alternating iteration 
using the information from the Jacobi-Trudi determinant. In the following example we show how this process can be done.

\begin{example} \label{E:extention}
Consider the weight $\lambda=\lambda_1+\lambda_2$ of $A_2$ on level 3, note that this weight is not of the form 
\eqref{E:weight}. If we want to compute the fusion product of $V^{\lambda}$ with itself on level 3 
($V^{\lambda}\otimes_3 V^{\lambda}$), we may use the fact that $S_{\lambda}=h_1h_2-h_3$ and compute:
\begin{align}
 V^{\lambda}\otimes_3 V^{\lambda} & =V^{\lambda}\otimes_3\left( V^{\lambda_1}\otimes_3 V^{2\lambda_1}-V^{3\lambda_1} \right) 
\notag \\ &=V^{\lambda}\otimes_3\left( V^{\lambda_1}\otimes_3 V^{2\lambda_1}\right)-V^{\lambda}\otimes_3V^{3\lambda_1} \notag
\\ &= \left( V^{\lambda}\otimes_3 V^{\lambda_1}\right)\otimes_3 V^{2\lambda_1}-V^{\lambda}\otimes_3V^{3\lambda_1} \notag
\end{align}

Now from Theorem~\ref{T:main} this can be done iteratively via orbits.
\end{example}

To ``fix" the orbits method, we construct an algebra $\mathcal{A}(\mathcal{O}(N,k))$ over $\mathbb{Q}$ whose basis is 
$\mathcal{O}(N,k)$ where the addition is formal addition and we define a new product of orbits $[a]\cdot[b]$ by using the 
information from the Jacobi-Trudi determinant \eqref{E:schur}. We define this ``new product" as follows.

Let $\hat a$, $\hat b$ be the standard form of two orbits $[a],[b]\in \mathcal{O}$. Let $(\lambda)=\tilde a$ and 
$(\mu)=\tilde b=(\tilde b_1,\dots,\tilde b_{N-1})$ be the conjugates of $\hat a$, $\hat b$ considered as partitions. The 
product $S_{(\lambda)}S_{(\mu)}$ can be computed using Pieri rules iteratively with the help of the Jacobi-Trudi determinant, 
\eqref{E:schur}. Let $H=(h_{\tilde b_i -i+j})_{i,j=1}^{N-1}=(h_{i,j})_{i,j=1}^{N-1}$ then 
$S_{(\mu)}=\det(H)\in\Lambda_N/I^{(N,k)}$, so
\begin{equation} \label{E:altite}
  S_{(\lambda)}S_{(\mu)} =S_{(\lambda)}\det(H). 
\end{equation}
By Theorem \ref{T:main} we know that each individual product $S_{(\lambda)}h_{\tilde b_i -i+j}$ can be computed via orbits by 
performing the product $[\lambda]\times [(1^{\tilde b_i -i+j},0^{k-\tilde b_i +i-j})]$, provided that $\tilde b_i -i+j\leq k$.
Therefore, we get that we can multiply any two Schur polynomials via orbits and by the results of this chapter, product 
\eqref{E:altite} is equivalent to compute the following orbit product and then translating the results via the Dictionary
\ref{E:dict1}
\begin{equation} \label{E:fixing}
 [a]\cdot[b]=
\begin{cases}
[a]\times[b], &\text{ if } [b] \text{ is an orbit}\\&\text{of the form \eqref{E:orbits}, }\\
\displaystyle{\sum_{\sigma\in S_{N-1}}}\text{sign}(\sigma)((([a]\cdot [h_{1,\sigma 1}])\cdot [h_{2,\sigma 2}])\cdots 
[h_{N-1,\sigma(N-1)}])  &\text{ otherwise,}
\end{cases} 
\end{equation}
where $[h_{i,j}]=[h_{\tilde b _i -i+j}]$, $\tilde b=(\tilde b_1,\dots,\tilde b_{N-1})$ is the conjugate of $\hat b$ 
considered as partition, and 
\[ [h_{\tilde b _i -i+j}]= \begin{cases}
[(1^{\tilde b_i -i+j},0^{k-\tilde b_i -i+j})] &\text{ if } 0\leq \tilde b_i -i+j\leq k \\
0 &\text{ otherwise,} 
\end{cases}\]
is of the form \eqref{E:orbits}. We extend this product linearly to all $\mathcal{A}(\mathcal{O}(N,k))$.

From Lemma \ref{L:conj} we know that the the standard form of an $S_k$-orbit is a partition whose conjugate is a partition 
contained in the $(N-1)\times k$ rectangle. The explicit correspondence between these two objects is given by
\begin{equation} \label{E:opart}
 [((N-1)^{a_N-1},\dots,1^{a_1},0^{a_0})]\longleftrightarrow \left(\sum_{i=1}^{N-1}a_i,\sum_{i=2}^{N-1}a_i,\dots,a_1\right). 
\end{equation}
This gives a dictionary between the algebras $\mathcal{A}(\mathcal{O}(N,k))$ and $\Lambda_N/I^{(N,k)}$ as follows
\begin{align} \label{E:dict1}
\mathcal{A}(\mathcal{O}(N,k))& \longleftrightarrow \Lambda_N/I^{(N,k)} \\
[\lambda] & \longleftrightarrow S_{(\lambda)},\notag\\
[h_m]=[(1^m,0^{k-m})]& \longleftrightarrow h_m,\notag\\
[\lambda]\cdot[\mu]& \longleftrightarrow S_{(\lambda)}S_{(\mu)},\notag
\end{align}
where the left hand side is indexed by orbits and the right hand side is indexed by partitions and the correspondence 
\eqref{E:opart} allows us to move from left to right.

According to the results in this chapter, we get the product \eqref{E:fixing} after translating to Schur functions using 
Dictionary \ref{E:dict1}, matching the associative product $S_{(\lambda)}S_{(\mu)}$. Therefore \eqref{E:fixing} is an 
associative product. Thus \eqref{E:fixing} gives a structure of associative algebra to $\mathcal{A}(\mathcal{O}(N,k))$ and 
this algebra is isomorphic to $\Lambda_N/I^{(N,k)}\cong\mathcal{F}(A_{N-1},k)$.

\begin{example}
Continuing with Example \ref{E:extention} we want to compute the fusion product $V^{\lambda}\otimes_3 V^{\lambda}$ where 
$\lambda=\lambda_1+\lambda_2$ of $A_2$ on level 3. 

The orbit corresponding to $\lambda$ is $[\lambda]=[(2,1,0)]$ and the conjugate of $\hat{\lambda}$ is $\tilde b=(2,1)$. Hence
Equation \eqref{E:fixing} gives
\[  [(2,1,0)]\cdot[(2,1,0)]=[(2,1,0)]\cdot \begin{vmatrix}
    [h_{\tilde b_1}] & [h_{\tilde b_1 +1}]  \\
[h_{\tilde b_2 -1}] & [h_{\tilde b_2 }]  \end{vmatrix}. \]

Since $\tilde b_1=2$ and $\tilde b_2=1$, we get 
\begin{align}
[h_{\tilde b_1}] &=[(1,1,0)], \notag\\
[h_{\tilde b_1+1}] &=[(1,1,1)], \notag\\
[h_{\tilde b_2}] &=[(1,0,0)], \notag\\
[h_{\tilde b_-1}] &=[(0,0,0)]. \notag
\end{align}

Therefore we get 
\begin{align} [(2,1,0)]\cdot[(2,1,0)]&=[(2,1,0)]\cdot\begin{vmatrix}
    [(1,1,0)] & [(1,1,1)] \\
   [(0,0,0)]   & [(1,0,0)]
 \end{vmatrix} \notag\\
&=[(2,1,0)]\cdot([(1,1,0)] \cdot [(1,0,0)] ) - [(2,1,0)]\cdot[(1,1,1)] \notag. 
\end{align}
From Example \ref{E:or1} we know that 
\[ [(2,1,0)]\cdot[(1,1,0)]= [(2,1,0)]\times [(1,1,0)]=  [(2,2,1)]+[(1,1,0)]+[(2,0,0)] \]
and it is readily checked that 
\begin{align}
[(2,2,1)]& \cdot [(1,0,0)]=[(2,2,1)] \times [(1,0,0)]= [(2,2,2)] + [(2,1,0)],\notag\\
[(1,1,0)]& \cdot [(1,0,0)]=[(1,1,0)] \times [(1,0,0)]= [(1,1,1)] + [(2,1,0)],\notag\\
[(2,0,0)]& \cdot [(1,0,0)]=[(2,0,0)] \times [(1,0,0)]= [(2,1,0)] + [(0,0,0)]\notag
\end{align}
and 
\[ [(2,1,0)]\cdot[(1,1,1)]=[(2,1,0)]\times[(1,1,1)]=[(2,1,0)]. \]
Therefore we get 
\[ [(2,1,0)]\cdot[(2,1,0)]=[(2,2,2)]+[(1,1,1)]+2[(2,1,0)]+[(0,0,0)], \]
which translates into the fusion product 
\[ V^{\lambda_1+\lambda_2}\otimes_3V^{\lambda_1+\lambda_2}= V^{3\lambda_2}\oplus V^{3\lambda_1} \oplus 
2V^{\lambda_1+\lambda_2} \oplus V^0. \]
\end{example}

\section{Rank-level duality from the orbits point of view }

The sub-algebra of simple currents of $\mathcal{F}(A_{N-1},k)$ is the sub-algebra with basis $\{ h_k^t \mid 0\leq t\leq N-1\}$.
The orbit corresponding to the partition $(k,0,...,0)$ is the orbit $[(1^k)]$ and by multiplying this orbit by itself 
iteratively we get that $[(1^k)]^t=[(t^k)]$. Therefore the orbits corresponding to the simple currents of 
$\mathcal{F}(A_{N-1},k)$ have the form $[(t^k)]$. Now, let $[b]=[(t^k)]$ be the orbit of a simple current and 
$[a]=[(m_1,...,m_k)]\in\mathcal{O}(N,k)$. Then the product $[a]\cdot [b]$ can be easily described as 
\[ [a]\cdot [b]=[(m_1+t,...,m_k+t)], \]
where every entry is considered modulo $N$. By abuse of notation we are going to denote the orbit on the right hand side by 
$[a+t]$.

 This product gives an action of $\mathbb{Z}_N$ on $\mathcal{O}(N,k)$. This action breaks the set $\mathcal{O}(N,k)$ into 
``orbits", that we call SC-orbits, to avoid confusion with $S_k$-orbits of $\mathbb{Z}_N^k$ and as an abbreviation for simple 
current orbits. 

Note that for $0\leq m\leq k$ the product of any element in the SC-orbit of $[(1^m,0^{k-m})]$ by any SC-orbit is independent of the choice of 
representative. In fact, let $0\leq t,t_1\leq N-1$ and $[a]\in \mathcal{O}(N,k)$, then $[a+t_1]$ is in the SC-orbit of 
$[a]$ 
and $[((t+1)^m,t^{k-m})]$ is in the SC-orbit of $[(1^m,0^{k-m})]$ and is one of the orbits of the form \eqref{E:orbits}, 
therefore
\begin{align} \label{E:simplecurrentproduct}
[a+t_1]\cdot [((t+1)^m,t^{k-m})]=& [a]\cdot [(t_1^k)]\cdot [(1^m,0^{k-m})] \cdot [(t^k)] \\
=& [a]\cdot [(1^m,0^{k-m})]\cdot [((t+t_1)^k)]. \notag
\end{align}
This defines a product on the set of SC-orbits of $\mathcal{O}(N,k)$, and it is equivalent to the product in the algebra
$\mathcal{A}'(\mathcal{O}(N,k))=\mathcal{A}(\mathcal{O}(N,k))/SC$, where SC is the ideal of 
$\mathcal{A}(\mathcal{O}(N,k))$ generated by $[(t^k)]-[(0^k)]$. This algebra is isomorphic to $\mathcal{F}'(A_{N-1},k)$.

To see the rank-level duality from the orbits point of view, first observe that any SC-orbit of $\mathcal{O}(N,k)$ has 
a representative of the form $[((N-1)^{a_{N-1}},...,0^{a_0})]$ where $a_0>0$, since for any orbit 
$[(m_1,...,m_k)]\in\mathcal{O}(N,k)$ we can add $(t^k)$ enough times to make one of the entries zero.

\begin{example}
Consider the set $\mathcal{O}(4,3)$. The SC-orbit of the element $[(2,2,1)]$ is the set $\{[(2,2,1)],[(3,3,2)],[(3,0,0)],
[(1,1,0)]\} $. Note that in fact this SC-orbit contains two elements with at least one zero entry.
\end{example}

Let $\mathcal{O}'(N,k)=\{[a]=[(a_1,\dots,a_k)]\in\mathcal{O}(N,k)\mid a_i=0\text{ for some } i, 1\leq i\leq k\}$, the set of 
$S_k$-orbits of $\mathbb{Z}_k^N$ with at least one zero entry. We have a one-to-one correspondence given by
\begin{align}
\mathcal{O}'(N,k) &\longrightarrow \mathcal{O}'(k,N)\notag\\
[((N-1)^{a_{(N-1)}},...,1^{a_1},0^{a_0})]&\longmapsto \left[\left( \sum_{i=1}^{N-1}a_i,\sum_{i=2}^{N-1}a_i,..., a_{(N-1)},0
\right)\right],
\notag
\end{align}
where $a_0>0$ and $\sum_{i=0}^{N-1}a_i=k$.

This map sends SC-orbits of $\mathcal{O}(N,k)$ to SC-orbits of $\mathcal{O}(k,N)$, and since the map is defined by 
conjugation of the orbits considered as partitions, it follows that the map is just a different way of defining the 
rank-level duality map defined in Section \ref{S:ranklevel} and Equation \eqref{E:ranklevelduality}.

\chapter{Application to tensor product }

\section{Classical Pieri rules }

Irreducible modules $V^{\lambda}$ for the finite dimensional Lie algebra $\mathfrak{g}=sl_N$, 
are indexed by weights 
$\lambda\in P^+$ and form a basis of an associative algebra over $\mathbb{Q}$, called the tensor algebra and denoted by
$T(A_{N-1})$, with associative product of $V^{\lambda}$ and $V^{\mu}$ defined by the tensor product of the two modules 
$V^{\lambda}\otimes V^{\mu}$. This product is a well defined product in $T(A_{N-1})$ since given modules $V^{\mu}$ and 
$V^{\lambda}$ for $sl_N$, with $\mu$ and $\lambda\in P^+$, the tensor product $V^{\mu}\otimes V^{\lambda} $ is completely 
reducible and it is given by
\begin{equation} \label{E:tensor}
 V^{\lambda}\otimes V^{\mu}=\sum_{\nu\in P^+}Mult_{\lambda,\mu}(\nu)V^{\nu} 
\end{equation}
where $Mult_{\lambda,\mu}(\nu)$ is given by the Racah-Speiser formula \eqref{E:racah}.

The Littlewood-Richardson coefficients are the structure constants for the algebra of symmetric polynomials $\Lambda_N$
with respect to the basis given by Schur polynomials indexed by partitions of length at most $N$. If $(\mu)$ and $(\lambda)$ 
are partitions of length at most $N$ and $S_{(\mu)}$ and $S_{(\lambda)}$ are the Schur polynomials as defined in 
\eqref{E:schur} (or in \eqref{E:schur1}), then we have
\begin{equation} \label{E:tensorschur}
  S_{(\lambda)}S_{(\mu)}=\sum_{(\nu) \text{ is a partition of length at most } N}c_{(\lambda),(\mu)}^{(\nu)}S_{(\nu)} .
\end{equation}
 
There is a relation between these two algebras. In fact, it is well known that there is an isomorphism between the tensor 
algebra $T(A_{N-1})$ and the quotient algebra $\Lambda_N/\langle S_{(1^N)}-1\rangle$, given by 
\begin{align} \label{E:tensordict}
 \Lambda_N/\langle S_{(1^N)}-1\rangle & \longleftrightarrow T(A_{N-1}) \\
 S_{(\mu)} & \longleftrightarrow V^{\mu} \notag\\
S_{(\lambda)}S_{(\mu)}=\sum_{(\nu)}c_{(\lambda),(\mu)}^{(\nu)}S_{(\nu)} & \longleftrightarrow 
V^{\mu}\otimes V^{\lambda}=\bigoplus_{\nu}Mult_{\lambda,\mu}(\nu)V^{\nu}, \notag
\end{align}
where the left hand side is indexed by partitions $(\mu)=(\mu_1,...,\mu_N)$ of length at most $N$, and the right hand side 
is indexed by weights $\mu=\sum_{j=1}^{N-1}(\mu_j-\mu_{j+1})\lambda_j \in P^+$ where $P^+$ is the set of dominant 
integral weights defined in \eqref{E:domintwts}. The correspondence \eqref{E:wp} allows one to move from one side to
the other and $c_{(\lambda),(\mu)}^{(\nu)}=Mult_{\lambda,\mu}(\nu)$.

Particular cases of the product \eqref{E:tensorschur} can be computed using Pieri rules. The following results are known as
classical Pieri rules.

\begin{theorem} (Classical Pieri rules for multiplication by $h_m$) \label{T:tprhm}

Let $(\mu)$ be a partition of length at most $N$ and $m$ a positive integer, then in $\Lambda_N/\langle S_{(1^N)}-1\rangle$
we have
\[  S_{(\mu)}h_m=\sum_{(\nu) /(\mu)\text{ is an m-row strip  } }S_{(\nu)}. \]
\end{theorem}

\begin{theorem} (Classical Pieri rules for multiplication by $e_m$) \label{T:tprem}

Let $(\mu)$ be a partition of length at most $N$ and $0\leq m\leq N$ an integer, then in 
$\Lambda_N/\langle S_{(1^N)}-1\rangle$ we have
\[  S_{(\mu)}e_m=\sum_{(\nu) /(\mu)\text{ is an m-column strip  } }S_{(\nu)}. \]
\end{theorem}

It is also known that the Littlewood-Richardson coefficients $c_{(\mu),(\lambda)}^{(\nu)}$ are at most one when 
$(\lambda)=(m)$ is a partition with a single part and when $(\lambda)=(1^m)$ is a partition with only ones. The 
characterization is given by 
\begin{equation} \label{E:htensortudose}
c_{(\mu),(m)}^{(\nu)}=
\begin{cases}
1, &\text{ if there exists a partition $(\bar\nu)$ of length at most N,  }\\
& \text{ such that $(\bar\nu)/(\mu)$ is an m-row strip and $(\bar\nu)\thicksim(\nu)$,}\\
0, &\text{ otherwise, }
\end{cases}
\end{equation}
where $\thicksim$ is the equivalence relation defined in Equation \eqref{E:er}, and 
\begin{equation} \label{E:etensortudose}
c_{(\mu),(1^m)}^{(\nu)}=
\begin{cases}
1, &\text{ if there exists a partition $(\bar\nu)$ of length at most N,  }\\
& \text{ such that $(\bar\nu)/(\mu)$ is an m-column strip and $(\bar\nu)\thicksim(\nu)$,}\\
0, &\text{ otherwise. }
\end{cases}
\end{equation}

We are going to use Equation \eqref{E:htensortudose} to prove that we can compute the weight decomposition of  
$V^{\mu}\otimes V^{m\lambda_1}$ via orbits, but since there is no level $k$ condition for tensor 
products, the sum $\sum_{j=1}^{N-1} a_j$ and $m$ can be arbitrarily large. Therefore, our 1-1 correspondence between weights
and orbits, given by \eqref{E:ow}, has to be modified, and we also have to modify our group $G$. 

\section{$S_\infty$-orbits of $\bigoplus_{i=1}^{\infty}\mathbb{Z}_N$ }

Let $\bar {G}=\bar {G}_N=\bigoplus_{i=1}^{\infty}\mathbb{Z}_N=\{ (x_1,x_2,\dots)\mid x_i\in\mathbb{Z}_N, x_j=0 \text{ for } 
j>>0\}$ and let $S_\infty$, the group of permutations of $\{1,2,\dots\}$ which fix all but finitely many numbers, act on 
$\bar {G}_N$ by permuting the coordinates. For $x\in\bar {G}_N$, define the element in standard form 
$S_\infty$-equivalent to $x$ to be 
\[ \hat x=( (N-1)^{a_{N-1}},...,1^{a_1},0,...)\]
where $i$ occurs $a_i$ times in $x$.

There is a 1-1 correspondence between dominant weights of $A_{N-1}$ and $S_\infty$-orbits of $\bar G$ as follows:
\begin{equation} \label{E:tensorwp}
a_1\lambda_1+...+a_{N-1}\lambda_{N-1} \in P^+\leftrightarrow \left[\left( (N-1)^{a_{N-1}},...,1^{a_1},0,...\right)\right],
\end{equation}
where the bracket on the right means the $S_\infty$-orbit of the element in standard form. 

The correspondence between partitions of length at most $N$ and weights of $A_{N-1}$ given by \eqref{E:wp} and the 
correspondence \eqref{E:tensorwp}, gives us a correspondence between partitions of length at most $N$ and orbits of $\bar G$.
In this chapter will use the same notation as before for weights and partitions, that is, if $\mu$ is a weight, then 
$(\mu)$ denotes its corresponding partition, but for this chapter, $[\mu]$ will denote the $S_\infty$-orbit in $\bar G$ 
corresponding to the weight $\mu$ and the representative in standard form of $[\mu]$ in $\bar G$ will be denoted by 
$\hat \mu$. The set of $S_\infty$-orbits of $\bar G$ will be denoted $\mathcal{O}_\infty=\mathcal{O}_\infty(N)$.

Now we extend the definition of the set $T([a],[b],[c])$ to $S_\infty$-orbits $[a]$, $[b]$ and $[c]$ of $\bar G$:
\[ T([a],[b],[c])=\{(x,y,z)\in [a]\times [b]\times [c] \mid x+y=z \}, \]
 and define
$M_{[a],[b]}^{[c]}$ as the number of $S_\infty$-orbits of $T([a],[b],[c])$. 

The number $M_{[a],[b]}^{[c]}$ can be also used to define a product of $S_\infty$-orbits of $\bar{G}_N$ as follows
\[ [a]\times [b]=\sum_{[c]\in\mathcal{O}_\infty} M_{[a],[b]}^{[c]}[c].\]

The procedure to compute the number $M_{[a],[b]}^{[c]}$ is similar to the one we used to compute the number
$M_{[a],[b]}^{(k)[c]}$. In fact, since elements in $\bar G$ have a finite number of non-zero entries, we can follow exactly
the same procedure to compute $M_{[a],[b]}^{(k)[c]}$, with $k>>0$, as follows.

\begin{definition}
Let $[a]$, $[b]\in\mathcal{O_{\infty}}$, and assume that 
$\hat{a}=( (N-1)^{a_{N-1}},...,1^{a_1},0,...)$ $\hat{b}=( (N-1)^{b_{N-1}},...,1^{b_1},0,...)$, and set $k$ to be any integer
such that $k\geq\sum_{i=1}^{N-1}(a_i+b_i)$. Let $\{y_1,\dots,y_t\}\subset [b]$ be the set of orbits of $\hat{b}$ obtained by
letting $S_k$ act on the first $k$ entries of $\hat{b}$, and for $1\leq i\leq t$ set:
\begin{equation} \label{E:inftylist}
 z_i=\hat a +y_i .
\end{equation}
We say that the equation $z_j=\hat a +y_j$ in the list \eqref{E:list} is \textbf{redundant}, if for some $i<j$ and 
$\sigma\in S_k$ we have
\[ \sigma\hat{a} =\hat{a}, \qquad\sigma y_j=y_i\quad\text{and}\quad \sigma z_j=z_i ,\]
that is, if the triples $(\hat{a} ,y_i,z_i)$ and $(\hat{a} ,y_j,z_j)$ are in the same $S_k$-orbit of $T([a],[b],[z_i])$.
\end{definition}

Then it is easily seen that
\begin{equation} \label{E:inftyorbits}
 M_{[a],[b]}^{[c]}\text{ equals the number of non-redundant equations of the form \eqref{E:inftylist} }. 
\end{equation}

From this, it follows that $ M_{[a],[b]}^{[c]}= M_{[a],[b]}^{(k)[c]}$, for $k>>0$.
\begin{example}  \label{E:tensorcomp}
Let $a=(2,1,0,\dots)$ and $b=(1,1,0,\dots)$ be elements in $\bar{G}_3$.  Then we have that 
$\hat a =(2,1,0,\dots)$ and $\hat b= (1,1,0,\dots)$. So the least $k$ we can choose to compute the orbit product 
$[a]\times[b]$ is $k=4$. By letting $S_4$ act on the first four entries of $\hat b$ we get the subset of $[b]$
\begin{align}
 \{ \hat b= &(1,1,0,\dots), (1,0,1,0,\dots), (1,0,0,1,0,\dots),\notag\\
& (0,1,1,0,\dots), (0,1,0,1,0,\dots), (0,0,1,1,0,\dots)\}.\notag
\end{align}

Now we have to remove all redundancies from the list:
\begin{align}
& (2,1,0,\dots) + (1,1,0,\dots) = (0,2,0,\dots) ,\notag\\
&(2,1,0,\dots)+ (1,0,1,0,\dots) = (0,1,1,0,\dots) ,\notag\\
&(2,1,0,\dots)+ (1,0,0,1,0,\dots) = (0,1,0,1,0,\dots) ,\notag\\
&(2,1,0,\dots)+ (0,1,1,0,\dots) = (2,2,1,0,\dots) ,\notag  \\
&(2,1,0,\dots)+ (0,1,0,1,0,\dots) = (2,2,0,1,0,\dots) ,\notag  \\
&(2,1,0,\dots)+ (0,0,1,1,0,\dots) = (2,1,1,1,0,\dots). \notag
\end{align}
The third and the fifth equations are redundant in that list. So we get:
\begin{align} 
 [(2,1,0,\dots)]\times [(1,1,0,\dots)]= & [(2,2,1,0,\dots)]+[(1,1,0,\dots)]\\& +[(2,0,0,\dots)] +[(2,1,1,1,0,\dots)] .
\notag
\end{align}

Later we will see how this corresponds to a tensor product decomposition for $A_2$ of 
$V^{\lambda_1+\lambda_2}\otimes V^{2\lambda_1}$.
\end{example}

\begin{example}
Now let $a=(2,2,1,0,\dots)$ and $b=(1,0,0,\dots)$ be elements in $\bar{G}_3$.  Then we have that 
$\hat a =(2,2,1)$ and $\hat b= (1,0,0,\dots)$. We can take $k=4$ to compute $[a]\times [b]$. By letting $S_4$ act on the 
first four entries of  $\hat b$ we get the subset of $[b]$
\[ \{ (1,0,0,\dots), (0,1,0,\dots), (0,0,1,0,\dots),\newline(0,0,0,1,0,\dots)\}. \]

Then we remove all redundancies from the list:
\begin{align}
  &(2,2,1,0,\dots)+ (1,0,0,\dots) = (0,2,1,0,\dots) , \notag\\
           &(2,2,1,0,\dots)+ (0,1,0,\dots) = (2,0,1,0,\dots) ,\notag\\
           &(2,2,1,0,\dots)+ (0,0,1,0,\dots) = (2,2,2,0,\dots),\notag  \\
&(2,2,1,0,\dots)+ (0,0,0,1,0,\dots) = (2,2,1,1,0,\dots).\notag  
\end{align}
We can see that the second equation is redundant, so we can remove it from the list and we get
\[ [(2,2,1,0,\dots)]\times [(1,0,\dots)]=  [(2,2,2,0,\dots)]+[(2,1,0,\dots)]+[(2,2,1,1,0,\dots)]. \]
\end{example}


\begin{theorem} \label{T:tensororbitproduct}
Let $[a],[c]\in\mathcal{O}$ and $[b]=[(1^m,0,\dots)]$. Suppose that $M_{[a],[b]}^{[c]}\ne 0$ then 
$M_{[a],[b]}^{[c]}=1$.
\end{theorem}
\begin{proof}
Since for $k>>0$ we have $M_{[a],[b]}^{[c]}=M_{[a],[b]}^{(k)[c]}$, the proof follows from Theorem \ref{T:orbitproduct}.
\end{proof}

We also get the following characterization of orbits products similar to Corollary \ref{C:char} in Chapter \ref{C:oua}.
\begin{corollary} \label{C:tensorchar}
Let $[a]\in\mathcal{O}$, assume that $\hat a=\left( (N-1)^{a_{N-1}},...,1^{a_1},0,\dots\right)$ and let 
$[b]=[(1^m,0,\dots)]$. Then for $[c]\in\mathcal{O}$ we have
\[ M_{[a],[b]}^{[c]}=
\begin{cases} 
1, &\text{if $\hat{c} =\left((N-1)^{a_{N-1}-m_{N-1}+m_{N-2}},...,1^{a_1-m_1+m_0},0,\dots\right)$, }\\ 
&\text{ for some integers $m_0$, $m_1$,...,$m_{N-1}$ such that $\displaystyle\sum_{i=0}^{N-1}m_i=m$}\\
&  \text{ and  } \quad 0\leq m_i\leq a_i\quad\text{, for  } \quad 1\leq i\leq N-1, \\
0, &\text{otherwise.}
\end{cases}
\]
\end{corollary}

It is well known that for $k>>0$ the coefficients $c_{(\mu),(\lambda)}^{(\nu)}=N_{(\mu),(\lambda)}^{(k)(\nu)}$, see 
\cite{Fe}, Theorem 6.1. Then we have the following result.

\begin{theorem} \label{T:tensormain}
For $1\leq m\in \mathbb{Z}$ let $\lambda=m\lambda_1$ be a multiple of the first fundamental weight for $A_{N-1}$, 
$\mu=a_1\lambda_1+...+a_{N-1}\lambda_{N-1}$ any other dominant weight and $[\lambda]$ and $[\mu]$ their corresponding orbits 
in $\bar G$. Then $c_{(\mu),(\lambda)}^{(\nu)}=M_{[\mu],[\lambda]}^{[\nu]}$ for any dominant weight $\nu$.
In other words, the product of the $S_\infty$-orbits $[\lambda]\times[\mu]$ matches classical Pieri rules for the 
multiplication $S_{(\mu)}h_m$, Theorem \ref{T:tprhm}. 
\end{theorem}
\begin{proof}
For $k>>0$ we have $M_{[\mu],[\lambda]}^{[\nu]}=M_{[\mu],[\lambda]}^{(k)[\nu]}$ and 
$c_{(\mu),(\lambda)}^{(\nu)}=N_{(\mu),(\lambda)}^{(k)(\nu)}$. Now from Theorem \ref{T:main} we have that 
$M_{[\mu],[\lambda]}^{(k)[\nu]}=N_{(\mu),(\lambda)}^{(k)(\nu)}$ from where the theorem follows. 
\end{proof}

\begin{example}
Let $(\mu)=(2,1)$ and $N=3$. Using Theorem~\ref{T:tprhm} to compute $S_{(\mu)}h_2$ we get
\vskip20pt
\begin{xy}
<0.7cm,0cm>:
(-.21,0); (0,0) **{ };
(0,0); (0,2) **@{-}, (2,0) **@{-};
(2,1); (0,1) **@{-}, (2,0) **@{-},
(1,2); (1,0) **@{-}, (0,2) **@{-};
(2.5,1) *{\times};
(3,0); (3,1) **@{-}, (5,0) **@{-};
(5,1); (3,1) **@{-}, (5,0) **@{-};
(4,1); (4,0) **@{-};
(5.5,1) *{=}; (8.5,1) *{+}; (12.5,1) *{+}; (16.5,1) *{+};
(6,0); (6,3) **@{-}, (8,0) **@{-};
(8,2); (6,2) **@{-}, (8,0) **@{-};
(7,3); (6,3) **@{-}, (7,0) **@{-};
(8,1); (6,1) **@{-};
(9,0); (9,3) **@{-}, (12,0) **@{-};
(12,1); (9,1) **@{-}, (12,0) **@{-};(12.5,1) *{+};
(10,3); (10,0) **@{-}, (9,3) **@{-};
(11,1); (11,0) **@{-};
(10,2); (9,2) **@{-};
(13,0); (13,2) **@{-}, (16,0) **@{-};
(16,1); (13,1) **@{-}, (16,0) **@{-};
(15,2); (15,0) **@{-}, (13,2) **@{-};
(14,2); (14,0) **@{-};
(17,0); (17,2) **@{-}, (21,0) **@{-};
(21,1); (17,1) **@{-}, (21,0) **@{-};
(18,2); (17,2) **@{-}, (18,0) **@{-};
(19,0); (19,1) **@{-};
(20,0); (20,1) **@{-};
\end{xy}
\vskip10pt
\noindent
that is, in the algebra $\Lambda_3/\langle S_{(1^3)}-1\rangle$, we have that:
\[ S_{(2,1)}h_2=S_{(2,2,1)}+S_{(3,1,1)}+S_{(3,2)}+S_{(4,1)}. \]
Now, the weight associated to the partition $(\mu)=(2,1)$ is $\mu =\lambda_1+\lambda_2$ and the weight associated to
$h_2=S_{(2,0,0)}$ is $2\lambda_1$. Then the above equation translates into the tensor product in $T(A_2)$
\[ V^{\lambda_1+\lambda_2}\otimes_3 V^{2\lambda_1}= V^{\lambda_2}\oplus V^{2\lambda_1}\oplus V^{\lambda_1+2\lambda_2}
\oplus V^{3\lambda_1+\lambda_2}. \]
Under the correspondence \eqref{E:ow} we have that the orbits corresponding to the weights $\lambda_1+\lambda_2$, 
$2\lambda_1$, $\lambda_2$, $\lambda_1+2\lambda_2$ and $3\lambda_1+\lambda_2$ are respectively 
[(2,1,0,\dots)], [(1,1,0,\dots)], [(2,0,\dots)], [(2,2,1,0\dots)] and [(2,1,1,1,0\dots)]. 
Looking at Example~\ref{E:tensorcomp} we see that this tensor product
agrees with the orbit calculation from that example. 
\end{example}

\section{Weight space decomposition of modules }

In this section we are going to describe how we can use Young tableaux to find the weight decomposition of an irreducible
module for $sl_{N}$. 

We start by noticing that the map \eqref{E:pw} between integral dominant weights of $A_{N-1}$ and partitions of length at 
most $N$ can be extended to a map from the set of sequences of non-negative integers of length $N$ to the weight lattice 
$P$ of $A_{N-1}$ 
\begin{equation} \label{E:sw}
 (\mu)=(\mu_1,\dots,\mu_N)\longmapsto  \mu=\sum_{i=1}^{N-1}(\mu_i-\mu_{i+1})\lambda_i. 
\end{equation}
The map is not one-to-one but it has a right inverse given by
\[ \mu=\sum_{i=1}^{N-1}a_i\lambda_i \longmapsto (\mu)=\left(c+\sum_{i=1}^{N-1}a_i,c+\sum_{i=2}^{N-1}a_i\dots,c+a_{N-1},
c\right), \]
where $c=\min\left\{ \sum_{i=1}^{N-1}a_i,\sum_{i=2}^{N-1}a_i,\dots,a_{N-1},0\right\}$.

Let $V=V^\lambda$ be an irreducible finite dimensional module for $sl_N$ where 
$\lambda=\sum_{i=1}^{N-1}a_i\lambda_i\in P^+$. Let 
\[ V^\lambda=\bigoplus_{\beta\in P}\text{Mult}_\lambda(\beta) V_\beta^\lambda\]
 be the weight space decomposition of $V^\lambda$ where 
\[ V_\beta^\lambda=\{ x \mid h\cdot x=\beta(h)x, \text{ for all } h\in\mathfrak{h}\}.\]
Let $(\lambda)=\left(\sum_{i=1}^{N-1}a_i,\dots,a_{N-1},0\right)$ be a partition associated to $\lambda$. The set of Young 
tableaux of shape $(\lambda)$ filled with numbers from 1 to $N$ is in 1-1 correspondence with the set of weight spaces of 
$V^\lambda$ including multiplicities. The correspondence is as follows:

The number of Young tableaux of shape $(\lambda)$ and content $(\mu)=(\mu_1,\dots,\mu_N)$ is equal to 
Mult$_\lambda(\beta)$ where $\beta=\sum_{i=1}^{N-1}(\mu_i-\mu_{i+1})\lambda_i$.

\begin{example} \label{E:tweightspace}
Consider the module $V^\lambda$ for $A_2$, where $\lambda=\lambda_1+\lambda_2$, the well known adjoint representation of 
$sl_3$. The partition associated to $\lambda$ is $(\lambda)=(2,1,0)$. Below we show the set of Young tableaux of shape 
$(\lambda)$ filled with numbers from 1 to 3, and underneath the each tableau we have the content and the weight associated 
to the content.
\vskip10pt
\begin{xy}
<0.7cm,0cm>:
(-5.38,0); (0,0) **{ };
(0,0); (0,2) **@{-}, (2,0) **@{-};
(1,2); (0,2) **@{-}, (1,0) **@{-};
(2,1); (0,1) **@{-}, (2,0) **@{-};
(.5,.5) *{1};
(1.5,.5) *{1};
(.5,1.5) *{2};
(1,-.5) *{(2,1,0)};
(1,-1.2) *{\lambda_1+\lambda_2};
(3,0); (3,2) **@{-}, (5,0) **@{-};
(4,2); (3,2) **@{-}, (4,0) **@{-};
(5,1); (3,1) **@{-}, (5,0) **@{-};
(3.5,.5) *{1};
(4.5,.5) *{1};
(3.5,1.5) *{3};
(4,-.5) *{(2,0,1)};
(4,-1.2) *{2\lambda_1-\lambda_2};
(6,0); (6,2) **@{-}, (8,0) **@{-};
(7,2); (6,2) **@{-}, (7,0) **@{-};
(8,1); (6,1) **@{-}, (8,0) **@{-};
(6.5,.5) *{1};
(7.5,.5) *{2};
(6.5,1.5) *{2};
(7,-.5) *{(1,2,0)};
(7,-1.2) *{-\lambda_1+2\lambda_2};
(9,0); (9,2) **@{-}, (11,0) **@{-};
(10,2); (9,2) **@{-}, (10,0) **@{-};
(11,1); (9,1) **@{-}, (11,0) **@{-};
(9.5,.5) *{1};
(10.5,.5) *{2};
(9.5,1.5) *{3};
(10,-.5) *{(1,1,1)};
(10,-1.2) *{0};
\end{xy}
\vskip10pt
\begin{xy}
<0.7cm,0cm>:
(6.61,0); (12,0) **{ };
(12,0); (12,2) **@{-}, (14,0) **@{-};
(13,2); (12,2) **@{-}, (13,0) **@{-};
(14,1); (12,1) **@{-}, (14,0) **@{-};
(12.5,.5) *{1};
(13.5,.5) *{3};
(12.5,1.5) *{2};
(13,-.5) *{(1,1,1)};
(13,-1.2) *{0};
(15,0); (15,2) **@{-}, (17,0) **@{-};
(16,2); (15,2) **@{-}, (16,0) **@{-};
(17,1); (15,1) **@{-}, (17,0) **@{-};
(15.5,.5) *{1};
(16.5,.5) *{3};
(15.5,1.5) *{3};
(16,-.5) *{(1,0,2)};
(16,-1.2) *{\lambda_1-2\lambda_2};
(18,0); (18,2) **@{-}, (20,0) **@{-};
(19,2); (18,2) **@{-}, (19,0) **@{-};
(20,1); (18,1) **@{-}, (20,0) **@{-};
(18.5,.5) *{2};
(19.5,.5) *{2};
(18.5,1.5) *{3};
(19,-.5) *{(0,2,1)};
(19,-1.2) *{-2\lambda_1+\lambda_2};
(21,0); (21,2) **@{-}, (23,0) **@{-};
(22,2); (21,2) **@{-}, (22,0) **@{-};
(23,1); (21,1) **@{-}, (23,0) **@{-};
(21.5,.5) *{2};
(22.5,.5) *{3};
(21.5,1.5) *{3};
(22,-.5) *{(0,1,2)};
(22,-1.2) *{-\lambda_1-\lambda_2};
\end{xy}
\vskip5pt
Therefore we have that (leaving off the superscript $\lambda_1+\lambda_2$ from the weight spaces on the right hand side)
\[V^{\lambda_1+\lambda_2}=V_{\lambda_1+\lambda_2}\oplus V_{2\lambda_1-\lambda_2}\oplus 
V_{-\lambda_1+2\lambda_2}\oplus 2V_0\oplus V_{\lambda_1-2\lambda_2}\oplus V_{-2\lambda_1+\lambda_2}\oplus 
V_{-\lambda_1-\lambda_2}.\]
\end{example}

Now, we discuss how we can use this to implement the Racah-Speiser algorithm.

The Weyl group $W\cong S_N$ of $sl_N$ acts on sequences of length $N$ as permutations, with the simple reflections acting as
transpositions 
\begin{equation} \label{E:weyl}
 r_i(\mu_1,\dots,\mu_N)=(\mu_1,\dots,\mu_{i+1},\mu_i,\dots,\mu_N),\qquad i=1,...,N-1. 
\end{equation}
This action allows one to get a version of the Racah-Speiser algorithm from the Young tableau point of view as follows.

Let $V^\lambda$ and $V^\mu$ be irreducible finite dimensional modules of type $A_{N-1}$, $(\lambda)$ and $(\mu)$ the 
partitions associated to $\lambda$ and $\mu$ respectively. The decomposition of the tensor product $V^\lambda\otimes V^\mu$
into irreducible modules can be computed by doing the following.

Step 1. List all contents of the Young tableaux of shape $(\lambda)$ whose fillings are with numbers from 1 to $N$.

Step 2. Add $(\mu)+(\rho)$ to all contents from step 1, where $(\rho)=(N-1,N-2,\dots,1,0)$.

Step 3. Apply the action of the Weyl group \eqref{E:weyl} to all sequences that are not in standard form from step 2 to write
them in standard form, with positive multiplicity if the number of transpositions is even and negative multiplicity if the 
number of transpositions is odd and we drop the sequences that are fixed by any transposition.

Step 4. Subtract $(\rho)$ from all partitions left in step 3 and use the correspondence \eqref{E:sw} to get the weights 
$\lambda$ associated to the partitions from step 3. The direct sum of the irreducible modules indexed by these weights equals
the tensor product decomposition of $V^\lambda\otimes V^\mu$.

\begin{example} 
Let $V^{\lambda_1+\lambda_2}$ and $V^{2\lambda_1}$ be irreducible modules for $A_2$. The weights of $V^{\lambda_1+\lambda_2}$
were given in Example~\ref{E:tweightspace} as well as all contents of tableaux of shape $(\lambda_1+\lambda_2)=(2,1)$. If we 
add $(2\lambda_1)+(\rho)=(2,0,0)+(2,1,0)=(4,1,0)$ to all these weights we get the sequences
\[ (6,2,0),\quad (6,1,1),\quad(5,3,0),\quad 2(5,2,1),\quad(5,1,2),\quad(4,3,1),\quad(4,2,2), \]
where the coefficient 2 in front of $(5,2,1)$ is the number of tableaux of shape $(\lambda_1+\lambda_2)=(2,1)$ and that 
content. 

The sequences $(6,1,1)$ and $(4,2,2)$ are fixed by the action of $r_2$, therefore they do not count for the tensor product.
The weight $(5,1,2)$ is not in standard form, but if we apply $r_2$ to it, we get $r_2\cdot(5,1,2)=(5,2,1)$, so it reduces 
the multiplicity of $(5,2,1)$ to 1. So we are left with the sequences
\[ (6,2,0),\qquad(5,3,0),\qquad (5,2,1),\qquad(4,3,1). \]
By subtracting $(\rho)=(2,1,0)$ from these weights we get
\[ (4,1,0),\quad(3,2,0),\quad (3,1,1),\quad(2,2,1) \]
and using the map \eqref{E:sw} we get that the weights associated to these sequences are 
\[ 3\lambda_1+\lambda_2,\qquad \lambda_1+2\lambda_2,\qquad 2\lambda_1,\qquad\lambda_2. \]
Therefore the tensor product decomposition of $V^{\lambda_1+\lambda_2}\otimes V^{2\lambda_1}$ is given by
\[ V^{\lambda_1+\lambda_2}\otimes V^{2\lambda_1}=V^{3\lambda_1+\lambda_2} \oplus V^{\lambda_1+2\lambda_2} \oplus 
V^{2\lambda_1} \oplus V^{\lambda_2}. \]
\end{example}

The advantage of this method over the geometrical one, presented in section \ref{S:racaspeicer}, is that we can compute 
tensor product decompositions by hand for ranks higher than 2 on a piece of paper. 

Now we present a version of the Kac-Walton algorithm from the tableau point of view. For a given $k>0$, the affine reflection
$r_0$ acts on the finite dimensional weight lattice $P$, and therefore on sequences of non-negative numbers of length at most
$N$. This action is defined by 
\begin{equation} \label{E:weyl1}
 r_0\cdot(a_0,\dots,a_{N-1})=\begin{cases}
(k+N+a_{N-1},a_1,\dots,a_{N-2},a_0-k-N), & \\ \qquad\qquad\qquad\qquad\qquad\qquad\text{if}\qquad a_0-a_{N-1} 
\geq k+N, \\
(k+N+a_{N-1}+c,a_1+c,\dots,a_{N-2}+c,a_0-k-N+c), &\\ \qquad\qquad\qquad \qquad\qquad\qquad\qquad\text{if}\qquad 
a_0-a_{N-1} <k+N, 
\end{cases}
\end{equation}
where $c=k+N-a_0$ if $a_0-k-N<0$ and $c=0$ otherwise. 

With the level $k$ action of the affine Weyl group $\widehat{W}$ defined this way, we get the following algorithm for fusion 
products. Note that the fundamental region for the action of $\widehat{W}$ defined by \eqref{E:weyl} and \eqref{E:weyl1} on 
length $N$ sequences of non-negative integers is the set of partitions $(a_0,\dots,a_{N-1})$ of length at most $N$  such that
$a_0-a_{N-1}\leq N+k$. We denote this set by $H_{N+k}$. That is
\[ H_{N+k}=\{(a_0,\dots,a_{N-1})\in\mathbb{Z}_{\geq 0}^N\mid a_0\geq\dots\geq a_{N-1}\text{ and }a_0-a_{N-1}\leq N+k \}. \]

Let $V^\lambda$ and $V^\mu$ be irreducible highest weight modules of type $A_{N-1}^{(1)}$ on level $k$, $(\lambda)$ and 
$(\mu)$ the partitions associated to $\lambda$ and $\mu$ respectively. The level $k$ fusion product 
$V^\lambda\otimes_k V^\mu$ into irreducible modules can be computed by doing the following.

Step 1. List all contents of the Young tableaux of shape $(\lambda)$ whose fillings are with numbers from 1 to $N$.

Step 2. Add $(\mu)+(\rho)$ to all contents from step 1, where $(\rho)=(N-1,N-2,\dots,1,0)$.

Step 3. Apply the action of the affine Weyl group defined by \eqref{E:weyl} and \eqref{E:weyl1} to all sequences from step 2 
to get them into the fundamental region $H_{N+k}$, with positive multiplicity if the number of reflections is even and 
negative multiplicity if the number of reflections is odd and drop the sequences that are fixed by any reflection.

Step 4. Subtract $(\rho)$ from all partitions left in step 3 and use the correspondence \eqref{E:sw} to get the weights 
$\lambda$ associated to the partitions from step 3. The direct sum of the irreducible modules indexed by these weights 
equals the level $k$ fusion product $V^\lambda\otimes_k V^\mu$.

\begin{example} 
Let $V^{\lambda_1+\lambda_2}$ and $V^{2\lambda_1}$ be irreducible modules for $A_2^{(1)}$ and $k=2$. The weights of 
$V^{\lambda_1+\lambda_2}$ were given in Example~\ref{E:tweightspace} as well as all contents of tableaux of shape 
$(\lambda_1+\lambda_2)=(2,1)$. If we add $(2\lambda_1)+(\rho)=(2,0,0)+(2,1,0)=(4,1,0)$ to all these weights we get the 
sequences
\[ (6,2,0),\quad (6,1,1),\quad(5,3,0),\quad 2(5,2,1),\quad(5,1,2),\quad(4,3,1),\quad(4,2,2), \]
where the coefficient 2 in front of $(5,2,1)$ is the number of tableaux of shape $(\lambda_1+\lambda_2)=(2,1)$ and that 
content. 

The sequences $(6,1,1)$ and $(4,2,2)$ are fixed by the action of $r_2$, and $(5,3,0)$ is fixed by $r_0$, therefore they do 
not count for the fusion product. The weights $(5,1,2)$ and  $(6,2,0)$ are outside the fundamental region $H_5$. Applying 
$r_2$ to $(5,1,2)$, we get $r_2\cdot(5,1,2)=(5,2,1)$, and applying $r_0$ to $(6,2,0)$ we get $r_0\cdot(6,2,0)=(5,2,1)$. 
Therefore the multiplicity of $(5,2,1)$ reduces to 0. So we are left with the sequence
\[ (4,3,1). \]
By subtracting $(\rho)=(2,1,0)$ from this sequence we get
\[ (2,2,1) \]
and using the map \eqref{E:sw} we get that the weight associated to this sequence is $\lambda_2$. Therefore we get the fusion
product 
\[ V^{\lambda_1+\lambda_2}\otimes_2 V^{2\lambda_1}=V^{\lambda_2}. \]
\end{example}

In future work we plan to investigate the dependence of the fusion coefficients on level $k$ by using these techniques.



\end{document}